\documentclass[11pt,reqno]{amsart}

\theoremstyle{plain}
\newtheorem{theorem} {Theorem}[section]
\newtheorem{lemma}[theorem] {Lemma}
\newtheorem{proposition}[theorem] {Proposition}
\newtheorem{corollary}[theorem] {Corollary}

\theoremstyle{definition}
\newtheorem{definition}[theorem] {Definition}

\newtheorem{example} [theorem]{Example}

\theoremstyle{remark}
\newtheorem{remark}[theorem] {Remark}

\numberwithin{equation}{section}

\usepackage{amsmath}
\usepackage{amsfonts}
\usepackage{amssymb}
\usepackage{amscd}
\DeclareMathOperator*{\esssup}{ess\,sup}

\newcommand{\R}{{\mathbb R}}
\newcommand{\Z}{{\mathbb Z}}
\newcommand{\N}{{\mathbb N}}

\newcommand{\PP}{{\mathcal P}}

\newcommand{\C}{{\mathcal C}}
\newcommand{\CC}{{\mathbb C}}

\newcommand{\TT}{{\mathcal T}}

\newcommand{\al}{{\alpha}}
\newcommand{\la}{{\lambda}}
\newcommand{\sa}{{\sigma}}

\newcommand{\iy}{{\infty}}
\newcommand{\vphi}{{\varphi}}
\newcommand{\vep}{{\varepsilon}}
\newcommand{\g}{{\gamma}}
\newcommand{\de}{{\delta}}

\newcommand{\Ome}{{\Omega}}
\newcommand{\om}{{\omega}}

\newcommand{\be}{{\beta}}

\newcommand{\BB}{{\mathfrak B}}
\newcommand{\bna}{\begin{eqnarray}}
\newcommand{\ena}{\end{eqnarray}}
\newcommand{\ba}{\begin{eqnarray*}}
\newcommand{\ea}{\end{eqnarray*}}
\newcommand{\beq}{\begin{equation}}
\newcommand{\eeq}{\end{equation}}

\begin{document}

\title[Sharp Constants of Approximation Theory]
{Sharp Constants of Approximation Theory. IV. Asymptotic Relations in General Settings}
\author{Michael I. Ganzburg}
 \address{Department of Mathematics\\ Hampton University\\ Hampton,
 VA 23668\\USA}
 \email{michael.ganzburg@gmail.com}
 \begin{abstract}
 In this paper we first introduce the
 unified definition of the sharp constant that includes constants in
three major problems of approximation theory, such as,  inequalities for approximating elements,
 approximation of individual elements, and approximation on
classes of elements. Second, we find sufficient conditions that imply limit relations between
 various sharp constants of approximation theory in general settings. Third,
 a number of examples from various areas of approximation theory illustrates the general approach.
\end{abstract}
 \maketitle

\section{Introduction}\label{S1}
\setcounter{equation}{0}
\noindent
We continue the study of the sharp constants
of approximation theory
that began in \cite{G1992,G2000,GT2017,G2017,G2018,G2019a,G2019b,G2021b,G2021c,G2021}.
In this paper we discuss limit relations between  sharp constants of approximation theory
in general settings.

Sharp constants of approximation theory have attracted much attention of approximation analysts
for more than a century (see monographs
\cite{ T1963, A1965, K1991, DL1993,  MMR1994, BE1995, KLB1996, BKKP2003} and references therein).
Our unified definition of the sharp constant (see \eqref{E2.5}) includes constants in
three major problems of approximation theory, such as,  inequalities for approximating elements
(Problem A), approximation of individual elements (Problem B), and approximation on
classes of elements (Problem C).

It was Bernstein who in 1913 and 1938 initiated the study of asymptotic relations in Problem B
for the errors of best polynomial approximation of $f_\la(\tau):=\vert \tau\vert^\la$
 by proving the following celebrated result
(see \cite{B1913} for $\la=1$ and \cite{B1938} for $\la>0$): for any $\la>0$, there is a constant
$\mu_\la\in(0,\iy)$ such that
\beq\label{E1.1}
\lim_{n\to\iy}n^\la \min_{P_n\in\PP_{n}}
\max_{\tau\in[-1,1]}\left\vert
f_\la(\tau)-P_n(\tau)\right\vert=\mu_\la,
\eeq
where $\PP_n$ is the class of all univariate algebraic
polynomials of degree at most $n$.
Nikolskii \cite{N1947} proved that an $L_1[-1,1]$-version of \eqref{E1.1}
is valid as well.
In addition to \eqref{E1.1}, Bernstein \cite{B1938} also showed that
$\mu_\la\sa^{-\la}$ coincides with the error of best approximation of $f_\la$ by
entire functions of exponential type $\sa,\,\sa>0$, on the real axis, that is,
\beq\label{E1.2}
\mu_\la=\sa^\la
\inf_{g\in B_\sa}
\sup_{\tau\in(-\iy,\iy)}\left\vert
f_\la(\tau)-g(\tau)\right\vert.
\eeq
Later Bernstein \cite{B1946, B1947(lim)}
(see also \cite[Sect. 2.6.22]{T1963} and \cite[Appendix, Sect. 84]{A1965})
extended the combination of relations
\eqref{E1.1} and \eqref{E1.2} to any continuous function $f$ of polynomial growth
in the following form:
\beq\label{E1.3}
\lim_{n\to\iy}\min_{P_n\in\PP_{n}}
\max_{\tau\in[-n/\sa,n/\sa]}\left\vert
f(\tau)-P_n(\tau)\right\vert
=\inf_{g\in B_\sa}
\sup_{\tau\in(-\iy,\iy)}\left\vert
f(\tau)-g(\tau)\right\vert.
\eeq
This relation holds true for almost all $\sa>0$,
while for certain functions (e.g., for $f_\la$) it is valid for all $\sa>0$.
In particular, \eqref{E1.3} implies \eqref{E1.1} and \eqref{E1.2}.
An $L_p$-version of \eqref{E1.3} for
$p\in [1,\iy)$ and all $\sa>0$ was obtained by Raitsin \cite{R1968}.
Certain periodic analogues of \eqref{E1.3} were given in \cite{B1946(trig)}.
The corresponding multivariate extensions of these results were established by the author
\cite{G1982, G1991,  G1992, G2000} (see also Examples \ref{Ex7.8} and
 \ref{Ex7.9} in Section \ref{S7}).
 Weighted $L_p$-versions of \eqref{E1.3}
 for $p\in (0,\iy]$ were obtained in
 \cite{G2008} (see also Example \ref{Ex7.10}).
 $L_p$-versions of \eqref{E1.1}
 and \eqref{E1.2} for complex $\la$ with $\mathrm{Re}\,\la>
 \max\{-1,-1/p\},\,p\in(0,\iy]$, or $\mathrm{Re}\,\la=0,\,p=\iy$,
  were recently proved in \cite{G2021a}.

 In addition, Bernstein \cite{B1946(Nik)} discussed limit relations in Problem C
 between the least upper bounds of the errors of best harmonic and polynomial approximation
 on H\"{o}lder classes $H^\la,\,\la\in(0,1]$. In particular, he proved the following relation:
 \beq\label{E1.4}
\lim_{n\to\iy}n^\la\sup_{f\in H^\la[-1,1]}\min_{P_n\in\PP_{n}}
\max_{\tau\in[-1,1]}\left\vert
f(\tau)-P_n(\tau)\right\vert
=\sup_{f\in H^\la(-\iy,\iy)}\min_{g\in B_1}
\sup_{\tau\in(-\iy,\iy)}\left\vert
f(\tau)-g(\tau)\right\vert.
\eeq
More general versions of \eqref{E1.4} were discussed by Bernstein \cite{B1947(const)} and by
the author \cite{G1992, G2003} (see also Examples \ref{Ex7.15} and \ref{Ex7.16}).

A general approach to limit (asymptotic) relations between sharp constants in normed spaces
in Problems B and C was developed by the author \cite{G1992}. This approach was later
extended to $\be$-normed spaces \cite{G2000}.

The first asymptotic relations between sharp constants in Nikolskii-type inequalities
of Problem A were established by Taikov \cite{T1965, T1993}, Garsia et al \cite{GRR1969},
 and Gorbachev \cite{G2005}.
 Levin and Lubinsky  \cite{LL2015a, LL2015b} obtained more general results.
 The latest asymptotics for the sharp constants in
 Nikolskii-Bernstein  type inequalities
 were given by the author and Tikhonov \cite{GT2017}; in particular,
 the following relation was proved in \cite{GT2017}:
 \beq\label{E1.5}
\lim_{n\to\iy}n^{-s-1/p}\sup_{Q\in \TT_n\setminus\{0\}}
\frac{\left\|Q^{(s)}\right\|_{L_\iy([0,2\pi))}}
{\left\|Q\right\|_{L_p([0,2\pi))}}
= \sa^{-s-1/p}\sup_{g\in \left(B_\sa\cap L_p(\R)\right)\setminus\{0\}}
\frac{\left\|g^{(s)}\right\|_{L_\iy((-\iy,\iy))}}
{\left\|g\right\|_{L_p((-\iy,\iy))}},
\eeq
where $\TT_n$ is the class of all univariate trigonometric
polynomials of degree at most $n$
and $p\in(0,\iy]$. A multivariate generalization of \eqref{E1.5}
was obtained in \cite{G2018}.
Univariate and multivariate  versions of \eqref{E1.5} for algebraic polynomials
were given in \cite{G2017, DGT2018, G2019a,  G2019b, DGT2019} (see also
Examples \ref{Ex7.2}, \ref{Ex7.3}, and \ref{Ex7.4}).

This paper was inspired by several publications
\cite{GT2017, G2017,
   G2018, G2019a, G2019b} about limit relations between sharp constants
   in univariate and multivariate Markov-Bernstein type inequalities of different metrics
   for trigonometric and algebraic polynomials and entire functions of exponential type.
   Analyzing these publications, we noticed that the proofs of the
   corresponding limit relations are based
   on verification of certain conditions that often
   coincide with the conditions introduced in
   \cite{G1992, G2000}. As a result of this observation, in this paper we develop
   a more general approach,
   compared with \cite{G1992, G2000}, to limit relations for sharp constants in vector
   spaces. In addition to Problems B and C, this approach also includes
   inequalities for approximating elements. However, in certain details the approach
   differs from the one in \cite{G1992, G2000} for Problems B and C.

   In Section 2 the definition of the sharp constant $\C$ is introduced,
   and special cases of $\C$ in approximation theory and analysis are discussed.
   In Section 3 we introduce conditions on sets, vector spaces, and operators and study
   their properties.
   In Section 4 we prove asymptotic relations between sharp constants
   in vector spaces that are based on conditions from Section 3.
   The proofs are simple but we believe that this general
   condition-based approach is useful in approximation theory.
   Special cases of asymptotic relations in Problems A, B, and C
   are discussed in Section 5.

   Function spaces, function classes,  special operators, and their properties
    are discussed in Section 6. In Section 7 for each of Problems A, B, and C,
    we discuss
    limit relations between sharp constants in function spaces,
     including three examples
    of Markov-Bernstein type inequalities of different metrics,
    three examples of limit theorems for polynomial and harmonic approximation,
    and two examples of limit relations between constants of approximation theory.
    In the capacity of approximating elements in the examples, we use univariate and
    multivariate trigonometric and algebraic polynomials
    and entire functions of exponential type. These examples cover most of known results
     in this area and Example \ref{Ex7.15} is new.

 \section{Definitions and Properties of Sharp Constants in Vector Spaces}\label{S2}
 \setcounter{equation}{0}
\noindent
Here, we first describe major spaces that are discussed in this paper. Next, we give
the unified definition of a
sharp constant in general settings and study some of its properties. Finally, we discuss
certain special cases of sharp constants of approximation theory and analysis.
\subsection{Spaces}\label{S2.1}
Let $F$ be a vector space over the field of complex numbers $\CC$
 or real numbers $\R$.  Let $F$
be equipped with a functional
$\|\cdot\|_F:F\to[0,\iy)$, which is not necessarily a norm though we use the same notation.
Occasionally, we will refer to $\|\cdot\|_F$ as a "norm." In particular, $\|\cdot\|_F$ could be
a norm, a seminorm, or more generally, a $\be$-norm or a $\be$-seminorm.
\begin{definition}\label{D2.1}
\cite[p. 1102]{K2003}
We say that $F$ is a \emph{$\beta$-normed space}, $\beta\in(0,1]$, if there exists
a   $\be$-\emph{norm} $\|\cdot\|_{F}$ on $F$ with the following properties:
\begin{itemize}
\item[(a)] $\|f\|_{F}\ge 0$ for all $f\in F$ and $\|f\|_{F}= 0$ if and only if $f=0$.
\item[(b)] $\|\alpha f\|_{F}=|\alpha|\,\|f\|_{F}$ for all $\alpha\in\CC$ (or $\alpha\in\R$).
\item[(c)] $\|f+g\|_{F}\le \left(\| f\|_{F}^\beta +\| g\|_{F}^\beta
\right)^{1/\beta}$ for all $f\in F,\,g\in F$.
\end{itemize}
\end{definition}
\noindent
Occasionally, we will refer to (b) and (c) as the homogeneity property and
 the $\be$-triangle inequality, respectively.
In case of $\be=1,\, \|\cdot\|_F$ is a norm and $F$ is a normed space.
Note that a   $\be$-norm is closely related to a $\g$-quasinorm $|||\cdot|||_F$ satisfying
conditions (a), (b), and (c*) $|||f+g|||_F\le \g (|||f|||_F+|||g|||_F),
\,\g\ge 1;\,f\in F,\,g\in F.$
Obviously, a   $\be$-norm is a $2^{1/\be-1}$-quasinorm. Conversely, the Aoki-Rolewicz theorem
\cite{A1942, R1957} states that if for some $\beta\in(0,1],\, |||\cdot|||_F$ is
a $2^{1/\be-1}$-quasinorm, then there exists a   $\be$-norm $\|\cdot\|_F$, which
 is equivalent to $|||\cdot|||_F$.
 The best known examples of $\beta$-normed spaces
 for $\be\in(0,1)$ are the sequence space $\ell_\be$,
 the function space $L_\be$, and the Hardy space $H_\be$.

If we replace property (a) with
\begin{itemize}
\item[(a*)] $\|f\|_{F}\ge 0$ for all $f\in F$ ,
\end{itemize}
then $\|\cdot\|_{F}$ and $F$ are called
a $\be$-\emph{seminorm} and a \emph{$\beta$-seminormed space},
respectively.
For $\beta=1,\,\|\cdot\|_{F}$
is a seminorm and $F$ is a seminormed space.

However, in this paper $\|\cdot\|_{F}$ is not necessarily a
$\beta$-norm or a $\beta$-seminorm.
For example, the trivial functional
\beq\label{E2.1}
\|f\|_{F,\mathrm{triv}}:=\left\{\begin{array}{ll}
0, &f=0,\\
1,&f\ne 0,
\end{array}\right.
\qquad f\in F,
\eeq
satisfies property (a) and also property (c)
 for any $\beta\in(0,1]$
but does not satisfy property (b).
Meantime, $\|\cdot\|_{F,\mathrm{triv}}$
satisfies the symmetry property
\begin{itemize}
\item[(b*)] $\|-f\|_{F}=\|f\|_{F}$ for all $f\in F$ .
\end{itemize}

Note that a vector space $F$ can be equipped with two (or more)
 different "norms;" we often denote them by $\|\cdot\|_{F}$
 and $\|\cdot\|_{F}^*$.
So the equality $F=G$ for vector spaces $F$ and $G$  means that they
 are equal as sets but can be equipped with different "norms,"
 while the equality $\|\cdot\|_F=\|\cdot\|_G$ also means that $F=G$.

In this paper we also discuss sequences of "norms"
$\left\{\|\cdot\|_{F_n}\right\}_{n=1}^\iy$
that possess either the continuity (C) or generalized continuity (GC) properties,
which can replace
the $\be$-triangle inequality in certain problems.

\begin{definition}\label{D2.2}
We say that $\left\{\|\cdot\|_{F_n}\right\}_{n=1}^\iy$ possesses the C-property if for any
$h_n\in F_n,\,g_n\in F_n,\,n\in\N$, such that
$\lim_{n\to\iy}\|g_n\|_{F_n}=0$,
the following relations hold true:
\bna
&&\liminf_{n\to\iy}\|h_n+g_n\|_{F_n}=\liminf_{n\to\iy}\|h_n\|_{F_n},\label{E2.2}\\
&&\limsup_{n\to\iy}\|h_n+g_n\|_{F_n}=\limsup_{n\to\iy}\|h_n\|_{F_n}.\label{E2.3}
\ena
\end{definition}
\noindent
\begin{proposition}\label{P2.3}
(i) If "norms" $\|\cdot\|_{F_n},\,n\in\N$, satisfy the symmetry property (b*)
and the $\be$-triangle inequality (c), where $\be\in(0,1]$ is independent of $n$, then
$\left\{\|\cdot\|_{F_n}\right\}_{n=1}^\iy$ possesses the C-property.\\
(ii) If $\|\cdot\|_{F_n}$ is a $\beta$-norm or a $\beta$-seminorm,
$n\in\N,$
where $\be\in(0,1]$ is independent of $n$, then
$\left\{\|\cdot\|_{F_n}\right\}_{n=1}^\iy$ possesses the C-property.\\
(iii) If $\|\cdot\|_{F_n}=\|\cdot\|_{F_n,\mathrm{triv}},\,n\in\N,$
then
$\left\{\|\cdot\|_{F_n}\right\}_{n=1}^\iy$ possesses the C-property.
\end{proposition}
\proof
Statement (i) follows immediately from the $\be$-triangle inequalities
\ba
\| h_n\|_{F_n}^\beta -\| g_n\|_{F_n}^\beta
\le \|h_n+g_n\|_{F_n}^\be\le \| h_n\|_{F_n}^\beta +\| g_n\|_{F_n}^\beta,
\qquad n\in\N,
\ea
while statements (ii) and (iii) are easy consequences of (i).\hfill $\Box$

\begin{definition}\label{D2.4}
Let $F_n$ be equipped with two "norms" $\|\cdot\|_{F_n}$ and
$\|\cdot\|_{F_n}^*,\,n\in\N$.
We say that $\left\{\|\cdot\|_{F_n}\right\}_{n=1}^\iy$ possesses the GC-property
with respect to the sequence 
$\left\{\|\cdot\|_{F_n}^*\right\}_{n=1}^\iy$ if for any
$h_n\in F_n,\,g_n\in F_n,\,n\in\N$, such that
$\lim_{n\to\iy}\|g_n\|_{F_n}^*=0$,
relations \eqref{E2.2} and \eqref{E2.3} hold true.
\end{definition}
\noindent
In case of $\|\cdot\|_{F_n} =\|\cdot\|_{F_n}^*,\,n\in\N $,
the GC-property is reduced to the C-property. In most of our
applications we use sequences of
"norms" with the C-property. However, in the following
typical example (related to approximation theory)
we present a sequence of "norms" that possesses the GC-property but does
not possess the C-property.
\begin{example}\label{Ex2.5}
Let $F_n\ne\{0\}$ and let  $\|\cdot\|_{F_n}^*,\,n\in\N,$ satisfy property (b)
of Definition \ref{D2.1}.
 Next, let
$\left\{\|\cdot\|_{F_n}^*\right\}_{n=1}^\iy$ possess the C-property. Let us set
\ba
\|f\|_{F_n}:=\|\vphi_n-f\|_{F_n}^*,\qquad f\in F_n,
\ea
 where $\vphi_n$ is a fixed element of $F_n,\,n\in\N$.
 Then $\left\{\|\cdot\|_{F_n}\right\}_{n=1}^\iy$ possesses the GC-property
with respect to $\left\{\|\cdot\|_{F_n}^*\right\}_{n=1}^\iy$.
Indeed, for any
$h_n\in F_n,\,g_n\in F_n,\,n\in\N$, such that
\ba
\lim_{n\to\iy}\|g_n\|_{F_n}^*
=\lim_{n\to\iy}\|-g_n\|_{F_n}^*=0,
\ea
we have
\bna\label{E2.4}
\liminf_{n\to\iy}\|h_n+g_n\|_{F_n}
&=&\liminf_{n\to\iy}\|(\vphi_n-h_n)+(-g_n)\|_{F_n}^*\nonumber\\
&=&\liminf_{n\to\iy}\|\vphi_n-h_n\|_{F_n}^*\nonumber\\
&=&\liminf_{n\to\iy}\|h_n\|_{F_n}.
\ena
Similar relations hold true if we replace $\liminf$
in \eqref{E2.4} with $\limsup$.
Hence \eqref{E2.2} and \eqref{E2.3} are valid.

On the other hand, since $F_n\ne\{0\},\,n\in\N$,
then by property (b) of Definition \ref{D2.1}, it is always possible to choose  $\vphi_n\in F_n,\,n\in\N$,
such that
$\lim_{n\to\iy}\|\vphi_n\|_{F_n}^*=a$, where $a\in(0,\iy)$ is a fixed number.
Then for $h_n:=(1/3)\vphi_n$ and $g_n:=(1-1/n)\vphi_n,\,n\in\N$, we have
$
\lim_{n\to\iy}\|g_n\|_{F_n}
=\lim_{n\to\iy}(1/n)\|\vphi_n\|_{F_n}^*=0
$
and
\ba
\lim_{n\to\iy}\|h_n+g_n\|_{F_n}
=\lim_{n\to\iy}\|\vphi_n-h_n-g_n\|_{F_n}^*=a/3;\quad
\lim_{n\to\iy}\|h_n\|_{F_n}
=\lim_{n\to\iy}\|\vphi_n-h_n\|_{F_n}^*=2a/3.
\ea
These relations show that in general $\left\{\|\cdot\|_{F_n}\right\}_{n=1}^\iy$
does not possess the C-property.
\end{example}

\subsection{Sharp Constants}\label{S2.2}
For vector spaces $F$ and $H$ equipped with the corresponding "norms" and for
an operator $L:F\to H$, we
define the "norm" of $L$ by the formula
\beq\label{E2.4a}
\|L\|=\|L\|_{F\to H}:=\sup_{h\in F,\,\|h\|_{F}> 0}
\frac{\|L(h)\|_{H}}{\|h\|_{F}}.
\eeq
If $B$ is a subset of $F$, we
define the "norm" $\|L\|$ of an operator $L:B\to H$ on $B$ by the right hand side of
\eqref{E2.4a} with $h\in F$ replaced by $h\in B$.

Let $F^{(j)}\ne\{0\}$ be a vector space over $\CC$ or $\R$ equipped with a "norm"
$\|\cdot\|_{F^{(j)}},\,j=1,\,2,$ and let $B$ be a nonempty subset of
 $F^{(2)}$. Next, let $D:B\to F^{(1)}$ be an operator.
In case of the \emph{imbedding} operator $D=I$, we assume that
$F^{(1)}\cap F^{(2)}\ne \emptyset$ and
$B\subseteq F^{(1)}\cap F^{(2)}$.
Certainly, if $F^{(1)}= F^{(2)}$ and $\|\cdot\|_{F^{(1)}}=\|\cdot\|_{F^{(2)}}$,
then $I$ is the identity operator.

Let us define the following \emph{sharp constant}:
\beq\label{E2.5}
\C=\C\left(D,B,F^{(1)},F^{(2)}\right):=\sup_{f\in B,\,\|f\|_{F^{(2)}}> 0}
\frac{\|D(f)\|_{F^{(1)}}}{\|f\|_{F^{(2)}}},
\eeq
which is the "norm" of $D$
on $B$. Clearly, $\C\in[0,\iy]$.

In certain cases the following formulae for $\C$ are valid:
\bna
&&\C=\sup_{f\in B,\,\|f\|_{F^{(2)}}=1}\|D(f)\|_{F^{(1)}}, \label{E2.6}\\
&&\C=\sup_{f\in B,\,\|f\|_{F^{(2)}}>0,\,\|D(f)\|_{F^{(1)}}=1}
\frac{1}{\|f\|_{F^{(2)}}}.\label{E2.7}
\ena
In particular, it is easy to verify the following proposition.
\begin{proposition}\label{P2.6}
(i) Let $B$ be a subspace of $F^{(2)},\, D$ be a $1$-homogeneous operator
(that is, $D(\al f)=\al D(f),\,\al>0$), and "norms" $\|\cdot\|_{F^{(j)}},\,j =1,\,2,$
satisfy the homogeneity property (b). Then relations \eqref{E2.6} and
\eqref{E2.7} hold true.\\
(ii) If $\|\cdot\|_{F^{(2)}}=\|\cdot\|_{F^{(2)},\mathrm{triv}}$,
then \eqref{E2.6} holds true.\\
(iii) If $\|\cdot\|_{F^{(1)}}=\|\cdot\|_{F^{(1)},\mathrm{triv}}$,
then \eqref{E2.7} holds true.
\end{proposition}
\subsection{Sharp Constants in Approximation Theory and Analysis}\label{S2.3}
It turns out that the problem of finding $\C=\C\left(D,B,F^{(1)},F^{(2)}\right)$
includes most classical problems
of approximation theory that are related to sharp constants. In addition, it includes
certain analysis problems as well. Below we first discuss three
major approximation problems A, B, and C  as special cases of the problem
of finding $\C$.\vspace{.15in}\\
\emph{A. Inequalities for Approximating Elements.}
Let $F^{(j)}$ be a $\beta_j$-normed functional space, $j=1,2;\, B$ be
a nontrivial subspace of $F^{(2)}$; and $D$ be a linear (often differential)
 operator.
Then $\C$ is the sharp constant in
various inequalities for approximating elements
that we call
\emph{Markov-Bernstein type
inequalities of different metrics.} This is one of the most studied topics
of approximation theory. Several special problems are discussed below.\vspace{.15in}\\
\emph{A1.  Markov-Bernstein Type Inequalities for
$\|\cdot\|_{F^{(1)}}=\|\cdot\|_{F^{(2)}}$.}\\
\emph{A1.1. A. A. Markov and  V. A. Markov type inequalities} for the set $B$
of algebraic polynomials of one or several variables. \\
\emph{A1.2. Bernstein type inequalities} for the set $B$
of trigonometric polynomials, splines, or entire functions of exponential type
of one or several variables.
\vspace{.15in}\\
\emph{A2.  Inequalities of Different Metrics for $D=I$ and
$\|\cdot\|_{F^{(1)}}\ne \|\cdot\|_{F^{(2)}}$.}\\
 \emph{A2.1. Nikolskii type inequalities} for the set $B$
of trigonometric or algebraic polynomials, splines, or  entire functions of
exponential type
of one or several variables.
\\
\emph{A2.2. Schur type inequalities} for  a nonweighted space $F^{(1)}$,
a weighted space $F^{(2)}$,
 and the set $B$
of trigonometric or algebraic polynomials or  entire functions of exponential type
of one or several variables.
\\
\emph{A2.3. Remez type inequalities} for a weighted  space $F^{(2)}$
when a weight vanishes outside of a measurable set
 and for the set $B$
of trigonometric, algebraic, or M\"{u}ntz polynomials.

There are numerous publications related to all these inequalities, see,
e.g.,
\cite{T1963, W1974, MMR1994, BE1995, KLB1996, KS1997,  G2001, BKKP2003,
G2005, DP2010,
EGN2015, AD2015, DP2016, G2017, G2018, ABD2018, G2019a, G2019b}
and references therein.
\vspace{.15in}\\
\emph{B. Approximation of Individual Elements.}
Let us define the operator $D:B\to F^{(1)}$ arbitrarily with the only restriction that
 $D(f)\ne 0$ for $f\in B$ and let us choose
$
\|\cdot\|_{F^{(1)}}=\|\cdot\|_{F^{(1)},\mathrm{triv}}.
$
We denote by
$F^{(2)}_*$ the vector space $F^{(2)}$ equipped with a certain $\be_2$-norm
$\|\cdot\|^{*}_{F^{(2)}}$ on $F^{(2)},\,\be_2\in(0,1]$.
For a fixed element $\vphi\in F^{(2)}_*\setminus B$
we choose $\|f\|_{F^{(2)}}=\|\vphi-f\|^*_{F^{(2)}},\,f\in F^{(2)}_*$.
 Then
\beq\label{E2.8*}
\C=
\sup_{f\in B,\,f\ne 0}\frac{\|f\|_{F^{(1)},
\mathrm{triv}}}{\|\vphi-f\|^{*}_{F^{(2)}}}
=\frac{1}{\inf_{f\in B,\,f\ne 0}\|\vphi-f\|^{*}_{F^{(2)}}}
=\frac{1}{E\left(\vphi,B,F_{*}^{(2)}\right)},
\eeq
where  $E\left(\vphi,B,F_{*}^{(2)}\right)$ is the error of best approximation by elements
from $B$ in the space
 $F^{(2)}_{*}$.
 To have $\C<\iy$, we assume that $E\left(\vphi,B,F_{*}^{(2)}\right)>0$.
 Note that the condition $f\ne 0$ in \eqref{E2.8*} is needed because
 $\|0\|_{F^{(1)},\mathrm{triv}}=0$.

So the problem of finding $\C$ is reduced to finding sharp or asymptotically
 sharp value of \linebreak
 $E\left(\vphi,B,F_{*}^{(2)}\right)$.
 These problems have attracted much attention in classical and contemporary
 approximation theory, having started with classical works by Chebyshev,
 A. A. Markov, V. A. Markov, Zolotarev, and Bernstein (see, e.g., \cite{T1963, A1965,
 DL1993, MMR1994, BE1995}).
 Much contributions to this topic,
 including asymptotic relations \eqref{E1.1} through \eqref{E1.4}
  and their generalizations,
  have been made
 in \cite{B1913, B1938, K1938, B1946, B1947(lim), T1963, A1965, G1982, G1991, G1992,
  G2000, S2003, G2010, CV2010, CL2013}
  and references therein.
 \vspace{.15in}\\
\emph{C. Approximation on Classes of Elements.}
Let $F^{(1)}=F^{(2)},\,D=I$, and let us choose
$
\|\cdot\|_{F^{(2)}}=\|\cdot\|_{F^{(2)},\mathrm{triv}}.
$
Let $\|\cdot\|_{F^{(1)}}$ be a $\beta_1$-seminorm of the form
\ba
\|f\|_{F^{(1)}}:=E\left(f,G,F^{(1)}_{*}\right)
=\inf_{g\in G}\|f-g\|^{*}_{F^{(1)}},\qquad f\in F^{(1)}_{*} ,
\ea
where $G$ is a subspace of $F^{(1)}$ and $F^{(1)}_{*}$ is a vector space
$F^{(1)}$ equipped with a certain
$\beta_1$-norm $\|\cdot\|^{*}_{F^{(1)}}$
on $F^{(1)},\,\be_1\in(0,1]$. Then
$
\C=\sup_{f\in B}E\left(f,G,F^{(1)}_{*}\right)
$
is the sharp constant of approximation by elements from $G$ in the
$\beta_1$-normed space $F_*^{(1)}$ on a class $B$.
Sharp or asymptotically sharp values of this constant of approximation
by trigonometric or algebraic polynomials, splines, and entire functions
of exponential type   have been found since the 1930s for various classes,
see, e.g.,
\cite{AK1937, N1946, B1946(Nik), B1947(const), T1963, A1965, K1971,
D1975, K1991, G1992,
 G2000, G2003, G2018a, G2018b} and references therein.
\vspace{.15in}\\
We also discuss the following important problem in analysis.\vspace{.15in}\\
\emph{D. Norms of Operators.}
If $F^{(j)}$ is a $\beta_j$-normed space, $j=1,\,2,$
and $B=F^{(2)}$, then
$\C=\|D\|$ is the "norm" of the operator
$D:F^{(2)}\to F^{(1)}$.
Finding or estimating $\C$ is one of the major problems in functional analysis,
Fourier analysis, real and complex analysis, and other areas.
Below we discuss just a few examples of such a problem.\vspace{.15in}\\
\emph{D1. Lebesgue constants in Fourier analysis}, see, e.g., \cite[Vol. I, Sect. 2.12]{Z1968}
and \cite{L2006}.\\
\emph{D2. Summability methods}, see, e.g., \cite[Ch. 8]{T1963}
and \cite[Ch. 8]{TB2004}.\\
\emph{D3. Hausdorf-Young-F. Riesz inequalities}, see, e.g.,
\cite[Vol. 2, Sects. 12.2 and 16.3]{Z1968}.\\
\emph{D4. Embedding and extension operators}, see, e.g., \cite{N1969, BG1983, BB2012}.\\
\emph{D5. Hilbert and maximal Hardy-Littlewood operators},
see, e.g., \cite[Sects 4.3 and 6.6]{G1987}.

\section{Problems and Conditions}\label{S3}
\setcounter{equation}{0}
\noindent
Here, we discuss major objectives of this paper and introduce conditions that
are needed for our results.
\subsection{Problems}\label{S3.1}
Let
\bna
&&\C:=\C\left(D,B,F^{(1)},F^{(2)}\right):=\sup_{f\in B,\,\|f\|_{F^{(2)}}> 0}
\frac{\|D(f)\|_{F^{(1)}}}{\|f\|_{F^{(2)}}},\label{E2.8a}\\
&&\C_n:=\C\left(D_n,B_n,F^{(1)}_n,F^{(2)}_n\right)
:=\sup_{f_n\in B_n,\,\|f_n\|_{F^{(2)}_n}> 0}
\frac{\|D_n(f_n)\|_{F^{(1)}_n}}{\|f_n\|_{F^{(2)}_n}}\label{E2.8b}
\ena
be two sharp constants over two families of components
$\left\{D,B,F^{(1)},F^{(2)}\right\}$ and
$\left\{D_n,B_n,F^{(1)}_n,F^{(2)}_n\right\}$
with a running parameter $n\in\N$ attached to the second family.

Our goal is to obtain asymptotic relations of the form:
\bna
&&\C\le \liminf_{n\to\iy}\C_n, \label{E2.8}\\
&&\C\ge \limsup_{n\to\iy}\C_n, \label{E2.9}\\
&&\C= \lim_{n\to\iy}\C_n.\label{E2.10}
\ena
In special cases when $\C$ and $\C_n$ are reduced to approximation
 of individual
elements or approximation on classes of elements, relations like
\eqref{E2.8}--\eqref{E2.10}
were obtained by Bernstein \cite{B1946(Nik), B1946(trig), B1946,
B1947(const), B1947(lim)}
 for approximation by
 univariate trigonometric or algebraic polynomials
and entire functions of exponential type.
The author extended these results to multivariate approximation
\cite{G1982, G1991} and to a more
 general setting \cite{G1992, G2000}.
Burenkov and Goldman \cite{BG1983} discussed relations like
\eqref{E2.8}--\eqref{E2.10} between norms of
 linear operators on periodic and nonperiodic spaces.
 Relations like \eqref{E2.8}--\eqref{E2.10} in Markov-Bernstein
 type inequalities of
 different metrics
 for trigonometric and algebraic polynomials and entire functions
 of exponential type
  were proved by Taikov \cite{T1965, T1993}, Gorbachev \cite{G2005, Gor2018},
   Levin and Lubinsky \cite{LL2015a,LL2015b},
  the author and Tikhonov \cite{GT2017},
  Dai, Gorbachev, and Tikhonov \cite{DGT2018, DGT2019},
   and the author
  \cite{G2017,G2018,G2019a,G2019b}.

 Below we introduce  some conditions on families
 $\left\{D,B,F^{(1)},F^{(2)}\right\}$ and
 $\left\{D_n,B_n,F^{(1)}_n,F^{(2)}_n\right\},\,n\in\N$,
 that are needed for the proof of relations \eqref{E2.8}--\eqref{E2.10}.

 \subsection{Conditions and their Properties}\label{S3.2}
The first group of conditions is needed for the proof of relation \eqref{E2.8}.
\vspace{.15in}\\
\emph{Condition C1).} A space $F^{(1)}_n$ is equipped with another "norm"
  $\|\cdot\|_{F^{(1)}_n}^*,\,n\in\N,$ and
  $\left\{\|\cdot\|_{F^{(1)}_n}\right\}_{n=1}^\iy$ possesses the GC-property
with respect to the sequence
$\left\{\|\cdot\|_{F^{(1)}_n}^*\right\}_{n=1}^\iy$ (see Definition \ref{D2.4}).
\vspace{.15in}\\
\emph{Condition C2).} There exists a sequence of operators
$L_n^{(1)}:F^{(1)}\to F_n^{(1)},\,n\in\N$.

Certain properties of operators
$L_n^{(1)},\,n\in\N$, are discussed in the next two conditions.
\vspace{.15in}\\
\emph{Condition C3).} For every $f\in B$ with $\|f\|_{F^{(2)}}> 0,$
 there exists a sequence $f_n\in B_n,\,n\in\N$,
such that
\beq \label{E2.11}
\limsup_{n\to\iy}\|f_n\|_{F_n^{(2)}}\le \|f\|_{F^{(2)}}
\eeq
and
\beq \label{E2.12}
\lim_{n\to\iy}\left\|L_n^{(1)}(D(f))-D_n(f_n)\right\|_{F_n^{(1)}}^*=0.
\eeq
\emph{Condition C4).} For every $h\in F^{(1)}$,
\beq \label{E2.13}
\liminf_{n\to\iy}\left\|L_n^{(1)}(h)\right\|_{F_n^{(1)}}\ge\|h\|_{F^{(1)}}.
\eeq

In certain cases a slightly different condition is needed.
\vspace{.15in}\\
\emph{Condition C4\textprime).} Let
$B\subseteq F^{(1)}\cap F^{(2)}$ and
$D(B)\subseteq B$.
Then for every $h\in B$,
inequality \eqref{E2.13} holds true.

A group of stronger conditions marked with asterisks
 is given below.\vspace{.15in}\\
\emph{Condition C1*).} A space $F^{(j)}_n$ is equipped with another "norm"
  $\|\cdot\|_{F^{(j)}_n}^*,\,n\in\N,$ and
  $\left\{\|\cdot\|_{F^{(j)}_n}\right\}_{n=1}^\iy$ possesses the GC-property
with respect to the sequence
$\left\{\|\cdot\|_{F^{(j)}_n}^*\right\}_{n=1}^\iy,\,j=1,\,2$
 (see Definition \ref{D2.4}).\vspace{.15in}\\
\emph{Condition C2*).} There exist two sequences of operators
$L_n^{(j)}:F^{(j)}\to F_n^{(j)},\,n\in\N,\, j=1,2$.

Certain properties of operators
$L_n^{(j)},\,n\in\N,\, j=1,2$, are discussed
in the next two conditions.
\vspace{.15in}\\
\textit{Condition C3*).} For every $f\in B$ with $\|f\|_{F^{(2)}}> 0,$
there exists a sequence $f_n\in B_n,\,n\in\N$,
such that
\bna
&&\lim_{n\to\iy}\left\|L_n^{(1)}(D(f))-D_n(f_n)\right\|_{F_n^{(1)}}^*=0,\label{E2.14}\\
&&\lim_{n\to\iy}\left\|f_n-L_n^{(2)}(f)\right\|_{F_n^{(2)}}^*=0.\label{E2.15}
\ena
\emph{Condition C4*).}  For every $h\in F^{(2)}$ with $\|h\|_{F^{(2)}}>0$,
\beq \label{E2.16}
\limsup_{n\to\iy}\left\|L_n^{(2)}(h)\right\|_{F_n^{(2)}}\le\|h\|_{F^{(2)}}.
\eeq

\begin{proposition}\label{P2.7}
Conditions C1*) through C4*) imply Conditions C1) through C3).
\end{proposition}
\proof
It suffices to show that C1*), C2*),  \eqref{E2.15} of C3*), and \eqref{E2.16}
of C4*) imply
\eqref{E2.11}. By C2*), C3*), and C4*), for every $f\in B$ there exists a sequence
 $\{f_n\}_{n=1}^\iy$ that satisfies \eqref{E2.15}, where  $\left\{L_n^{(2)}\right\}_{n=1}^\iy$
 satisfies \eqref{E2.16}. Then using C1*), \eqref{E2.15}, and \eqref{E2.16}, we obtain
 \ba
 \limsup_{n\to\iy}\|f_n\|_{F_n^{(2)}}
 =\limsup_{n\to\iy}\left\|L_n^{(2)}(f)+\left(f_n-L_n^{(2)}(f)\right)\right\|_{F_n^{(2)}}
 =\limsup_{n\to\iy}\left\|L_n^{(2)}(f)\right\|_{F_n^{(2)}}
 \le \|f\|_{F^{(2)}}.
 \ea
 Hence \eqref{E2.11} follows.\hfill $\Box$\vspace{.15in}\\
 Next, we introduce the second group of conditions that is
  needed for the proof of relation
  \eqref{E2.9}.
 \vspace{.15in}\\
  \emph{Condition C5).} A space $F^{(j)}_n$ is equipped with another "norm"
  $\|\cdot\|_{F^{(j)}_n}^*,\,n\in\N,$ and
  $\left\{\|\cdot\|_{F^{(j)}_n}\right\}_{n=1}^\iy$ possesses the GC-property
with respect to the sequence
$\left\{\|\cdot\|_{F^{(j)}_n}^*\right\}_{n=1}^\iy,\,\,j=1,\,2$
 (see Definition \ref{D2.4}).\vspace{.15in}\\
 \emph{Condition C6).} There exist two sequences of operators
$l_n^{(j)}:F^{(j)}\to F_n^{(j)},\,n\in\N,\,j=1,2$, such that the
following two properties hold:\vspace{.15in}\\
\emph{Condition C6a).}
For every $h\in F^{(1)}$ with $\|h\|_{F^{(1)}}>0$,
\beq \label{E2.17}
\limsup_{n\to\iy}\left\|l_n^{(1)}(h)\right\|_{F_n^{(1)}}\le\|h\|_{F^{(1)}}.
\eeq
\emph{Condition C6b).} For every $h\in F^{(2)}$,
\beq\label{E2.18}
\liminf_{n\to\iy}\left\|l_n^{(2)}(h)\right\|_{ F^{(2)}_n}
\ge \|h\|_{F^{(2)}}.
\eeq
\emph{Condition C7).}
There is a constant $C_1\in[1,\iy)$ such that
for every $h\in F^{(2)}$ with $\|h\|_{F^{(2)}}>0$,
\beq \label{E2.24}
\limsup_{n\to\iy}\left\|l_n^{(2)}(h)\right\|_{F_n^{(2)}}\le C_1\|h\|_{F^{(2)}}.
\eeq
Next, we need one of the following two conditions associated with relations
like \eqref{E2.6} and \eqref{E2.7}.\vspace{.15in}\\
\emph{Condition C8).} There is a set
$B_n^\prime \subseteq\left\{f_n\in B_n:\|f_n\|_{F^{(2)}_n}=1\right\}$
such that
the following relations hold:
\bna
\C_n=\sup_{f_n\in B_n^\prime}\|D_n(f_n)\|_{F^{(1)}_n}
=\sup_{f_n\in B_n^\prime}\frac{\|D_n(f_n)\|_{F_n^{(1)}}}{\|f_n\|_{F_n^{(2)}}},
\qquad n\in\N.\label{E2.23}
\ena
\emph{Condition C9).}
There is a set
$B_n^\prime \subseteq\left\{f_n\in B_n:\|f_n\|_{F^{(2)}_n}>0,\,
\|D_n(f_n)\|_{ F^{(1)}_n}=1\right\}$
such that
the following relations hold:
\bna
\C_n=\sup_{f_n\in B_n^\prime}
\frac{1}{\|f_n\|_{F^{(2)}_n}}
=\sup_{f_n\in B_n^\prime}\frac{\|D_n(f_n)\|_{F_n^{(1)}}}{\|f_n\|_{F_n^{(2)}}},
\qquad n\in\N.\label{E2.26}
\ena
\emph{Condition C10).} (Compactness condition)
For any sequence $\{f_n\}_{n=1}^\iy,\,f_n\in B_n^\prime\setminus \{0\},\,n\in\N$,
there exist a sequence $\{n_k\}_{k=1}^\iy$ of natural numbers and an element $f\in B$
such that
\bna
&&\lim_{k\to\iy}\left\|D_{n_k}\left(f_{n_k}\right)
-l_{n_k}^{(1)}(D(f))\right\|_{F_{n_k}^{(1)}}^*=0,
\label{E2.20}\\
&&\lim_{k\to\iy}\left\|f_{n_k}-l_{n_k}^{(2)}(f)\right\|_{F_{n_k}^{(2)}}^*=0.
\label{E2.21}
\ena
Further, we discuss a different group of conditions with a weak compactness
condition that replaces  C10).\vspace{.15in}\\
\emph{Condition C11).} There exist two families of operators
$l_{n,s}^{(j)}:F_n^{(j)}\to F_s^{(j)},\,n\in\N,\,s\in\N,\,j=1,2$,
 such that the following property holds:\vspace{.15in}\\
\emph{Condition C11a).}
\beq\label{E2.27}
\sup_{n\in\N,\,s\in\N}\left\|l_{n,s}^{(2)}\right\|\le 1.
\eeq
 We recall that the "norm" of an operator is defined by \eqref{E2.4a}.
 \vspace{.15in}\\
\emph{Condition C12).} (Weak compactness condition)
For any sequence $\{f_n\}_{n=1}^\iy,\,f_n\in B_n^\prime,\,n\in\N$,
there exists $f\in B$ such that
for any fixed $s\in\N$,
there exists a sequence $\{n_k\}_{k=1}^\iy$ of natural numbers
 ($n_k=n_k(s),\,k\in\N$)
such that
\bna
&&\lim_{k\to\iy}\left\|l_{n_k,s}^{(1)}\left(D_{n_k}
\left(f_{n_k}\right)\right)-l_{s}^{(1)}
(D(f))\right\|_{F_{s}^{(1)}}^*=0,\label{E2.29}\\
&&\lim_{k\to\iy}\left\|l_{n_k,s}^{(2)}(f_{n_k})-l_{s}^{(2)}(f)
\right\|_{F_{s}^{(2)}}^*=0,\label{E2.30}
\ena
and, in addition,
\beq\label{E2.31}
\left\|D_{n_k}\left(f_{n_k}\right)\right\|_{F_{n_k}^{(1)}}
=\left\|l_{n_k,s}^{(1)}\left(D_{n_k}
\left(f_{n_k}\right)\right)\right\|_{F_{s}^{(1)}},
 \quad k\in\N.
\eeq
\emph{Condition C13).}
There is a constant $C_2\in(0,\iy)$ such that
for every $h_n\in F_n^{(2)}$,
\beq \label{E2.24a}
\limsup_{n\to\iy}\|h_n\|_{F_n^{(2)}}
\le C_2\limsup_{s\to\iy}\limsup_{n\to\iy}
\left\|l_{n,s}^{(2)}(h_n)\right\|_{F_s^{(2)}}.
\eeq
In Section \ref{S5.3} we will need modified conditions
C6a), C10), and C12) marked with asterisks below. \vspace{.15in}\\
\emph{Condition C6*a).}
There exists a sequence of operators
$l_n^{(1)}:F^{(1)}\to F_n^{(1)},\,n\in\N$, such that
for every $h\in F^{(1)}$ with $\|h\|_{F^{(1)}}>0$,
\beq \label{E3.12a}
\limsup_{n\to\iy}\left\|l_n^{(1)}(h)\right\|_{F_n^{(1)}}\le\|h\|_{F^{(1)}}.
\eeq
\emph{Condition C10*a).}
For any sequence $\{f_n\}_{n=1}^\iy,\,f_n\in B_n^\prime\setminus \{0\},\,n\in\N$,
there exists a sequence $\{n_k\}_{k=1}^\iy$ of natural numbers and an element $f\in B$
such that
\beq\label{E3.17a}
\lim_{k\to\iy}\left\|D_{n_k}\left(f_{n_k}\right)
-l_{n_k}^{(1)}(D(f))\right\|_{F_{n_k}^{(1)}}^*=0.
\eeq
\emph{Condition C10*b).}  For
any sequence $\{f_n\}_{n=1}^\iy,\,f_n\in B_n^\prime\setminus \{0\},\,n\in\N$,
for
the sequence $\{n_k\}_{k=1}^\iy$ of natural numbers from Condition C10*a),
 and for the element $f\in B$ from Condition C10*a),
there exists a sequence of operators
$l_n^{(2)}:F^{(2)}\to F_n^{(2)},\,n\in\N$, such that
Conditions C6b) and C7) are satisfied. In addition,
the following relation holds:
\beq\label{E3.18a}
\lim_{k\to\iy}\left\|f_{n_k}-l_{n_k}^{(2)}(f)\right\|_{F_{n_k}^{(2)}}^*=0.
\eeq
\emph{Condition C12*a).}
There exists a family of operators
$l_{n,s}^{(1)}:F_n^{(1)}\to F_s^{(1)},\,n\in\N,\,s\in\N$, such that
for any sequence $\{f_n\}_{n=1}^\iy,\,f_n\in B_n^\prime\setminus \{0\},\,n\in\N$,
there exists $f\in B$ such that
for any fixed $s\in\N$,
there exists a sequence $\{n_k\}_{k=1}^\iy$ of natural numbers
 ($n_k=n_k(s),\,k\in\N$)
such that
\bna
&&\lim_{k\to\iy}\left\|l_{n_k,s}^{(1)}
\left(D_{n_k}\left(f_{n_k}\right)\right)-l_{s}^{(1)}
(D(f))\right\|_{F_{s}^{(1)}}^*=0,\label{E3.20a}\\
&&\left\|D_{n_k}\left(f_{n_k}\right)\right\|_{F_{n_k}^{(1)}}
=\left\|l_{n_k,s}^{(1)}\left(D_{n_k}
\left(f_{n_k}\right)\right)\right\|_{F_{s}^{(1)}},
 \quad k\in\N.\label{E3.22a}
\ena
\emph{Condition C12*b).} For
any sequence $\{f_n\}_{n=1}^\iy,\,f_n\in B_n^\prime\setminus \{0\},\,n\in\N$,
for
the sequence $\{n_k\}_{k=1}^\iy$ of natural numbers from Condition C12*a)
 and the element $f\in B$
from Condition C12*a),
there exists a sequence of operators
$l_n^{(2)}:F^{(2)}\to F_n^{(2)},\,n\in\N$, such that
Conditions C6b) and C7) are satisfied and
there exists a sequence of operators
$l_{n,s}^{(2)}:F_n^{(2)}\to F_s^{(2)},\,n\in\N,\,s\in\N$, such that
Conditions C11a) and C13) are satisfied.
In addition, the following relation holds:
\beq\label{E3.21a}
\lim_{k\to\iy}\left\|l_{n_k,s}^{(2)}(f_{n_k})-l_{s}^{(2)}(f)
\right\|_{F_{s}^{(2)}}^*=0,\qquad s\in\N.
\eeq

\begin{remark}\label{R2.7d}
Condition C4\textprime) includes inequality \eqref{E2.13}
with a weaker requirement $h\in B$,
so if $B\subseteq F^{(1)}\cap F^{(2)}$ and
$D(B)\subseteq B$
(for example, if $D=I$ is an imbedding operator),
then Condition C4\textprime) can efficiently replace Condition C4),
see Theorem \ref{T4.1aa} and Corollaries \ref{C2.11},
\ref{C2.14a} (i), \ref{C2.19} (i)  for details.
\end{remark}

\begin{remark}\label{R2.7c}
Modified compactness conditions C10*b) and C12*b)
 are characterized
 by the dependence of operators $\left\{l_n^{(2)}\right\}_{n\in\N}$
  and $\left\{l_{n,s}^{(2)}\right\}_{n\in\N,s\in\N}$ on a
 sequence $\{f_n\}_{n=1}^\iy$
 ($f_n\in B_n^\prime\setminus \{0\},\,n\in\N$),
 the sequence $\{n_k\}_{k=1}^\iy$, and $f\in B$,
 while these operators are defined in C6) and C11)
 independently of C10) and C12).
 \end{remark}

 \begin{remark}\label{R2.7b}
The set $B_n^\prime,\,n\in\N$, mentioned in the Compactness Conditions
C10), C12), C10*a), C10*b) C12*a), and C12*b)
 is defined in C8) or C9). In other words,
if any of these Compactness Conditions is satisfied,
 then either C8) or C9)
must be satisfied as well.
Actually, Conditions C8) and C9) state only the existence
of $B_n^\prime,\,n\in\N$. So these sets should be defined for each concrete
problem separately (see Chapters \ref{S5} and \ref{S7}).
\end{remark}

\begin{remark}\label{R2.7a}
  Conditions C4*), C6a), C6*a), and C7) can be replaced by stronger
   conditions for the "norms" of the corresponding operators. For example, if
   the condition $\limsup_{n\to \iy}\left\|L_{n}^{(2)}\right\|\le 1$
   is satisfied, then
   C4*) is satisfied as well. However in general, C4*) does not imply this condition.
   \end{remark}
  \noindent
Several conditions defined above are satisfied for certain special "norms."
\begin{proposition}\label{P2.8}
Let $j=1$ or $j=2$ and let $\be_j\in(0,1]$ be independent of $n\in\N$.
 Next, let $F_n^{(j)}$ be equipped with two "norms"
 $\left\|\cdot\right\|_{F_{n}^{(j)}}$ and
 $\left\|\cdot\right\|_{F_{n}^{(j)}}^*,\,n\in\N$
 (see Conditions C1), C1*), and C5)). In addition, let
 $\left\|\cdot\right\|_{F_{n}^{(j)}},\,n\in \N$,
be one of the following "norms:"
\begin{itemize}
\item[a)]
$\left\|\cdot\right\|_{F_{n}^{(j)}}
=\left\|\cdot\right\|_{F_{n}^{(j)}}^*$
is a $\be_j$-norm or a $\be_j$-seminorm,
$n\in\N$.
\item[b)]
 $\left\|\cdot\right\|_{F_{n}^{(j)}}
 =\left\|\cdot\right\|_{F_{n}^{(j)}}^*
:=\|\cdot\|_{F_n^{(j)},\mathrm{triv}},\,n\in\N$.
\item[c)]
$\left\|f_n\right\|_{F_{n}^{(j)}}
=\|\vphi_n-f_n\|_{F_n^{(j)}}^*,\,f_n\in F_n^{(j)}$,
where $\vphi_n $ is a fixed element
of $F_n^{(j)}$, and
$\|\cdot\|_{F_n^{(j)}}^*$ is a
$\be_j$-norm or a $\be_j$-seminorm,
$n\in\N$.
\end{itemize}
Then the following statements are valid:\\
(i) If $j=1$, then Condition C1) is satisfied.\\
(ii) If $j=1,\,2$, then Conditions C1*) and C5) are satisfied.
\end{proposition}
\proof
In case of "norms" a) and b), the corresponding sequence
$\left\{\|\cdot\|_{F_n^{(j)}}\right\}_{n=1}^\iy
=\left\{\|\cdot\|_{F_n^{(j)}}^*\right\}_{n=1}^\iy$
possesses
the C-property by statements (ii) and (iii)
of Proposition \ref{P2.3}.
In case of "norm" c), the sequence
$\left\{\|\cdot\|^*_{F_n^{(j)}}\right\}_{n=1}^\iy$ possesses
the C-property by statement (ii)
of Proposition \ref{P2.3}. Then
$\left\{\|\cdot\|_{F_n^{(j)}}\right\}_{n=1}^\iy$ possesses
the GC-property with respect to
$\left\{\|\cdot\|^*_{F_n^{(j)}}\right\}_{n=1}^\iy$.
Indeed, like in Example \ref{Ex2.5},
for any
$h_n\in F_n^{(j)},\,g_n\in F_n^{(j)},\,n\in\N$, such that
\ba
\lim_{n\to\iy}\|g_n\|_{F_n^{(j)}}^*
=\lim_{n\to\iy}\|-g_n\|_{F_n^{(j)}}^*=0,
\ea
we have
\bna\label{E2.4c}
\liminf_{n\to\iy}\|h_n+g_n\|_{F_n^{(j)}}
&=&\liminf_{n\to\iy}\|(\vphi_n-h_n)+(-g_n)\|_{F_n^{(j)}}^*\nonumber\\
&=&\liminf_{n\to\iy}\|\vphi_n-h_n\|_{F_n^{(j)}}^*\nonumber\\
&=&\liminf_{n\to\iy}\|h_n\|_{F_n^{(j)}}.
\ena
Similar relations hold true if we replace $\liminf$
in \eqref{E2.4c} with $\limsup$.
Hence \eqref{E2.2} and \eqref{E2.3}
 are valid for $F_n$ replaced with $F_n^{(j)}$.
 \hfill $\Box$
\begin{proposition}\label{P2.9}
(i) Let $F_{n}^{(j)},\,n\in\N$,  be a $\be_j$-normed or
$\be_j$-seminormed space, where
$\be_j$ is independent of $n,\,j=1,2.$
In addition, let $B_n$ be a subspace of $F_{n}^{(2)}$
and let $D_n$ be a $1$-homogeneous operator
(that is, $D_n(\al f)=\al D_n(f),\,\al>0$),
$n\in\N$. Then Conditions C8) and C9)
are satisfied.\\
(ii) If  $\left\|\cdot\right\|_{F_{n}^{(2)}}
=\|\cdot\|_{F_n^{(2)},\mathrm{triv}},\,n\in\N$,
then C8) is satisfied.\\
(iii) If  $\left\|\cdot\right\|_{F_{n}^{(1)}}
=\|\cdot\|_{F_n^{(1)},\mathrm{triv}},\,n\in\N$,
then C9) is satisfied.
\end{proposition}
\noindent
\proof
Let us define a set $G_n,\,n\in\N$, by
\ba
G_n:=\left\{\begin{array}{ll}
\left\{f_n\in B_n:\|f_n\|_{F^{(2)}_n}=1\right\} & \mbox{for Condition C8)},\\
\left\{f_n\in B_n:\|f_n\|_{F^{(2)}_n}>0,\,
\|D_n(f_n)\|_{ F^{(1)}_n}=1\right\}&\mbox{for Condition C9)}.
\end{array}\right.
\ea
It follows immediately from Proposition \ref{P2.6} that relation
\eqref{E2.23} or relation \eqref{E2.26} holds true
with $B_n^\prime$ replaced by $G_n,\,n\in\N$.
Therefore, there exists $B_n^\prime\subseteq G_n,\,n\in\N$, such that
\eqref{E2.23} or \eqref{E2.26} is valid
for all premises of statements (i), (ii), and (iii) of Proposition \ref{P2.9}
 \hfill $\Box$

 \begin{remark}\label{P2.10}
  For the convenience of the reader, we present here
  certain explanations of the conditions of this section.
All the conditions of this section are divided into
 groups I and II that are needed
 for the proofs of inequalities \eqref{E2.8} and \eqref{E2.9},
 respectively.
 Group I contains C1)--C4) and certain modified conditions
  C1*)--C4*) with the "leading" conditions
 C3) and C3*) that, roughly speaking,
 provide "good" approximation $f_n\in B_n,\,n\in\N,$ to
 a given element $f\in B$.
 Group II contains C5)--C13) and certain modified conditions
 C6*a), C10*a), C10*b), and C12*b) with the "leading" conditions
 C10) and C12). Roughly speaking,
 the  compactness and weak compactness conditions
 guarantee the existence of element $f\in B$ and a
 subsequence $\{f_{n_k}\}_{k=1}^\iy$ of a given sequence
 $\{f_{n}\}_{n=1}^\iy$ with "good" approximation to $f$.
 In addition, note that Condition C5), which is identical to C1*),
 is included into the condition list to make groups I and II
 disjoint.

 Next, the "norms"
 $\|\cdot\|_{F^{(j)}_n}$ and
 $\|\cdot\|_{F^{(j)}_n}^*,\,n\in\N,\,j=1,2$,
 coincide for inequalities for approximating elements
  (Problem A).
 However, both "norms" are different for
 approximation Problems B and C because
 the following definition is needed
 to discuss the error of best approximation:
 $\left\|f_n\right\|_{F_{n}^{(j)}}
:=\|\vphi_n-f_n\|_{F_n^{(j)}}^*,\, n\in\N$
(cf. Proposition \ref{P2.8}).
 \end{remark}

 \section{Asymptotic Relations between Sharp Constants in Vector Spaces}\label{S4}
 \setcounter{equation}{0}
 \noindent
 Here, we prove
 asymptotic relations
 \eqref{E2.8}, \eqref{E2.9}, and \eqref{E2.10}.
 We first discuss asymptotic relations of the form
 \beq\label{E2.32}
\C=\C\left(D,B,F^{(1)},F^{(2)}\right)\le \liminf_{n\to\iy}
\C\left(D_n,B_n,F^{(1)}_n,F^{(2)}_n\right)=\liminf_{n\to\iy}\C_n.
\eeq
\begin{theorem}\label{T2.10}
If Conditions C1) through C4) are satisfied, then \eqref{E2.32} holds true.
\end{theorem}
\proof
Let $f\in B$ be a fixed element with $\|f\|_{F^{(2)}}> 0$.
Then $D(f)\in F^{(1)}$ and using Condition C2) and
 inequality \eqref{E2.13} of Condition C4),
we have
\beq \label{E2.33}
\|D(f)\|_{F^{(1)}}\le \liminf_{n\to\iy}\left\|L_n^{(1)}(D(f))\right\|_{F_n^{(1)}}.
\eeq
Next by C2) and C3), there exists a sequence $f_n\in B_n, n\in\N,$
such that  \eqref{E2.11} and \eqref{E2.12} hold. Then using C1), we obtain from
\eqref{E2.33}, \eqref{E2.12}, and \eqref{E2.11}
\ba
&&\|D(f)\|_{F^{(1)}}
\le \liminf_{n\to\iy}\left\|D_n(f_n)+\left(L_n^{(1)}
(D(f))-D_n(f_n)\right)\right\|_{F_n^{(1)}}
=\liminf_{n\to\iy}\left\|D_n(f_n)\right\|_{F_n^{(1)}}\\
&&\le \liminf_{n\to\iy}\left(\|f_n\|_{F^{(2)}_n} \C_n\right)
\le \limsup_{n\to\iy}\|f_n\|_{F^{(2)}_n} \liminf_{n\to\iy}\C_n
\le \|f\|_{F^{(2)}} \liminf_{n\to\iy}\C_n.
\ea
Therefore, $\|D(f)\|_{F^{(1)}}/\|f\|_{F^{(2)}}\le \liminf_{n\to\iy}\C_n$
and \eqref{E2.32} is established.\hfill $\Box$
\vspace{.15in}\\
The following theorem can be proved similarly to Theorem \ref{T2.10}.
\begin{theorem}\label{T4.1aa}
If Conditions C1) through C3) and C4\textprime) are satisfied,
then \eqref{E2.32} holds true.
\end{theorem}
\begin{corollary}\label{C2.11}
If Conditions C4) (or C4\textprime)) and C1*) through C4*)
are satisfied, then \eqref{E2.32} holds
true.
\end{corollary}
\noindent
The corollary follows immediately from Theorems \ref{T2.10},
\ref{T4.1aa},
and Proposition \ref{P2.7}.

Next, we discuss asymptotic relations of the form
 \beq\label{E2.34}
\C=\C\left(D,B,F^{(1)},F^{(2)}\right)\ge \limsup_{n\to\iy}
\C\left(D_n,B_n,F^{(1)}_n,F^{(2)}_n\right)=\limsup_{n\to\iy}\C_n.
\eeq
In the following Theorems \ref{T2.12}, \ref{T2.13}, and  \ref{T4.6}
we take account of Remark \ref{R2.7b}.

\begin{theorem}\label{T2.12}
Let $\C_n$ be a finite number for every $n\in\N$, and let either
Conditions C7), C8) or Condition C9) be satisfied.
 In addition, let  C5), C6), and C10) be satisfied.
   Then \eqref{E2.34} holds true.
\end{theorem}
\proof
Assume without loss of generality that $\C<\iy$.
Next, we assume that Conditions C7) and C8) are satisfied.
Since $\C_n<\iy$, we see from C8) that for any fixed $\vep>0$
there exists $f_n=f_{n,\vep}\in B_n^\prime$ with $\|f_n\|_{F^{(2)}_n}=1$
such that
\beq\label{E2.35}
\C_n<\frac{\|D_n(f_n)\|_{F_n^{(1)}}}{\|f_n\|_{F_n^{(2)}}}+\vep,\qquad n\in\N.
\eeq
If Condition C9) is satisfied, then for any fixed $\vep>0$
there exists $f_n=f_{n,\vep}\in B_n^\prime$ with
 $\|f_n\|_{F^{(2)}_n}>0$ and $\|D_n(f_n)\|_{F^{(1)}_n}=1$
such that \eqref{E2.35} holds true as well.

Let $\{n_p\}_{p=1}^\iy$ be a sequence such that
\beq\label{E2.36}
\limsup_{n\to\iy} \C_n=\lim_{p\to\iy} \C_{n_p}.
\eeq
 Next, we apply Conditions C6), C6a), C6b) and C10) to a sequence
  $\{f_{n_p}\}_{p=1}^\iy$. Then there exist
  a subsequence $\{f_{n_{p_k}}\}_{k=1}^\iy$ and an element $f\in B$
such that
\bna
&&\lim_{k\to\iy}\left\|D_{n_{p_k}}\left(f_{n_{p_k}}\right)-l_{n_{p_k}}^{(1)}(D(f))
\right\|_{F_{n_{p_k}}^{(1)}}^*=0,
\label{E2.37}\\
&&\lim_{k\to\iy}\left\|f_{n_{p_k}}-l_{n_{p_k}}^{(2)}(f)\right\|_{F_{n_{p_k}}^{(2)}}^*=0,
\label{E2.38}
\ena
where $l_{n}^{(j)},\,n\in\N,\,j=1,\,2$, satisfy C6a) and C6b).
In particular, $\|f\|_{F^{(2)}}<\iy$ since $B\subseteq F^{(2)}$.

If C7) and C8) are satisfied, then we have from \eqref{E2.38},  C5), and
\eqref{E2.24} of C7)
\ba
1=\lim_{k\to\iy}\left\|f_{n_{p_k}}\right\|_{F_{n_{p_k}}^{(2)}}
&=&\lim_{k\to\iy}\left\|l^{(2)}_{n_{p_k}}(f)
 +\left(f_{n_{p_k}}-l^{(2)}_{n_{p_k}}(f)\right)\right\|_{F_{n_{p_k}}^{(2)}}\\
&=&\lim_{k\to\iy}\left\|l^{(2)}_{n_{p_k}}(f)\right\|_{F_{n_{p_k}}^{(2)}}
\le C_1\|f\|_{F^{(2)}}.
\ea
Hence there exist $k_0\in\N$ and constants $C_3>0$ and $C_4>0$ such that
\beq\label{E2.39}
\|f\|_{F^{(2)}}\ge C_3,\qquad
\inf_{k\ge k_0}\left\|l^{(2)}_{n_{p_k}}(f)\right\|_{F_{n_{p_k}}^{(2)}}\ge C_4.
\eeq

If C9) is satisfied, then we have from \eqref{E2.37} and
Conditions C5), C6a), and C6b)
\ba
1&=&\lim_{k\to\iy}\left\|D_{n_{p_k}}\left(f_{n_{p_k}}\right)\right\|_{F_{n_{p_k}}^{(1)}}
=\lim_{k\to\iy}\left\|l^{(1)}_{n_{p_k}}(D(f))
+\left(D(f)-l^{(1)}_{n_{p_k}}(D(f))\right)\right\|_{F_{n_{p_k}}^{(1)}}\\
&=&\lim_{k\to\iy}\left\|l^{(1)}_{n_{p_k}}(D(f))\right\|_{F_{n_{p_k}}^{(1)}}
\le \|D(f)\|_{F^{(1)}}
\le \C\|f\|_{F^{(2)}}
\le \C\liminf_{k\to\iy}\left\|l^{(2)}_{n_{p_k}}(f)\right\|_{F_{n_{p_k}}^{(2)}}.
\ea
Hence \eqref{E2.39} holds true in this case as well with $C_3=1/\C$.

Next using \eqref{E2.37} and \eqref{E2.38} again, we obtain
by \eqref{E2.35}, \eqref{E2.36}, and C5),
\bna\label{E2.40}
\lim_{k\to\iy} \C_{n_{p_k}}
\le \frac{\limsup_{k\to\iy}\left\|D_{n_{p_k}}\left(f_{n_{p_k}}\right)
\right\|_{F_{n_{p_k}}^{(1)}}}
{\liminf_{k\to\iy}\left\|f_{n_{p_k}}\right\|_{F_{n_{p_k}}^{(2)}}}+\vep
=  \frac{\limsup_{k\to\iy}\left\|l^{(1)}_{n_{p_k}}(D(f))
\right\|_{F_{n_{p_k}}^{(1)}}}
{\liminf_{k\to\iy}\left\|l^{(2)}_{n_{p_k}}(f)\right\|_{F_{n_{p_k}}^{(2)}}}+\vep,
\ena
where $\liminf_{k\to\iy}\left\|l^{(2)}_{n_{p_k}}(f)\right\|_{F_{n_{p_k}}^{(2)}}>0$
by \eqref{E2.39}.
Since $\|f\|_{F^{(2)}}> 0$ by \eqref{E2.39}, it follows from \eqref{E2.36}, \eqref{E2.40},
 and Conditions C6a) and C6b) that
 \beq\label{E2.41}
 \limsup_{n\to\iy} \C_n\le \frac{\|D(f)\|_{F^{(1)}}}{\|f\|_{F^{(2)}}}+\vep
 \le \C+\vep.
 \eeq
 Finally letting $\vep\to 0+$ in \eqref{E2.41}, we arrive at \eqref{E2.34}.
 \hfill $\Box$

 \begin{theorem}\label{T2.13}
Let $\C_n$ be a finite number for every $n\in\N$,
and let either Conditions C7), C8), C13) or
Condition C9) be satisfied.
In addition, let  C5), C6), C11), and C12)
 be satisfied. Then \eqref{E2.34} holds true.
\end{theorem}
\proof
The proof is similar to that of Theorem \ref{T2.12}.
We can assume without loss of generality  that $\C<\iy$.
Next, we assume that Conditions C7), C8),
and C13) are satisfied.
Since $\C_n<\iy$, we see from C8) that for any fixed $\vep>0$
there exists $f_n=f_{n,\vep}\in B_n^\prime$ with $\|f_n\|_{F^{(2)}_n}=1$
such that
\beq\label{E2.42}
\C_n<\frac{\|D_n(f_n)\|_{F_n^{(1)}}}{\|f_n\|_{F_n^{(2)}}}+\vep,
\qquad n\in\N.
\eeq
If Condition C9) is satisfied, then for any fixed $\vep>0$
there exists $f_n=f_{n,\vep}\in B_n^\prime$ with
 $\|f_n\|_{F^{(2)}_n}>0$ and $\|D_n(f_n)\|_{F^{(1)}_n}=1$
such that \eqref{E2.42} holds true as well.

Let $\{n_p\}_{p=1}^\iy$ be a sequence such that
\beq\label{E2.43}
\limsup_{n\to\iy} \C_n=\lim_{p\to\iy} \C_{n_p}.
\eeq
 Next, we apply Conditions C6), C6a), C6b), C11), and C12) to a sequence
  $\{f_{n_p}\}_{p=1}^\iy$. Then
  there exists an element $f\in B$ such that
   for any fixed $s\in\N$ there exists
  a subsequence $\{f_{n_{p_k}}\}_{k=1}^\iy$
  ($n_{p_k}=n_{p_k}(s),\,k\in\N$)
such that
\bna
&&\lim_{k\to\iy}\left\|l_{n_{p_k},s}^{(1)}
\left(D_{n_{p_k}}\left(f_{n_{p_k}}\right)\right)-l_{s}^{(1)}(D(f))
\right\|_{F_{s}^{(1)}}^*=0,
\label{E2.44}\\
&&\lim_{k\to\iy}\left\|l_{n_{p_k},s}^{(2)}\left(f_{n_{p_k}}\right)
-l_{s}^{(2)}(f)\right\|_{F_{s}^{(2)}}^*=0,
\label{E2.45}
\ena
and in addition,
\beq\label{E2.46}
\left\|D_{n_{p_k}}\left(f_{n_{p_k}}\right)\right\|_{F_{n_{p_k}}^{(1)}}
=\left\|l_{n_{p_k},s}^{(1)}\left(D_{n_{p_k}}\left(f_{n_{p_k}}\right)\right)
\right\|_{F_{s}^{(1)}},
 \quad k\in\N.
\eeq
Here, $l_{n}^{(j)}$ and $l_{n,s}^{(j)},\,n\in\N,\,s\in\N,\,j=1,\,2$,
 satisfy C6a), C6b), and C11a).
In particular, $\|f\|_{F^{(2)}}<\iy$ since $B\subseteq F^{(2)}$.
If C7), C8), and C13) are satisfied, then
we have from \eqref{E2.24a} of C13), Condition C5), \eqref{E2.45},
and \eqref{E2.24} of C7)
\ba
1&=&\lim_{k\to\iy}\left\|f_{n_{p_k}}\right\|_{F_{n_{p_k}}^{(2)}}
\le C_2\limsup_{s\to\iy}\limsup_{k\to\iy}
\left\|l_{n_{p_k},s}^{(2)}(f_{n_{p_k}})\right\|_{F_s^{(2)}}\\
&=&C_2 \limsup_{s\to\iy}\limsup_{k\to\iy}
\left\|l^{(2)}_{s}(f)+\left(l_{n_{p_k},s}^{(2)}(f_{n_{p_k}})
-l^{(2)}_{s}(f)\right)\right\|_{F_s^{(2)}}\\
&=& C_2\limsup_{s\to\iy}\left\|l^{(2)}_{s}(f)\right\|_{F_{s}^{(2)}}
\le C_1C_2\|f\|_{F^{(2)}}.
\ea
Hence there exist $s_0\in\N$ and constants
$C_5=1/(C_1C_2)>0$ and $C_6>0$
such that
\beq\label{E2.47}
\|f\|_{F^{(2)}}\ge C_5,
\qquad
\inf_{s\ge s_0}\left\|l^{(2)}_{s}
(f)\right\|_{F_{s}^{(2)}}\ge C_6.
\eeq
Note that the second inequality in \eqref{E2.47}
immediately follows from the first one
and inequality \eqref{E2.18} of Condition C6b).

If C9) is satisfied, then
by C5) we obtain from \eqref{E2.46}, \eqref{E2.44}, and Conditions
C6a) and C6b)
\ba
1&=&\lim_{k\to\iy}\left\|D_{n_{p_k}}
\left(f_{n_{p_k}}\right)\right\|_{F_{n_{p_k}}^{(1)}}
=\lim_{k\to\iy}\left\|l^{(1)}_{n_{p_k},s}\left(D_{n_{p_k}}
\left(f_{n_{p_k}}\right)\right)\right\|_{F_{s}^{(1)}}\\
&=&\left\|
l^{(1)}_{n_{p_k},s}\left(D_{n_{p_k}}
\left(f_{n_{p_k}}\right)\right)
+\left(l_s^{(1)}(D(f))
-l^{(1)}_{n_{p_k},s}\left(D_{n_{p_k}}
\left(f_{n_{p_k}}\right)\right)\right)\right\|_{F^{(1)}_s}\\
&=&\left\|l_s^{(1)}(D(f))\right\|_{F^{(1)}_s}
\le \|D(f)\|_{F^{(1)}}
\le \C\|f\|_{F^{(2)}}.
\ea
Hence \eqref{E2.47} holds true in this case as well
with $C_5=1/\C$, by \eqref{E2.18} of Condition C6b).

Further, we use  \eqref{E2.42}, \eqref{E2.46}, \eqref{E2.44},  \eqref{E2.45},
and Conditions C5) and C11a)   to obtain
the following relations for $s\ge s_0$:
\ba
\lim_{k\to\iy} \C_{n_{p_k}}
&\le &\frac{\limsup_{k\to\iy}\left\|D_{n_{p_k}}\left(f_{n_{p_k}}\right)
\right\|_{F_{n_{p_k}}^{(1)}}}
{\liminf_{k\to\iy}\left\|f_{n_{p_k}}\right\|_{F_{n_{p_k}}^{(2)}}}+\vep\\
&\le & \frac{\limsup_{k\to\iy}\left\|l_{n_{p_k},s}^{(1)}
\left(D_{n_{p_k}}\left(f_{n_{p_k}}\right)\right)
\right\|_{F_{s}^{(1)}}}
{\liminf_{k\to\iy}\left\|l_{n_{p_k},s}^{(2)}
\left(f_{n_{p_k}}\right)\right\|_{F_{s}^{(2)}}}+\vep\\
&=&\frac{\left\|l^{(1)}_{s}(D(f))
\right\|_{F_{s}^{(1)}}}
{\left\|l^{(2)}_{s}(f)\right\|_{F_{s}^{(2)}}}+\vep.
\ea
Therefore for some $f\in B$, satisfying \eqref{E2.47}, we have
\beq\label{E2.48}
\lim_{k\to\iy} \C_{n_{p_k}}
\le \frac{\limsup_{s\to\iy}\left\|l^{(1)}_{s}(D(f))
\right\|_{F_{s}^{(1)}}}
{\liminf_{s\to\iy}\left\|l^{(2)}_{s}(f)\right\|_{F_{s}^{(2)}}}+\vep,
\eeq
where $\liminf_{s\to\iy}\left\|l^{(2)}_{s}(f)\right\|_{F_{s}^{(2)}}>0$
by \eqref{E2.47}.
Since $\left\|f\right\|_{F^{(2)}}>0$ by \eqref{E2.47},
it follows from
 \eqref{E2.43},
\eqref{E2.48} and Conditions C6a) and C6b) that
\beq\label{E2.49}
\limsup_{n\to\iy} \C_{n}
\le \frac{\|D(f)\|_{F^{(1)}}}{\|f\|_{F^{(2)}}}+\vep
 \le \C+\vep.
 \eeq
 Finally letting $\vep\to 0+$ in \eqref{E2.49}, we arrive at \eqref{E2.34}.
 \hfill $\Box$ \vspace{.12in}

 Combining Theorems \ref{T2.10}, \ref{T2.12}, \ref{T2.13}
 and Corollary \ref{C2.11}, we can obtain the relation
 \beq\label{E2.49a}
\C=\C\left(D,B,F^{(1)},F^{(2)}\right)= \lim_{n\to\iy}
\C\left(D_n,B_n,F^{(1)}_n,F^{(2)}_n\right)=\lim_{n\to\iy}\C_n.
\eeq
\begin{corollary}\label{C2.13a}
If all the conditions of Theorems \ref{T2.10}
or \ref{T4.1aa} or Corollary \ref{C2.11}
are satisfied and, in addition, all the conditions of Theorem \ref{T2.12}
or Theorem \ref{T2.13} are satisfied, then relation \eqref{E2.49a} holds true.
\end{corollary}
Special versions of Theorems \ref{T2.12}
and \ref{T2.13} are presented below.
\begin{theorem}\label{T4.6}
Let $\C_n$ be a finite number for every $n\in\N$.\\
(a) Let
Conditions C5), C6*a), C8), C10*a), and C10*b) be satisfied.
Then \eqref{E2.34} holds true.\\
(b) Let
Conditions C5), C6*a), C8), C12*a), and C12*b) be satisfied.
Then \eqref{E2.34} holds true.
\end{theorem}
\proof
The proof of Theorem \ref{T4.6} follows that of Theorem
\ref{T2.12} (or Theorem \ref{T2.13}) if Condition C10)
(or Condition C12)) is replaced with Conditions
 C10*a) and C10*b) (or Conditions C12*a)
and C12*b)).
The outline of the proof with a new component of the proof,
compared with that of Theorem \ref{T2.12}
(or Theorem \ref{T2.13}), is presented below.

Given $\vep>0$ let  $\{f_n\}_{n=1}^\iy,\,
f_n=f_{n,\vep}\in B_n^\prime$ with $\|f_n\|_{F^{(2)}_n}=1$
be a sequence such that \eqref{E2.35}  (or \eqref{E2.42})
holds true.
In addition, let $\{n_p\}_{p=1}^\iy$ be a sequence such that
\eqref{E2.36}  (or \eqref{E2.43}) is valid.
Next, by Condition C10*a) (or C12*a), there exist a subsequence
 $\{n_{p_k}\}_{k=1}^\iy$ and $f\in B$ such that
 \eqref{E2.37}  (or \eqref{E2.44} and \eqref{E2.46}) hold(s) true.
 Further,  by Condition C10*b) (or C12*b), there exists a sequence
 of operators $\left\{l_n^{(2)}\right\}_{n=1}^\iy$
 (or, in addition, a sequence of operators
 $\left\{l_{n,s}^{(2)}\right\}_{n,s=1}^\iy$)
 such that \eqref{E2.38}  (or \eqref{E2.45}) is valid.

 The rest of the proof can be copied from the proof
 of Theorem \ref{T2.12} (or Theorem \ref{T2.13}).\hfill $\Box$

\begin{remark}\label{R2.13b}
 It is not difficult to extend asymptotic relations of this section to
 families $D_\al,\,B_\al,\,F_\al^{(1)},$ and $F_\al^{(2)}$ defined on a set $E$
 filtering with respect to a filter $\Phi$. The corresponding upper
 and lower limits are considered with respect to the filtering set $E$
 (cf. \cite[Ch. 4]{Bou1947}).
 \end{remark}
 \begin{remark}\label{R2.13c}
 We apply the results of this section to Problems A, B, and C
 of approximation theory (see Sections 2.3, 5, and 7).
 However, it is possible to apply them to Problem D of
 analysis as well (see Section 2.3). The similar condition-based approach
  to norms of linear operators in periodic and nonperiodic spaces
  was discussed in \cite{BG1983}.
  \end{remark}

\section{Asymptotic Relations in Approximation Theory}\label{S5}
\setcounter{equation}{0}
\noindent
We recall that three major Problems A, B, and C of
approximation theory related
 to sharp constants are described in Section \ref{S2.3}.
 Here, we apply Theorems  \ref{T2.10}, \ref{T4.1aa}
 \ref{T2.12}, \ref{T2.13}, and \ref{T4.6}
 and Corollary \ref{C2.11} to these problems of approximation theory.
 In Problem A, in general, we cannot essentially reduce
 the number of conditions
   in asymptotic results of  Section \ref{S4}.
  However, we can do that in Problems B and C.
  \subsection{Inequalities for Approximating Elements (Problem A)}\label{S5.1}
  In this case we assume that
  $F^{(j)}$ and $F^{(j)}_n$ are  $\beta_j$-seminormed spaces,
  where $\be_j$ is independent of $n,\,j=1,2,\,n\in\N$.
  In addition, we set
  $\|\cdot\|_{F^{(j)}_n}^*=\|\cdot\|_{F^{(j)}_n},\,j=1,2,\,n\in\N$.
  Next, let $B$
  and $B_n$ be
nontrivial subspaces of $F^{(2)}$ and $F^{(2)}_n$, respectively,
 and let $D_n$ be a $1$-homogeneous operator (that is,
 $D_n(\al f)=\al D_n(f),\,\al>0$), $n\in\N$.

 Then Conditions C1), C1*), and C5) are satisfied by
 Proposition \ref{P2.8}.
 In addition, Conditions C8) and C9) are satisfied
  by
 Proposition \ref{P2.9} (i).
 However, it appears that in Problem A Conditions C1) through C3)
 are more efficient than Conditions C1*) through C4*)
 and Condition C9)
 is more efficient than Condition C8)
 (see Corollary \ref{C7.1} and Examples \ref{Ex7.2},
 \ref{Ex7.3}, and \ref{Ex7.4}).

 The following corollary follows from  Theorems
 \ref{T2.10}, \ref{T4.1aa},
 \ref{T2.12}, and \ref{T2.13}.
 \begin{corollary}\label{C2.14a}
 (i) If Conditions C2), C3), C4) (or C4\textprime))
 are satisfied, then
 \beq\label{E2.54a}
 \sup_{f\in B,\,\|f\|_{F^{(2)}}> 0}
\frac{\|D(f)\|_{F^{(1)}}}{\|f\|_{F^{(2)}}}
\le\liminf_{n\to\iy}
 \sup_{f\in B_n,\,\|f\|_{F^{(2)}_n}> 0}
\frac{\|D_n(f)\|_{F^{(1)}_n}}{\|f\|_{F^{(2)}_n}}.
 \eeq
 (ii) Let $\C_n$
  be a finite number for every $n\in\N$, and let
Conditions C6) and C10) be satisfied. Then
 \beq\label{E2.55a}
 \sup_{f\in B,\,\|f\|_{F^{(2)}}> 0}
\frac{\|D(f)\|_{F^{(1)}}}{\|f\|_{F^{(2)}}}
 \ge \limsup_{n\to\iy}
 \sup_{f\in B_n,\,\|f\|_{F^{(2)}_n}> 0}
\frac{\|D_n(f)\|_{F^{(1)}_n}}{\|f\|_{F^{(2)}_n}}.
\eeq
(iii) Let $\C_n$
be a finite number for every $n\in\N$, and let
Conditions C6), C11), and C12) be satisfied.
 Then \eqref{E2.55a} holds
true.
 \end{corollary}

 \begin{remark}\label{R5.2a}
 Unlike Problems B and C, inequalities for
 approximating elements were not included into
 the general approach to limit relations between sharp constants
  developed in \cite{G1992,G2000}.
 \end{remark}

 \subsection{Approximation of Individual Elements
  (Problem B)}\label{S5.2}
 In this case we set
 \beq\label{E2.50}
 F^{(1)}_n=F^{(1)};\,D(f):=f_0,\,f\in B;
 \,D_n(f_n):=f_0,\,f_n\in B_n;\,n\in\N,
 \eeq
 where $f_0\ne 0$ is
 a fixed element of $F^{(1)}$. In addition, we set
 \beq\label{E2.51}
 \|\cdot\|_{F^{(1)}_n}^*=\|\cdot\|_{F^{(1)}_n}=\|\cdot\|_{F^{(1)}}
 :=\|\cdot\|_{F^{(1)},\mathrm{triv}},\qquad n\in\N,
 \eeq
 where $\|\cdot\|_{F^{(1)},\mathrm{triv}}$ is defined by \eqref{E2.1}.
 We also assume that $B$ and $B_n$ are proper nonempty subsets of
 $F^{(2)}$ and $F^{(2)}_n$, respectively, $n\in\N$.
 It follows from \eqref{E2.50} and \eqref{E2.51} that
 \beq\label{E2.51a}
 \|D(f)\|_{F^{(1)}}=\left\|D_n(f_n)\right\|_{F^{(1)}_n}=1,\qquad f\in B,\quad
 f_n\in B_n,\quad n\in\N.
 \eeq
 Next, for certain $\be_2$-norms  $\|\cdot\|_{F^{(2)}}^*$ and
 $\|\cdot\|_{F^{(2)}_n}^*,\,n\in\N$, where $\be_2\in(0,1]$ is independent of $n$,
 we define the following "norms"
 \beq\label{E2.52}
 \|h\|_{F^{(2)}}:=\|\vphi-h\|_{F^{(2)}}^*,\quad h\in F^{(2)};\qquad
 \|h_n\|_{F^{(2)}_n}:=\|\vphi_n-h_n\|_{F^{(2)}_n}^*,\quad h_n\in F^{(2)}_n,
 \eeq
 where $\vphi\in F^{(2)}\setminus B$ and $\vphi_n\in F^{(2)}_n\setminus B_n,\,n\in\N$,
 are fixed elements. Then by \eqref{E2.5}, \eqref{E2.51a}, and \eqref{E2.52},
 \beq\label{E2.53}
 \C=1/E\left(\vphi,B,F^{(2)}_*\right),
 \qquad \C_n
 =1/E\left(\vphi_n,B_n,F^{(2)}_{*n}\right),
 \eeq
 where $F^{(2)}_{*}$ and $F^{(2)}_{*n}$ are vector spaces
 $F^{(2)}$ and $F^{(2)}_{n}$ equipped with the "norms" $\|\cdot\|_{F^{(2)}}^*$ and
 $\|\cdot\|_{F^{(2)}_n}^*,\,n\in\N$, respectively.
 In addition, we recall that
 \ba
 E\left(\vphi,B,F^{(2)}_*\right)
 =\inf_{f \in B}\|\vphi-f\|^*_{F^{(2)}},\quad
 E\left(\vphi_n,B_n,F^{(2)}_{*n}\right)
 =\inf_{f_n \in B_n}\|\vphi_n-f_n\|^*_{F^{(2)}_n},\quad n\in\N,
 \ea
 are the errors of best approximation (see  \eqref{E2.8*}).
 We also assume that
 $E\left(\vphi,B,F^{(2)}_*\right)>0$ and
 $E\left(\vphi_n,B_n,F^{(2)}_{*n}\right)>0$, so $\C<\iy$
 and $\C_n<\iy,\,n\in\N$.
 Further, let
 \beq\label{E2.53a}
 L_n^{(1)}=l_n^{(1)}=l_{n,s}^{(1)}=I,\,n\in\N,\,s\in\N,
 \eeq
 be the corresponding identity
 operators (because of \eqref{E2.50}
 and \eqref{E2.51}).
 Since by \eqref{E2.51a},
 $B_n\setminus\{0\}=\left\{f_n\in B_n:\|f_n\|_{F^{(2)}_n}>0,\,
\|D_n(f_n)\|_{ F^{(1)}_n}=1\right\}$,
 we choose $B_n^\prime$
 from Condition C9) as a subset of $B_n\setminus\{0\}$,
  satisfying the equality
\ba
E\left(\vphi_n,B_n,F^{(2)}_{*n}\right)
=E\left(\vphi_n,B_n^\prime,F^{(2)}_{*n}\right),\qquad n\in\N.
\ea
 Then Conditions C1), C1*), and C5) are satisfied by
 Proposition \ref{P2.8}. In addition, Condition C9) is satisfied
 by the choice of $B_n^\prime,\,n\in\N$, and by Proposition \ref{P2.9} (iii)
 if we take into account \eqref{E2.26}, \eqref{E2.51a}, and \eqref{E2.53}.
  Next, relations
 \eqref{E2.12}, \eqref{E2.13}, \eqref{E2.14}, \eqref{E2.17}, \eqref{E2.20},
 \eqref{E2.29}, and \eqref{E2.31}
   of Conditions C3), C3*), C4), C4\textprime), C6a), C10),
 and C12) are trivially satisfied by \eqref{E2.50}, \eqref{E2.51},
 \eqref{E2.51a}, and \eqref{E2.53a}.

 In particular, Conditions C1) through C3) can be reduced to
 inequality  \eqref{E2.11}  from C3),
 which coincides with \eqref{E2.54} below.
 Since \eqref{E2.54} is one of the inequalities we wish to prove,
 Conditions C1) through C3) are obviously inefficient in
 Problem B.
 However, Conditions C1*) through C4*)
 can be reduced to more efficient relations \eqref{E2.15} and \eqref{E2.16}
 from C3*) and C4*), respectively
 (see Corollary \ref{C7.6} (i) and Examples \ref{Ex7.8},
 \ref{Ex7.9}, and \ref{Ex7.10}).

 The following corollary follows from relations
 \eqref{E2.50} through \eqref{E2.53a}, Corollary \ref{C2.11}, and Theorems
 \ref{T2.12} and \ref{T2.13}.
 \begin{corollary}\label{C2.14}
 (i) If there exists a sequence of operators
$L_n^{(2)}:F^{(2)}\to F_n^{(2)},\,n\in\N$, such that
 relations \eqref{E2.15} of C3*) and \eqref{E2.16} of C4*)
 are satisfied, then
 \beq\label{E2.54}
 \limsup_{n\to\iy} E\left(\vphi_n,B_n,F^{(2)}_{*n}\right)
 \le E\left(\vphi,B,F^{(2)}_*\right).
 \eeq
 (ii) If there exists a sequence of operators
$l_n^{(2)}:F^{(2)}\to F_n^{(2)},\,n\in\N$, such that
 C6b) and relation \eqref{E2.21} of C10) are satisfied, then
 \beq\label{E2.55}
 \liminf_{n\to\iy} E\left(\vphi_n,B_n,F^{(2)}_{*n}\right)
 \ge E\left(\vphi,B,F^{(2)}_*\right).
 \eeq
 (iii) If there exist families of operators
$l_n^{(2)}:F^{(2)}\to F_n^{(2)},\,n\in\N$,
and $l_{n,s}^{(2)}:F^{(2)}_n\to F_s^{(2)},\,n\in\N,\,s\in\N$,
 such that
 C6b), C11a), and relation \eqref{E2.30} of C12) are satisfied, then
\eqref{E2.55} holds true.
\end{corollary}
\begin{remark}\label{R2.17a}
Special cases of Corollary \ref{C2.14} were discussed in \cite{G1992, G2000}.
In case of linear operators $L_n^{(2)}:F_{*}^{(2)}\to F_{*n}^{(2)}$
with
$\limsup_{n\to \iy}\left\|L_{n}^{(2)}\right\|_{F_{*}^{(2)}\to F_{*n}^{(2)}}
\le 1$ and $\vphi_n=L_n^{(2)}(\vphi),\,n\in\N$, relation \eqref{E2.16}
of C4*) is satisfied. Indeed, for every $h\in F^{(2)}$,
\beq\label{E2.56}
\limsup_{n\to \iy}\left\|L_{n}^{(2)}(h)\right\|_{F_{n}^{(2)}}
=\limsup_{n\to \iy}\left\|L_{n}^{(2)}(\vphi)-L_{n}^{(2)}(h)
\right\|_{F_{n}^{(2)}}^*
\le \|\vphi-h\|^*_{F^{(2)}}
=\|h\|_{F^{(2)}}.
\eeq
Next, in case of linear operators $l_n^{(2)}:F_{*}^{(2)}\to F_{*n}^{(2)}$
with
$\liminf_{n\to \iy}\left\|l_{n}^{(2)}(h)\right\|^*_{F_{n}^{(2)}}
\ge \|h\|^*_{F^{(2)}}$ and $\vphi_n=l_n^{(2)}(\vphi),\,n\in\N$, relation \eqref{E2.18}
of C6b) is satisfied. The proof of this fact is similar to \eqref{E2.56}.

Further, in case of linear operators $l_{n,s}^{(2)}:F_{*n}^{(2)}\to F_{*s}^{(2)}$
with
$l_s^{(2)}=l_{n,s}^{(2)}\,l_n^{(2)},\,n\in\N,\,s\in\N;\linebreak
\sup_{n\in\N,s\in\N}\left\|l_{n,s}^{(2)}\right\|_{F_{*n}^{(2)}\to F_{*s}^{(2)}}
\le 1$,
 and $\vphi_n=l_n^{(2)}(\vphi),\,n\in\N$, relation \eqref{E2.27}
of C11a) is satisfied.
Indeed, for every $h\in F^{(2)}$,
\ba\label{E2.57}
\left\|l_{n,s}^{(2)}(h)\right\|_{F_{s}^{(2)}}
=\left\|l_{s}^{(2)}(\vphi)-l_{n,s}^{(2)}(h)\right\|_{F_{s}^{(2)}}^*
=\left\|l_{n,s}^{(2)}\left(l_{n}^{(2)}(\vphi)-h\right)\right\|_{F_{s}^{(2)}}^*
\le  \|l_{n}^{(2)}(\vphi)-h\|^*_{F^{(2)}_n}
=\|h\|_{F^{(2)}_n}.
\ea
In all these special cases,
Corollary \ref{C2.14} was established in \cite[Theorems 2.1a) and 2.2a)]{G1992} for
normed  spaces $F_{*}^{(2)}$ and $F_{*n}^{(2)},\,n\in\N$. For
$\be_2$-normed  spaces
$F_{*}^{(2)}$ and $F_{*n}^{(2)},\,n\in\N$, the corresponding results were proved in
\cite[Theorems 1.1a) and 1.2a)]{G2000}.
\end{remark}

\subsection{Approximation on Classes of Elements (Problem C)}\label{S5.3}
 We set
 \beq\label{E2.58}
 F^{(2)}=F^{(1)},\quad F^{(2)}_n=F^{(1)}_n,\quad
 D=D_n=I,\quad\,n\in\N,
 \eeq
 where $I$ is the corresponding imbedding   operator,
 and
 \beq\label{E2.59}
 \|\cdot\|_{F^{(2)}}
 :=\|\cdot\|_{F^{(2)},\mathrm{triv}},\qquad
 \|\cdot\|_{F^{(2)}_n}^*=\|\cdot\|_{F^{(2)}_n}
 :=\|\cdot\|_{F^{(2)}_n,\mathrm{triv}},\quad
  n\in\N.
 \eeq
 It follows from \eqref{E2.59} that
 \beq\label{E2.59a}
 \|f\|_{F^{(2)}}=\left\|f_n\right\|_{F^{(2)}_n}=1,\qquad f\in B\setminus\{0\},\quad
 f_n\in B_n\setminus\{0\},\quad n\in\N.
 \eeq
 Next, for $\be_1$-norms $\|\cdot\|_{F^{(1)}}^*$ and
 $\|\cdot\|_{F^{(1)}_n}^*,\,n\in\N$, where $\be_1\in(0,1]$ is independent of $n$,
 we define the following "norms"
 \bna
 &&\|h\|_{F^{(1)}}:=E\left(h,G,F^{(1)}_*\right)
 =\inf_{g\in G}\|h-g\|^*_{F^{(1)}},\qquad h\in F^{(1)},\label{E2.60}\\
 &&\|h_n\|_{F^{(1)}_n}:=E\left(h_n,G_n,F^{(1)}_{*n}\right)
 =\inf_{g_n\in G_n}\|h_n-g_n\|^*_{F^{(1)}_n},\qquad h_n\in F^{(1)}_n,\label{E2.60a}
 \quad n\in\N.
 \ena
 Here, $F^{(1)}_{*}$ and $F^{(1)}_{*n}$ are vector spaces
 $F^{(1)}$ and $F^{(1)}_{n}$ equipped with the "norms" $\|\cdot\|_{F^{(1)}}^*$ and
 $\|\cdot\|_{F^{(1)}_n}^*$, respectively, and, in addition,
 $G$ and $G_n$ are subspaces of $F^{(1)}_{*}$ and $F^{(1)}_{*n}$,
 respectively, \,$n\in\N$.
 Note that
 $\|\cdot\|_{F^{(1)}}$ and $\|\cdot\|_{F^{(1)}_n},\,n\in\N$,
 defined by \eqref{E2.60} and \eqref{E2.60a}
  are $\be_1$-seminorms.
  Then by \eqref{E2.5}, \eqref{E2.59a}, \eqref{E2.60}, and \eqref{E2.60a},
 \beq\label{E2.61}
 \C=\sup_{f\in B\setminus\{0\}}E\left(f,G,F^{(1)}_*\right),\quad
  \C_n=\sup_{f_n\in B_n\setminus\{0\}}E\left(f_n,G_n,F^{(1)}_{*n}\right)
 =\sup_{f_n\in B_n}E\left(f_n,G_n,F^{(1)}_{*n}\right)
 \eeq
 are sharp constants of best approximation on classes of elements
 $B\setminus\{0\}$ and $B_n,\,n\in\N$.
 In addition,
 since by \eqref{E2.59a},
 $B_n\setminus\{0\}=\left\{f_n\in B_n:
\|f_n\|_{ F^{(2)}_n}=1\right\}$, we
 choose $B_n^\prime$ from Condition C8) as a subset of $B_n\setminus\{0\}$,
 satisfying the equality
 \ba
 \sup_{f_n\in B_n}E\left(f_n,G_n,F^{(1)}_{*n}\right)
 =\sup_{f_n\in B_n^\prime}E\left(f_n,G_n,F^{(1)}_{*n}\right), n\in\N.
 \ea
Let us also define $l_n^{(2)}:F^{(2)}\to F_n^{(2)}$
and $l_{n,s}^{(2)}:F^{(2)}_n\to F_s^{(2)}$
 for each sequence $\{f_n\}_{n=1}^\iy,\,
f_n\in B_n\setminus\{0\},\,n\in\N,\,s\in\N$, and $f\in B$
 (mentioned in C10*a)
and C12*a)), by the formulae
\beq\label{E2.63}
l_n^{(2)}(h):=\left\{\begin{array}{ll}f_n,&h=f,\\
f_{n,0},&h\in F^{(2)}\setminus \{f\}, \end{array}\right.\quad
l_{n,s}^{(2)}(h_n):=\left\{\begin{array}{ll}f_s,&h_n=f_n,\\
f_{s,0},&h_n\in F^{(2)}_n\setminus \{f_n\}, \end{array}\right.
\eeq
where $f_{n,0}\ne 0$ is
a fixed element from $F_n^{(2)},\,n\in\N$.
It follows from \eqref{E2.59} and \eqref{E2.63} that
for  every $h\in F^{(2)}\setminus \{f\}$ and every
$h_n\in F^{(2)}_n\setminus \{f_n\}$,
\beq\label{E2.63a}
\left\|l_n^{(2)}(h)\right\|_{F_n^{(2)}}
=\left\|l_{n,s}^{(2)}(h_n)\right
\|_{F_s^{(2)}}
=1;\qquad l_n^{(2)}(f)=f_n,\quad l_{n,s}^{(2)}(f_n)=f_s,\quad
n\in\N,\quad s\in\N.
\eeq
Then Conditions C1), C1*), C5) are satisfied by
 Proposition \ref{P2.8}. In addition, Condition  C8)
 is satisfied by the choice of $B_n^\prime,\,n\in\N$, and
 by Proposition \ref{P2.9} (ii) if we take into account \eqref{E2.23},
 \eqref{E2.59a}, and \eqref{E2.61}.
 Next, relations
 \eqref{E2.11}, \eqref{E2.15},
 \eqref{E2.16}, \eqref{E2.18}, \eqref{E2.24} (with $C_1=1$),
 \eqref{E2.21},
 \eqref{E2.27}, \eqref{E2.30}, \eqref{E2.24a} (with $C_2=1$),
  \eqref{E3.18a}, and \eqref{E3.21a}
   of Conditions C3), C3*), C4\textprime), C6b), C7),
    C10), C13),   C10*b), and C12*b)
  are trivially satisfied by \eqref{E2.59a} and \eqref{E2.63a}.

  Since inequality \eqref{E2.11} of Condition C3) is trivially satisfied,
  it appears that in Problem C Conditions C1) through C3)
 are more efficient than Conditions C1*) through C4*)
 (see Corollary \ref{C7.12} (i)  and Examples \ref{Ex7.15}
 and \ref{Ex7.16}).

 The following corollary follows from relations
  \eqref{E2.58} through \eqref{E2.63a} and Theorems \ref{T4.1aa},
  \ref{T4.6} (a), and \ref{T4.6} (b).
  \begin{corollary}\label{C2.19}
 (i) If relations \eqref{E2.12} of C3) and \eqref{E2.13} of C4\textprime)
 are satisfied, then
 \beq\label{E2.64}
 \sup_{f\in B\setminus\{0\}} E\left(f,G,F^{(1)}_{*}\right)
 \le \liminf_{n\to\iy}\sup_{f_n\in B_n}E\left(f_n,G_n,F^{(1)}_{*n}\right).
 \eeq
 (ii) Let $\sup_{f_n\in B_n}E\left(f_n,G_n,F^{(1)}_{*n}\right)$ be a
 finite number for every $n\in \N$.
 If Conditions C6*a) and C10*a) are satisfied, then
 \beq\label{E2.65}
 \sup_{f\in B\setminus\{0\}} E\left(f,G,F^{(1)}_{*}\right)
 \ge \limsup_{n\to\iy}\sup_{f_n\in B_n}E\left(f_n,G_n,F^{(1)}_{*n}\right).
 \eeq
 (iii)  Let $\sup_{f_n\in B_n}E\left(f_n,G_n,F^{(1)}_{*n}\right)$ be a
 finite number for every $n\in \N$.
 If Conditions C6*a) and C12*a) are satisfied, then \eqref{E2.65} holds true.
 \end{corollary}
 \begin{remark}\label{R2.19a}
 Note that relations
 $B\subseteq F^{(1)}\cap F^{(2)}$ and
$D(B)\subseteq B$ of Condition C4\textprime) are satisfied because
 $F^{(2)}=F^{(1)}$ and $D=I$ is an imbedding operator by \eqref{E2.58}.
 \end{remark}
 \begin{remark}\label{R2.20}
 For linear operators $L_n^{(1)},\,l_n^{(1)}$ and
linear operators $l_{n,s}^{(1)}$
with
$l_s^{(1)}=l_{n,s}^{(1)}\,l_n^{(1)},\,n\in\N,\,s\in\N,$ and
$\sup_{n\in\N,s\in\N}\left\|l_{n,s}^{(1)}
\right\|_{F_{*n}^{(1)}\to F_{*s}^{(1)}}
\le 1$,
 Corollary \ref{C2.19} with slightly different conditions was obtained
 in \cite[Theorems 2.1 (b) and 2.2 (b)]{G1992} for normed
 spaces $F_*^{(1)}$ and $F_{*n}^{(1)},\,n\in\N$.
 For $\be_2$-normed  spaces
  $F_*^{(1)}$ and $F_{*n}^{(1)},\,n\in\N$,
   the corresponding results were proved in
 \cite[Theorems 1.1 (b) and 1.2 (b)]{G2000}.
 \end{remark}
 \begin{remark}\label{R5.7}
 Summarizing the applications of the results from Section \ref{S4}
 to Problems A, B, and C presented in this section,
 we conclude that Corollary \ref{C2.14a} is new, while weaker versions
 of Corollaries \ref{C2.14} and \ref{C2.19} were obtained in \cite{G1992, G2000}
 (see also Remarks \ref{R5.2a}, \ref{R2.17a} and \ref{R2.20}).
 \end{remark}

 \section{Function Spaces, Function Classes, and Operators}\label{S6}
 \setcounter{equation}{0}
\noindent
Applications of Corollaries \ref{C2.14a}, \ref{C2.14}, and \ref{C2.19} to
concrete problems of approximation theory (see Section \ref{S7})
 require more specified spaces.
 Here, we introduce some function spaces, function classes,
 and operators and study their certain properties.
 \setcounter{equation}{0}
\noindent
 \subsection{Notation}\label{S6.1}
Let $\R^m$ be the Euclidean $m$-dimensional space with elements
$x=(x_1,\ldots,x_m),\, y=(y_1,\ldots,y_m),
\,t=(t_1,\ldots,t_m)$,
the inner product $t\cdot y:=\sum_{j=1}^mt_jy_j$,
and the norm $\vert t\vert:=\sqrt{t\cdot t}$.
Next, $\CC^m:=\R^m+i\R^m$ is the $m$-dimensional
complex space with elements
$z=(z_1,\ldots, z_m)=x+iy$
and the norm $\vert z\vert:=\sqrt{\vert x\vert^2+\vert y\vert^2}$;
$\Z^m$ denotes the set of all integral lattice points
$k=(k_1,\ldots,k_m)$ in $\R^m$;
and $\Z^m_+$ is a subset of $\Z^m$
of all points with nonnegative coordinates.
We also use multi-indices $\be=(\be_1,\ldots,\be_m)\in \Z^m_+$
and $\al=(\al_1,\ldots,\al_m)\in \Z^m_+$
with
 \ba
 \vert\be\vert:=\sum_{j=1}^m\be_j,\qquad
 \vert\al\vert:=\sum_{j=1}^m\al_j,\qquad
 y^\be:=y_1^{\be_1}\cdot\cdot\cdot y_m^{\be_m}, \qquad
 D^\al:=\frac{\partial^{\al_1}}{\partial y_1^{\al_1}}\cdot\cdot\cdot
 \frac{\partial^{\al_m}}{\partial y_m^{\al_m}}.
 \ea

 Let $V$ be a centrally symmetric (with respect to the origin)
 closed
 convex body in $\R^m$ and
 $V^*:=\{y\in\R^m: \forall\, t\in V, \vert t\cdot y\vert \le 1\}$
 be the polar of $V$.
 It is well known that $V^*$ is a centrally symmetric (with respect to the origin)
 closed
 convex body in $\R^m$ and $V^{**} =V$ (see, e.g., \cite[Sect. 14]{R1970}).
 The set $V$ generates the following dual norms
 on $\R^m$ and $\CC^m$ by
 \ba
 \|y\|_V^*:=\sup_{t\in V}\vert t\cdot y\vert,\quad y\in\R^m;\qquad
 \|z\|_V^*:=\sup_{t\in V}\left\vert\sum_{j=1}^m t_jz_j\right\vert,\quad z\in\CC^m.
 \ea
 Note also that $V^*$ is the unit ball in the norm $\|\cdot\|_V^*$ on $\R^m$.
 Examples of sets $V$
 and $V^*$ and dual norms $\|\cdot\|_V^*$ are given below.
 Given $\sa\in\R^m,\,\sa_j>0,\,1\le j\le m$, and $M>0$, let
$\Pi^m_\sa:=\{t\in\R^m: \vert t_j\vert\le\sa_j, 1\le j\le m\},\,
Q^m_M:=\{t\in\R^m: \vert t_j\vert\le M, 1\le j\le m\}$,
and
$\BB^m_M:=\{t\in\R^m: \vert t\vert\le M\}$
be the $m$-dimensional parallelepiped, cube, and ball, respectively.

 For example, if $\sa\in\R^m,\,\sa_j>0,\,1\le j\le m$, and $V=\left\{t\in\R^m:
 \left(\sum_{j=1}^m\vert t_j/\sa_j\vert^{\mu}\right)^{1/\mu}\le 1\right\}$, then
 for $y\in\R^m,\,
 \|y\|_V^*=\left(\sum_{j=1}^m\vert \sa_j y_j\vert^\la\right)^{1/\la}$
 and $V^*=\left\{y\in\R^m:
 \left(\sum_{j=1}^m\vert \sa_j y_j\vert^\la\right)^{1/\la}\le 1\right\}$,
 where $\mu\in[1,\iy],\,\la\in[1,\iy]$, and $1/\mu+1/\la=1$.
 In particular,
 $\|y\|_{\Pi^m_\sa}^*=\sum_{j=1}^m\sa_j\vert y_j\vert,\,
 \|y\|_{Q^m_M}^*=M\sum_{j=1}^m\vert y_j\vert$,
 and $\|y\|_{\BB^m_M}^*=M\vert y\vert$.

  We also use the floor function $\lfloor \cdot\rfloor$.
 \noindent
\subsection{Function spaces}\label{S6.2}
In the capacity of spaces $F^{(1)},\,F^{(2)},\,F^{(1)}_n,$ and
$F^{(2)}_n$ from Sections  \ref{S3}, \ref{S4}, and \ref{S5}
 or their $\ast$-modifications from Section  \ref{S5},
 we mostly use the following function spaces.

Let $\mu$ be a $\sa$-finite positive measure on a $\sa$-algebra
 $\Sigma$ over $\R^m$
and let $S_\mu(\R^m)$ be a vector space of all
$\mu$-measurable functions
 $f:\R^m\to \CC^1$.
 Next, let $F(\R^m)=F(\R^m,\mu)$ be a $\be$-seminormed space,
 $\be\in(0,1]$,
 of functions from $S_\mu (\R^m)$.
 We assume that the  $\be$-seminorm of $F(\R^m)$ possesses
 properties of monotonicity and $\boldsymbol{\om}$-continuity.

 \begin{definition}\label{D6.1}
 We say that the  $\be$-seminorm $\|\cdot\|_{F(\R^m)}$ is
 \emph{monotone} if relations
 $g\in S_\mu (\R^m),\, f\in F(\R^m)$,
  and $\vert g(x)\vert\le \vert f(x)\vert,\,x\in\R^m$, imply
  $g\in F(\R^m)$ and
  $\|g\|_{F(\R^m)}\le \|f\|_{F(\R^m)}$.
 \end{definition}
 Given a $\mu$-measurable set $\Ome\subset\R^m$ and a space
 $F(\R^m)$,
 we define a space $F(\Ome)$ of the restrictions
 $f\left|_\Ome\right.$
 of functions $f\in S_\mu (\R^m)$ to $\Ome$ with the finite
 $\be$-seminorm  $\|f\left|_\Ome\right.\|_{F(\Ome)}
 :=\|f\cdot \chi(\Ome)\|_{F(\R^m)}$,
 where $\chi(\Ome)$ is the characteristic function of $\Ome$.
 If there is no confusion, we use the notation $f$
 instead of $f\left|_\Ome\right.$.
 For a monotone $\be$-seminorm $\|\cdot\|_{F(\R^m)}$,
 \beq\label{E6.1}
 \|f\|_{F(\Ome)}\le \|f\|_{F(\R^m)},\qquad f\in F(\R^m),
 \eeq
 so if
$f\in F(\R^m)$, then $f=f\left|_\Ome\right.\in F(\Ome)$.
In a sense, the next definition allows us to reverse
the inequality sign in \eqref{E6.1}.
\begin{definition}\label{D6.2}
 We say that the  monotone $\be$-seminorm $\|\cdot\|_{F(\R^m)}$ is
 $\boldsymbol{\om}$-\emph{continuous} if for any increasing sequence of $\mu$-measurable
 sets $\boldsymbol{\om}:=\{\Ome_n\}_{n=1}^\iy$
 in $\R^m$ with $\cup_{n=1}^\iy\Ome_n=\R^m$,
 the following relation is valid:
 \beq\label{E6.2}
 \lim_{n\to\iy}\|f\|_{F(\Ome_n)}
 = \|f\|_{F(\R^m)},\qquad f\in F(\R^m).
 \eeq
 \end{definition}
 A typical example  of a space
  with a monotone and
 $\boldsymbol{\om}$-continuous seminorm is
  the seminormed space $C_0(\R^m)$ of all complex-valued
  continuous functions $f$ on $\R^m$ with the "point" seminorm
  $\|f\|_{C_0(\R^m)}:=\vert f(0)\vert$.
  For a compact set $\Ome\subset \R^m,\,0\in\Ome,$
  we can also define the seminormed space $C_0(\Ome)$.
  We use these seminorms in Examples \ref{Ex7.3} and \ref{Ex7.4}.

 All spaces but $C_0(\Omega)$ that we discuss in Chapter \ref{S7}
 are normed or $\be$-normed.
 Examples of spaces with monotone and
 $\boldsymbol{\om}$-continuous norms include certain
 Orlicz, Lorentz, and Marcinkiewicz spaces (see, e.g.
  \cite[Sect. 4.3]{KA1982}). In particular, normed
separable spaces $F(\R^m)$ or normed spaces,
  satisfying the Fatou condition (named C-condition in
 \cite[Sect. 4.3]{KA1982}), have $\boldsymbol{\om}$-continuous norms.

 A typical example of a space with a monotone and
 $\boldsymbol{\om}$-continuous $\be$-norm is
  the space $L_{r,\mu}(\R^m), r\in(0,\iy]$,
   of all $\mu$-measurable complex-valued functions $f$
 on $\R^m$  with the finite $\be$-norm
 \beq\label{E6.2a}
 \|f\|_{L_{r,\mu}(\R^m)}:=\left\{\begin{array}{ll}
 \left(\int_{\R^m}\vert f(x)\vert^r d\mu(x)\right)^{1/r}, & 0<r<\iy,\\
 \esssup_{x\in \R^m} \vert f(x)\vert, &r=\iy,
 \end{array}\right.
 \eeq
 where $\be=\min\{1,r\}$.
 For a $\mu$-measurable set $\Ome\subset\R^m$, we also use
 the $\be$-normed space $L_{r,\mu}(\Ome)$ ($r\in(0,\iy],
 \,\be=\min\{1,r\}$).

 In all examples of Chapter \ref{S7},
 $\mu$ is a weighted measure on $\R^m$, i.e.,
 $d\mu(x)=W(x)dx,\,x\in\R^m$, where $W$ is a locally Lebesgue-integrable
 weight on $\R^m$.
 In particular, we use the weighted measure with a univariate
 power weight
 $W(x)=\vert x\vert^{2\nu+1}$ in Example \ref{Ex7.4}.
 In case of the Lebesgue measure $\mu$ (i.e., $W=1$), we use the notation
 $L_{r}(\R^m)$ and $L_{r}(\Ome),\,r\in(0,\iy]$, in Examples
 \ref{Ex7.2}, \ref{Ex7.3}, \ref{Ex7.8}, \ref{Ex7.9},
 \ref{Ex7.15}, and \ref{Ex7.16}.

 In Example \ref{Ex7.10} we use a univariate exponential weight
 $W(x)=W_r(x):=\exp(-rQ(x)),\,r\in(0,\iy)$,
 where $Q:\R^1\to [0,\iy)$ is an even differentiable
    convex function (see, e.g., \cite{LL2001, G2008}
    for applications of $W_r$ in analysis).
 A typical example of $W_r$ is
    a Freud weight $W_{r,\al}(x):=\exp(-r\vert x\vert^\al)),\,\al>1$.
 The corresponding $\be$-norm is a version of \eqref{E6.2a}
 in the following form:
 \beq\label{E6.2c}
 \|f\|_{L_{r}(W_r,\R^1)}:=
 \left(\int_{\R^1}\vert f(x)\exp(-Q(x))\vert^r dx\right)^{1/r}, \qquad 0<r<\iy,
 \eeq
 where  $\be=\min\{1,r\}$.
 An extension of \eqref{E6.2c} to $r=\iy$ is given by
 \beq\label{E6.2d}
 \|f\|_{L_{\iy}(W_\iy,\R^1)}
 :=\esssup_{x\in R^1} \vert f(x)\exp(-Q(x))\vert.
 \eeq
  The space $L_{r}(W_r,\R^1)$ ($r\in(0,\iy],\,\be=\min\{1,r\}$) is the space
   of all Lebesgue-measurable complex-valued functions $f$
 on $\R^1$  with the finite
 monotone and
 $\boldsymbol{\om}$-continuous $\be$-norms \eqref{E6.2c} and \eqref{E6.2d}.
 Note that norm \eqref{E6.2d} is different, compared with
 norm \eqref{E6.2a} for $r=\iy$.

   Given $a>0$, the space of all functions $f$ from $L_\iy(\R^m)$
 that are $2\pi a$-periodic in each variable is denoted by
 $\tilde{L}_{\iy,2\pi a,m}$ with the norm
 $\|f\|_{\tilde{L}_{\iy,2\pi a,m}}=\|f\|_{L_\iy(R^m)}
 =\|f\|_{L_\iy(Q^m_{\pi a})}$ (see Example \ref{Ex7.15}).

 \subsection{Function classes}\label{S6.3}
 In the capacity of sets $B$ and
$B_n$ from Sections  \ref{S3}, \ref{S4}, and \ref{S5},
 we use the following function classes.

  Let  $\PP_{n,m}$ be a set of all
  polynomials $P(y)=\sum_{\vert\be\vert\le n}c_\be y^\be$
   in $m$ variables of degree at most $n,\,n\in \Z^1_+$,
   with complex coefficients.
Given $a>0$, the set of all trigonometric polynomials
 $Q(x)=\sum_{k\in aV\cap \Z^m}c_ke^{(ik)\cdot x}$ with complex
 coefficients is denoted by $\TT_{aV}$
 and the set of all trigonometric polynomials
 $T(x)=\sum_{k\in aV\cap \Z^m}c_ke^{(ik)\cdot (x/a)}$ with complex
 coefficients is denoted by $\TT_{a,V}$.

  \begin{definition}\label{D6.3}
 We say that an entire function $f:\CC^m\to \CC^1$ has exponential type $V$
 if for any $\vep>0$ there exists a constant $C_0(\vep,f)$ such that
 for all $z\in \CC$,
 $\vert f(z)\vert\le C_0(\vep,f)\exp\left((1+\vep)\|z\|_V^*\right)$.
 \end{definition}
  The set of all entire function of exponential type $V$ is denoted
  by $B_V$.
The set $B_V$ was defined by Stein and Weiss
  \cite[Sect. 3.4]{SW1971}.
  For $V=\Pi^m_\sa,\,V=Q^m_M,$ and $V=\BB^m_M$,
   similar
  sets of entire functions were
  defined by Bernstein \cite{B1948} and Nikolskii
  \cite{N1951}, \cite[Sects. 3.1, 3.2.6]{N1969}, see also
  \cite[Definition 5.1]{DP2010}.
  Properties of functions from $B_V$ have been investigated in numerous
  publications (see, e.g., \cite{B1948, N1951, N1969, SW1971, NW1978,
  G1982, G1991, G2001} and
  references therein).

  Throughout the remaining sections of the paper,
  if no confusion may occur, the same notation
   is applied to functions
  $f\in B_V,\,f\in\TT_{aV},\,f\in\TT_{a,V}$, and $f\in\PP_{n,m}$
   and their restrictions to $\R^m$
   or to a subset $\Ome$ of $\R^m$ (e.g., in the form
  $f\in  B_V\cap L_p(\R^m))$.
  In addition, $\TT_{aV}\subset B_{aV}$ and
  $\TT_{a,V}\subset B_{V}$.

  Let $M_{C,N}$ be the set of all Lebesgue measurable functions $f$,
  satisfying the inequality
  \ba
  \vert f(x)\vert\le C(1+\vert x\vert)^N
  \ea
  for a. a. $x\in\R^m$,
  where $C>0$ and $N\ge 0$ are constants independent of $x$.
  In particular, we discuss classes $\TT_{a,V}\cap M_{C,N}$
   in Examples \ref{Ex7.8} and \ref{Ex7.15}.

  In addition, let $\om:[0,\iy)\to [0,\iy)$ be a modulus of continuity,
  i.e., $\om(0)=0$ and $\om$ is a continuous, nondecreasing, and
  subadditive function on $[0,\iy)$.
  Next, we define  the class $H_\om(\Omega)$
   of all functions $f$ on  $\Omega\subseteq \R^m$,
  satisfying the inequality
  \ba
  \vert f(x)-f(y)\vert\le \om(\vert x-y\vert),
  \qquad x\in\Omega,\quad y\in\Omega.
  \ea
  We also need a subclass
  $H_{\omega,0}(\R^m):=\{f\in H_\om(\R^m):f(0)=0\}$ of $H_\om(\R^m)$.
  Note that
  \beq\label{E6.5aa}
  H_{\omega,0}(\R^m)\subseteq M_{\om(2),1}.
  \eeq
  Indeed, if $f\in H_{\omega,0}(\R^m)$, then
  $\vert f(x)\vert\le\om(\vert x\vert),\,x\in\R^m$. Next,
  $\max_{\tau\in[0,2]}\om(\tau)=\om(2)$, and
  if $2^k<\tau\le 2^{k+1},\,k\in\N$,
  then $\om(\tau)\le \om(2^{k+1})< \tau\om(2)$. Therefore,
  $\vert f(x)\vert\le \om(2)(1+\vert x\vert),\, x\in\R^m$,
  and \eqref{E6.5aa} holds true.

  In case of $\om(\tau)=\tau^\la,\,\tau\in[0,\iy),\, 0<\la\le 1$, the class
  $H^\la(\Omega):=H_\om(\Omega)$ is called the H\"{o}lder class.
  For $m=1$ we also define the class $W^sH_\om(\R^1)$
  of all $s$-differentiable functions on $\R^1$ with
  $f^{(s)}\in H_\om(\R^1),\,s\in \Z^1_+$.

  We also recall that the error of best approximation of
  $\vphi\in F(\Ome)$ by elements from a subspace $B$ of $F(\Ome)$
  is defined by
  \beq\label{E6.2a1}
  E(\vphi,B,F(\Ome)):=\inf_{f\in B}\|\vphi-f\|_{F(\Ome)}
  \eeq
  (see \eqref{E2.8*}).

  Certain properties of functions from defined classes
   are discussed below. We first prove multivariate versions
    of two known univariate results
    (see, e.g., \cite[Sects. 4.10.2 and 2.6.21]{T1963}).

    \begin{lemma}\label{L6.3a}
    If a function $f\in B_V$ has a period of $2\pi a$ in each variable,
    then $f\in\TT_{a,V},\,a>0$.
    \end{lemma}
    \proof
    The entire function $f$ coincides with its Fourier series
    $\sum_{k\in \Z^m}c_ke^{(ik)\cdot (x/a)}$.
    Since $f\in L_\iy(\R^m)$,
    this function is a tempered distribution on $\R^m$.
    Then by the generalized Paley-Wiener theorem \cite[p. 13]{NW1978},
    the support of the Fourier transform of $f$ is a subset of $V$.
    So $c_k=0$ for $k\in\Z^m\setminus aV$. \hfill$\Box$

     \begin{lemma}\label{L6.3b}
     If $\vphi\in \tilde{L}_{\iy,2\pi a,m}$, then
     \beq\label{E6.2a2}
     E(\vphi,B_V,L_\iy(\R^m))
     =E(\vphi,\TT_{a,V},L_\iy(\R^m)).
     \eeq
     \end{lemma}
     \proof
     Since $\TT_{a,V}\subseteq B_V$, we have
     \beq\label{E6.2a3}
     E(\vphi,B_V,L_\iy(\R^m))
     \le E(\vphi,\TT_{a,V},L_\iy(\R^m)).
     \eeq
     Next given $\vep\in (0,\iy)$, let $g_\vep\in B_V\cap L_\iy(\R^m)$
     satisfy the relation
     \beq\label{E6.2a4}
     \|\vphi-g_\vep\|_{L_\iy(\R^m)}
     <E(\vphi,B_V,L_\iy(\R^m))+\vep.
     \eeq
     Then the function
     \ba
     G_N(x):=(2N+1)^{-m}
     \sum_{k\in Q_N^m}g_\vep(x+2\pi a k),\qquad N\in\N,
     \ea
     belongs to $B_V$ and by \eqref{E6.2a4},
     \beq\label{E6.2a5}
     \|\vphi-G_N\|_{L_\iy(\R^m)}
     <E(\vphi,B_V,L_\iy(\R^m))+\vep.
     \eeq
     Moreover,
     $\sup_{N\in\N}\|G_N\|_{L_\iy(\R^m)}\le \|g_\vep\|_{L_\iy(\R^m)}$.
      Applying the compactness theorem
     (see \cite[Lemma 2.3]{G2018}) to the sequence
     $\{G_N\}_{N=1}^\iy$ we see that there exist a sequence
     of natural numbers $\{N_s\}_{s=1}^\iy$ and a function
     $G_0\in B_V$ such that
     \beq\label{E6.2a6}
     \lim_{s\to\iy} G_{N_s}(x)=G_0(x)
     \eeq
     uniformly on any compact $K\subset\R^m$. Then
     \eqref{E6.2a5} and \eqref{E6.2a6} imply that
      \beq\label{E6.2a7}
     \|\vphi-G_0\|_{L_\iy(\R^m)}
     \le E(\vphi,B_V,L_\iy(\R^m))+\vep.
     \eeq
     Next, denoting by $e_j$ the $j$th unit vector
      in $\R^m,\,1\le j\le m$,
     we have
     \beq\label{E6.2a8}
     \max_{1\le j\le m}\left\|G_N(\cdot+2\pi a e_j)
     -G_N(\cdot)\right\|_{L_\iy(R^m)}\le C_7/N,
     \eeq
     where a constant $C_7>0$ is independent of $N$. Setting $N=N_s$ and
     letting $s\to\iy$ in \eqref{E6.2a8},
     we see from \eqref{E6.2a6} and \eqref{E6.2a8}
      that $G_0$ is $2\pi a$-periodic
      in each variable.
      Therefore, $G_0\in\TT_{a,V}$ by Lemma \ref{L6.3a},
      and it follows from \eqref{E6.2a7} that
      \beq\label{E6.2a9}
     E(\vphi,\TT_{a,V},L_\iy(\R^m))
     \le \|\vphi-G_0\|_{L_\iy(\R^m)}
     \le E(\vphi,B_V, L_\iy(\R^m))+\vep.
     \eeq
     Letting $\vep\to 0+$ in \eqref{E6.2a9} and taking account of
     \eqref{E6.2a3}, we arrive at \eqref{E6.2a2}.\hfill$\Box$
     \vspace{.12in}\\
     We also need the following extension theorem.

     \begin{lemma}\label{L6.3c}
    Given function $f:Q_M^m\to\R^1,\,M>0$, there exists a function
    $\tilde{f}:\R^m\to\R^1$ such that $\tilde{f}=f$ on $Q_M^m$ and
    $\tilde{f}$ is $4M$-periodic in each variable. In addition,
    if $f\in H_\om(Q_M^m)$, then $\tilde{f}\in H_\om(\R^m)$.
    \end{lemma}
    \proof
    Let us set
    \beq\label{E6.2a10}
    h_p(x):=(-1)^{k_p}(x_p-2Mk_p),
    \quad x\in Q_M^m+2Mk, \quad 1\le p\le m,\quad k=(k_1,\ldots, k_m)\in\Z^m.
    \eeq
    Then $h(x)=(h_1(x),\ldots,h_m(x)):\R^m\to Q_M^m$
    is a vector-valued function. Indeed, this statement is obviously
    true for any interior point $x$ of $Q_M^m+2Mk,\,k\in\Z^m$.
    Let $x$ belong to the boundary of $Q_M^m+2Mk,\,k\in\Z^m$.
    Then without loss of generality we can assume that
    there exists $s\in \N,\,1\le s\le m$, such that
    \ba
    &&x_p=-M+2Mk_p=M+2M(k_p-1),\qquad 1\le p\le s,\\
    &&x_p=x_p^*+2Mk_p,\qquad \vert x_p^*\vert<M,
    \qquad s+1\le p\le m.
    \ea
    It follows from \eqref{E6.2a10} that $h(x_p)=(-1)^{k_p+1}M,
    \,1\le p\le s$, and $h(x_p)$
    for $s+1\le p\le m$ is uniquely defined by \eqref{E6.2a10}
    as well. Therefore, $h$ is a function.

    Next, $h$ is $4M$-periodic in each variable.
    Indeed, if $x\in Q_M^m+2Mk,\,k\in\Z^m$, then for every
    $l=(l_1,\ldots, l_m) \in\Z^m$,
    \ba
    h_p(x+4Ml)=(-1)^{k_p+2l_p} (x_p+4Ml_p-2M(k_p+2l_p))=h_p(x),
    \qquad 1\le p\le m.
    \ea
    In addition, $h(x)=x$ on $Q_M^m$. Then
    $\vert h(x)-h(x^*)\vert
    \le \vert x-x^*\vert,\,x\in\R^m,\,x^*\in\R^m.$
    Indeed, it suffices to show that $y_p=h_p(x_p)\in H^1(\R^1),
    \,1\le p\le m$.
    Since $y_p$ is $4M$-periodic on $\R^1$ and the curve has
    the slope of $1$ on $[-M,M]$ and  the slope of $-1$ on
    $[-2M,-M]$ and $[M,2M]$, we see that $y_p\in H^1(\R^1),
    \,1\le p\le m$.

    Finally, the function $\tilde{f}(x):=f(h(x))$ satisfies
    all conditions of the lemma.\hfill$\Box$
     \vspace{.12in}\\
     Note that for $m=1$ this result can be found in
     \cite[Sect. 3.5.71]{T1963},
     while for $m>1$ Shvedov \cite[Proof of Theorem 4]{S1982}
     presented a geometric construction of $\tilde{f}$
     with no proof provided.

  \subsection{Operators}\label{S6.4}
  In the capacity of operators $D$ and $D_n$ from Sections \ref{S3},
  \ref{S4}, and \ref{S5},
 we  use linear operators, in particular, linear
 differential operators of the form
 $d_N:=\sum_{\vert\al\vert=N}b_\al D^\al$
 with constant coefficients $b_\al\in\CC^1,\,
 \al\in\Z^m_+,\,\vert\al\vert=N,\,
N\in \Z^1_+$.
In case of $N=0$ we assume that $d_N$
is the corresponding imbedding (or identity) operator $I$.

 In the capacity of operators $L_n^{(1)},\,L_n^{(2)},\,l_n^{(1)},
 \,l_n^{(2)},\,l_{n,s}^{(1)}$ and $l_{n,s}^{(2)}$
  from the conditions  of Section \ref{S3.2}
  and also from Sections \ref{S4} and \ref{S5}, we use the
  following linear restriction operator.
  Let $\Ome_1$ and $\Ome_2$ be $\mu$-measurable subsets of $\R^m$ such that
  $\Ome_2\subseteq \Ome_1$ and let
  $R=R(\Ome_1,\Ome_2):F(\Ome_1)\to F(\Ome_2)$ be the restriction operator.
  For $\Omega\subseteq\R^m$ we also use the notation
  $D_n\left\vert_{\Omega}\right.:=D_nR(\R^m,\Omega)$.
  If no confusion may occur, we often use the notation $D_n$
  instead of $D_n\left\vert_{\Omega}\right.$.

  The following proposition is a direct consequence of \eqref{E6.1} and
  \eqref{E6.2}.
  \begin{proposition}\label{P6.4}
  (i) If $F(\R^m)$ is a space with a monotone  $\be$-seminorm, then
  $\|R(\Ome_1,\Ome_2)\|\le 1$, i.e., for any $h\in F(\Ome_1)$,
  \ba
  \|R(\Ome_1,\Ome_2)\,(h)\|_{F(\Ome_2)}\le  \|h\|_{F(\Ome_1)}.
  \ea
  (ii) Let $\boldsymbol{\om}=\{\Ome_n\}_{n=1}^\iy$ be
  any increasing sequence of $\mu$-measurable
 sets
 in $\R^m$ with $\cup_{n=1}^\iy\Ome_n=\R^m$.
  If $F(\R^m)$ is a space with a $\boldsymbol{\om}$-continuous
   $\be$-seminorm, then
  for any $h\in F(\R^m)$,
  \ba
  \lim_{n\to\iy}
  \|R(\R^m,\Ome_n)\,(h)\|_{F(\Ome_n)}=  \|h\|_{F(\R^m)}.
  \ea
  \end{proposition}
  In addition, for a set $A$ of functions
  $f: \R^m \to\CC^1$ we define by
  $A\left\vert_\Omega\right.:=\{R(\R^m,\Omega)(f):f\in A\}$
  the set of the restrictions to $\Omega\subseteq\R^m$
   of functions from $A$.
   In Section \ref{S7}, if no confusion may occur,
     we often do not use this restriction
   notation for sets
    $B_V,\,\TT_{aV},\,\TT_{a,V},\,\PP_{n,m}$,
    and their subsets.

  \section{Examples of Asymptotic Relations between Sharp
 Constants in Function Spaces}\label{S7}
\setcounter{equation}{0}
\noindent
 In this section we discuss specified conditions and examples of
 relations between sharp constants
 in function spaces separately for each major Problem A, B, and C
 of approximation theory described in Sections \ref{S2.3} and \ref{S5}.

 Throughout the section  $\{\Ome_n\}_{n=1}^\iy$ and
 $\{\Ome_n^*\}_{n=1}^\iy$ are
  increasing sequences of $\mu$-measurable
 sets
 in $\R^m$ with $\cup_{n=1}^\iy\Ome_n=\R^m$
 and $\cup_{n=1}^\iy\Ome_n^*=\R^m$.

 \subsection{Inequalities for Approximating Elements
  (Problem A)}\label{S7.1}
  Here, we discuss a special case of Corollary \ref{C2.14a} and
  three examples.\vspace{.12in}\\
 \textbf{General result.}
 We assume that all assumptions of Section \ref{S5.1} are held
 in the following special cases.
 First,
 we assume that $F^{(j)}=F^{(j)}(\R^m)
 =F^{(j)}(\R^m,\mu)$ is
 a $\be_j$-seminormed functional space
 with the monotone and $\boldsymbol{\om}$-continuous
 $\be_j$-seminorm
 (see Definitions \ref{D6.1} and \ref{D6.2}),
 $\be_j\in(0,1]$, and
 $F^{(j)}_n=F^{(j)}(\Ome_n)$
 or $F^{(j)}_n=F^{(j)}(\Ome_n^*),\, j=1,2,
 \,n\in\N $.
 Next, we define the restriction operators $L_n^{(1)}
 =R(\R^m,\Ome_n^*),\,
 l_n^{(1)}=l_n^{(2)}=R(\R^m,\Ome_n),
 \,n\in\N$, and
 \ba
 l_{n,s}^{(1)}=l_{n,s}^{(2)}
 =\left\{\begin{array}{ll}
 R(\Ome_n,\Ome_s), &n\ge s,\\
 R(\Ome_s,\Ome_n), &n< s,
 \end{array}\right.
\qquad n\in\N,\,s\in\N.
 \ea

We also assume that $B,\,B_n$, and $B_n^*$
are nontrivial subspaces of $F^{(2)}(\R^m),\,
 F^{(2)}(\Ome_n)$, and $F^{(2)}(\Ome_n^*)$,
 respectively, $n\in\N$.
 In addition, let $D:B\to F^{(1)}(\R^m),\,
 D_n:B_n\to F^{(1)}(\Ome_n)$,
 and  $D_n^*:B_n^*\to F^{(1)}(\Ome_n^*),\,n\in\N$,
  be linear operators.
 Finally, let
$B_n^\prime\subseteq\left\{f_n\in B_n\setminus\{0\}:
\|D_n(f_n)\|_{ F^{(1)}(\Ome_n)}=1\right\}$ be a set
from Condition C9) such that
the following relations hold:
\ba
\sup_{f_n\in B_n^\prime}
\frac{1}{\|f_n\|_{F^{(2)}(\Ome_n)}}
=\sup_{f_n\in B_n^\prime}
\frac{\|D_n(f_n)\|_{F^{(1)}(\Ome_n)}}{\|f_n\|_{F^{(2)}(\Ome_n)}}
=\sup_{f\in B_n\setminus\{0\}}
\frac{\|D_n(f)\|_{F^{(1)}(\Ome_n)}}{\|f\|_{F^{(2)}(\Ome_n)}},
\qquad n\in\N.
\ea
The existence of $B_n^\prime,\,n\in\N$, follows immediately
from Proposition \ref{P2.9} (i).

  Below we introduce some conditions on families
 $\left\{D,B,F^{(1)}(\R^m),F^{(2)}(\R^m)\right\},\linebreak
 \left\{D_n,B_n,F^{(1)}(\Ome_n),F^{(2)}(\Ome_n)\right\}$,
 and $\left\{D_n^*,B_n^*,F^{(1)}(\Ome_n^*),F^{(2)}
 (\Ome_n^*)\right\},\,
 n\in\N$,
 that
 are needed for the validity of relations like \eqref{E2.54a}
 and/or \eqref{E2.55a} in function spaces.\vspace{.12in}\\
 \emph{Condition C3.1).} For every $f\in B\setminus\{0\}$
 there exists a sequence $f_n^*\in B_n^*,\,n\in\N$,
such that
\beq \label{E7.1}
\limsup_{n\to\iy}\|f_n^*\|_{F^{(2)}(\Ome_n^*)}
\le \|f\|_{F^{(2)}(\R^m)}
\eeq
and
\beq \label{E7.2}
\lim_{n\to\iy}
\left\|D(f)-D_n^*(f_n^*)\right\|_{F^{(1)}(\Ome_n^*)}=0.
\eeq
\emph{Condition C3.1.1).} Let $\Ome_n^*\subseteq \Ome_n,\,n\in\N$.
 For every $f\in B\setminus\{0\}$
 there exists a sequence $f_n\in B_n,\,n\in\N$,
such that
\beq \label{E7.2a}
\limsup_{n\to\iy}\|f_n\|_{F^{(2)}(\Ome_n)}
\le \|f\|_{F^{(2)}(\R^m)}
\eeq
and
\beq \label{E7.2b}
\lim_{n\to\iy}
\left\|D(f)-D_n(f_n)\right\|_{F^{(1)}(\Ome_n^*)}=0.
\eeq
\emph{Condition C10.1).} (Compactness condition)
For any sequence $\{f_n\}_{n=1}^\iy$ such that
$f_n\in B_n^\prime,\,n\in\N$,
there exist a sequence $\{n_k\}_{k=1}^\iy$ of natural numbers and an element $f\in B$
such that
\bna
&&\lim_{k\to\iy}\left\|D_{n_k}(f_{n_k})-D(f)\right\|_{F^{(1)}
\left(\Ome_{n_k}\right)}=0,
\label{E7.3}\\
&&\lim_{k\to\iy}\left\|f_{n_k}-f\right\|_{F^{(2)}
\left(\Ome_{n_k}\right)}=0.
\label{E7.4}
\ena
\emph{Condition C12.1).} (Weak compactness condition)
For any sequence $\{f_n\}_{n=1}^\iy$
such that
$f_n\in B_n^\prime,\,n\in\N$,
there exists $f\in B$ such that
for any fixed $s\in\N$,
there exists a sequence $\{n_k\}_{k=1}^\iy$ of natural numbers
 ($n_k=n_k(s),\,k\in\N$)
such that
\bna
&&\lim_{k\to\iy}\left\|D_{n_k}(f_{n_k})-
D(f)\right\|_{F^{(1)}\left(\Ome_{s}\right)}=0,\label{E7.5}\\
&&\lim_{k\to\iy}\left\|f_{n_k}-f
\right\|_{F^{(2)}\left(\Ome_{s}\right)}=0,\label{E7.6}
\ena
and, in addition,
\beq\label{E7.7}
\left\|D_{n_k}(f_{n_k})\right\|_{F^{(1)}\left(\Ome_{n_k}\right)}
=\left\|D_{n_k}(f_{n_k})\right\|_{F^{(1)}\left(\Ome_{s}\right)},
 \quad k\in\N.
 \eeq
 The following corollary holds true.
 \begin{corollary}\label{C7.1}
 (i) If Condition C3.1)
 is satisfied, then
 \beq\label{E7.8}
 \sup_{f\in B\setminus\{0\}}
\frac{\|D(f)\|_{F^{(1)}(\R^m)}}{\|f\|_{F^{(2)}(\R^m)}}
\le\liminf_{n\to\iy}
 \sup_{f\in B_n^*\setminus\{0\}}
\frac{\|D_n^*(f)\|_{F^{(1)}(\Ome_n^*)}}{\|f\|_{F^{(2)}(\Ome_n^*)}}.
 \eeq
 (ii) If Condition C3.1.1)
 is satisfied, then
 \beq\label{E7.8a}
 \sup_{f\in B\setminus\{0\}}
\frac{\|D(f)\|_{F^{(1)}(\R^m)}}{\|f\|_{F^{(2)}(\R^m)}}
\le\liminf_{n\to\iy}
 \sup_{f\in B_n\setminus\{0\}}
\frac{\|D_n(f)\|_{F^{(1)}(\Ome_n)}}{\|f\|_{F^{(2)}(\Ome_n)}}.
 \eeq
 (iii) Let
 \ba
 \C\left(D_n,B_n,F^{(1)}(\Ome_n),F^{(2)}(\Ome_n)\right)
 =\sup_{f\in B_n\setminus\{0\}}
\frac{\|D_n(f)\|_{F^{(1)}(\Ome_n)}}{\|f\|_{F^{(2)}(\Ome_n)}}
\ea
 (see \eqref{E2.8b})
  be a finite number for every $n\in\N$, and let
Condition C10.1) be satisfied. Then
 \beq\label{E7.9}
 \sup_{f\in B\setminus\{0\}}
\frac{\|D(f)\|_{F^{(1)}(\R^m)}}{\|f\|_{F^{(2)}(\R^m)}}
 \ge \limsup_{n\to\iy}
 \sup_{f\in B_n\setminus\{0\}}
\frac{\|D_n(f)\|_{F^{(1)}(\Ome_n)}}{\|f\|_{F^{(2)}(\Ome_n)}}.
\eeq
(iv) Let $\C\left(D_n,B_n,F^{(1)}(\Ome_n),F^{(2)}(\Ome_n)\right)$
be a finite number for every $n\in\N$, and let
Condition C12.1) be satisfied.
 Then \eqref{E7.9} holds
true.\\
(v)  Let $\Ome_n^*\subseteq \Ome_n$ and let $B_n^*$ be a
set of the restrictions $(f_n)^*:=R(\R^m,\Omega_n^*)(f_n)$
 to $\Ome_n^*$
 of all elements $f_n$ from $B_n,\,n\in\N$.
 In addition, let
 $D_n^*((f_n)^*):=(D_n(f_n))^*,\,n\in\N$.
 Next,
let \linebreak
$\C\left(D_n,B_n,F^{(1)}(\Ome_n),F^{(2)}(\Ome_n)\right)$
be a finite number for every $n\in\N$, and let
Condition C3.1) and either C10.1) or C12.1)  be satisfied.
 Then
\bna\label{E7.8b}
\limsup_{n\to\iy}
 \sup_{f\in B_n\setminus\{0\}}
\frac{\|D_n(f)\|_{F^{(1)}(\Ome_n)}}{\|f\|_{F^{(2)}(\Ome_n)}}
&\le& \sup_{f\in B\setminus\{0\}}
\frac{\|D(f)\|_{F^{(1)}(\R^m)}}{\|f\|_{F^{(2)}(\R^m)}}\nonumber\\
 &\le&\liminf_{n\to\iy}
 \sup_{f\in B_n^*\setminus\{0\}}
\frac{\|D_n^*(f)\|_{F^{(1)}(\Ome_n^*)}}{\|f\|_{F^{(2)}(\Ome_n^*)}}.
\ena
(vi) Let
 $\C\left(D_n,B_n,F^{(1)}(\Ome_n),F^{(2)}(\Ome_n)\right)$
be a finite number for every $n\in\N$, and let
Condition C3.1.1) and either C10.1) or C12.1)  be satisfied.
 Then
 \beq\label{E7.10}
 \sup_{f\in B\setminus\{0\}}
\frac{\|D(f)\|_{F^{(1)}(\R^m)}}{\|f\|_{F^{(2)}(\R^m)}}
 = \lim_{n\to\iy}
 \sup_{f\in B_n\setminus\{0\}}
\frac{\|D_n(f)\|_{F^{(1)}(\Ome_n)}}{\|f\|_{F^{(2)}(\Ome_n)}}.
\eeq
 \end{corollary}
 \proof
 Note first that in this proof all elements in discussed
 Conditions C1) through C13) are replaced without further notice
 by the corresponding elements that are used in Conditions
 C3.1), C3.1.1), C10.1), C12.1) and relations \eqref{E7.8}
 through \eqref{E7.10}. For example, in Condition
  C3) $B_n$ is replaced by $B_n^*$.

 By the assumptions about the functional spaces and operators
 introduced above, Conditions C2), C4), C6), C6a),
  C6b),  C11), and C11a)
 of Section \ref{S3.2} are satisfied. Indeed,
 Conditions C2), C6), and C11) are satisfied by the direct
 definitions of the operators whose existence is required
 in these conditions.
 In addition,
 Conditions C4), C6a), C6b), and C11a) are satisfied
  by Proposition \ref{P6.4}.
  Next, Conditions C3.1), C10.1), and C12.1) of this section
  are equivalent to Conditions C3), C10), and C12) of Section
  \ref{S3.2}, respectively.
  Therefore, statements (i), (iii), and (iv) of Corollary \ref{C7.1}
   follow immediately
  from statements (i), (ii), and (iii) of
  Corollary \ref{C2.14a}, respectively.

  To prove statement (ii), we recall that
the $\be_1$-seminorm $\|\cdot\|_{F^{(1)}(\R^m)}$ is monotone and
$\boldsymbol{\om}$-continuous. So for
 any $f\in B\setminus\{0\}$, we obtain by
 Proposition \ref{P6.4} and relation \eqref{E7.2b}
 of Condition C3.1.1)
 \bna\label{E7.8c}
&&\|D(f)\|_{F^{(1)}(\R^m)}
=\lim_{n\to\iy}\|D(f)\|_{F^{(1)}(\Ome_n^*)}\nonumber\\
&&\le \liminf_{n\to\iy}\left(\|D(f)-D_n(f_n)\|
_{F^{(1)}(\Ome_n^*)}^{\be_1}
+\|D_n(f_n)\|_{F^{(1)}(\Ome_n^*)}^{\be_1}
\right)^{1/\be_1}\nonumber\\
&&=\liminf_{n\to\iy}\left\|D_n(f_n)\right\|_{F^{(1)}(\Ome_n^*)}
\le \liminf_{n\to\iy}\left\|D_n(f_n)\right\|_{F^{(1)}(\Ome_n)}\nonumber\\
&&\le \liminf_{n\to\iy}
\sup_{f\in B_n\setminus\{0\}}
\frac{\|D_n(f)\|_{F^{(1)}(\Ome_n)}}{\|f\|_{F^{(2)}(\Ome_n)}}
\limsup_{n\to\iy}\left(\|f_n\|_{F^{(2)}(\Ome_n)} \right).
\ena
Then \eqref{E7.8a} follows from \eqref{E7.8c} and \eqref{E7.2a}.
Finally, statement (v) follows from (i), (iii), and (iv),
and statement (vi) follows from (ii), (iii), and (iv).
  \hfill$\Box$\vspace{.12in}\\
  \textbf{Examples.}
  There is a number of publications
  \cite{T1965, T1993, G2005, LL2015a, LL2015b, GT2017, G2017,
   G2018, G2019a, G2019b, G2021b, G2021c}
   about special cases of Corollary \ref{C7.1};
    here, we call them Bernstein-Nikolskii type
    inequalities or inequalities of different metrics
    (see also Section \ref{S2.3} A).
    Examples \ref{Ex7.2}, \ref{Ex7.3}, and \ref{Ex7.4} below
    discuss inequalities \eqref{E7.8a},
    \eqref{E7.8b} and \eqref{E7.10}
    proved in these papers in various situations.

    \begin{example}\label{Ex7.2}
    (Multivariate Bernstein-Nikolskii type
    inequalities).
   Let
    $F^{(1)}(\R^m)=L_q(\R^m),\,
     \be_1=\min\{1,q\},\,
    F^{(2)}(\R^m)=L_p(\R^m),\,
    \be_2=\min\{1,p\},\,
    B=B_V\cap L_p(\R^m),\,
    0<p\le q\le\iy,\,
    \Ome_n=Q_{\pi n}^m,\,
    \Ome_n^*=Q_{n^\de}^m,\, \de=\de(p,m)\in(0,1),\,
    B_n=\TT_{n,V}\left\vert_{\Omega_n}\right.
    =\TT_{n,V},\,
    B_n^*=\TT_{n,V}\left\vert_{\Omega_n^*}\right.
    =\TT_{n,V},\,n\in\N.$
    Next, let
    $d_N=\sum_{\vert\al\vert=N}b_\al D^\al$,
    where
    $b_\al\in\CC^1$ ($\al\in \Z^m_+,\,
    \vert \al\vert=N$), and $N\in\Z^1_+$,
    are independent of
    $n$, and let
    $D=d_N\left\vert_{\R^m}\right.=d_N,\,
    D_n=d_N\left\vert_{\Omega_n}\right.=d_N,\,
    D_n^*=d_N\left\vert_{\Omega_n^*}\right.=d_N,
    \,n\in\N$.

    Inequality \eqref{E7.8a} for $0<p\le q\le\iy$ and
    equality \eqref{E7.10} for $0<p\le\iy$ and $q=\iy$
    were proved by the author
    \cite{G2018} in the following form:
    \bna
    && \sup_{f\in (B_{V}\cap L_p(\R^m))
    \setminus\{0\}}\frac{\|d_N(f)\|_{L_q(\R^m)}}
{\|f\|_{L_p(\R^m)}}
    \le \liminf_{n\to\iy}
\sup_{T\in\TT_{n,V}
\setminus\{0\}}\frac{\|d_N(T)\|_{L_q(Q_{\pi n}^m)}}
{\|T\|_{L_p(Q_{\pi n}^m)}}, \label{E7.10a}\\
&&\sup_{f\in (B_{V}\cap L_p(\R^m))
    \setminus\{0\}}\frac{\|d_N(f)\|_{L_\iy(\R^m)}}
{\|f\|_{L_p(\R^m)}}
    = \lim_{n\to\iy}
\sup_{T\in\TT_{n,V}
\setminus\{0\}}\frac{\|d_N(T)\|_{L_\iy(Q_{\pi n}^m)}}
{\|T\|_{L_p(Q_{\pi n}^m)}}. \label{E7.10b}
\ena
    The following asymptotic results
    (see \cite[Theorems 1.2 and 1.3]{G2018})
    are easy corollaries of \eqref{E7.10a} and
    \eqref{E7.10b} ($0<p\le q\le\iy$):
    \bna
   && \sup_{f\in (B_{V}\cap L_p(\R^m))
    \setminus\{0\}}\frac{\|d_N(f)\|_{L_q(\R^m)}}
{\|f\|_{L_p(\R^m)}}
    \le \liminf_{n\to\iy} n^{-N-m/p+m/q}
\sup_{Q\in\TT_{nV}
\setminus\{0\}}\frac{\|d_N(Q)\|_{L_q(Q_\pi^m)}}
{\|Q\|_{L_p(Q_\pi^m)}}, \label{E7.11}\\
&&\sup_{f\in (B_{V}\cap L_p(\R^m))
    \setminus\{0\}}\frac{\|d_N(f)\|_{L_\iy(\R^m)}}
{\|f\|_{L_p(\R^m)}}
    = \lim_{n\to\iy} n^{-N-m/p}
\sup_{Q\in\TT_{nV}
\setminus\{0\}}\frac{\|d_N(Q)\|_{L_\iy(Q_\pi^m)}}
{\|Q\|_{L_p(Q_\pi^m)}}. \label{E7.12}
\ena
     For $m=1$ and $V=[-\sa,\sa]$ relations
     \eqref{E7.11} and \eqref{E7.12} were
    established by the author and Tikhonov
    \cite[Theorems 1.4 and 1.5]{GT2017}.
    Similar relations for $m=1,\,V=[-\sa,\sa]$, and
    $N=0$ were obtained earlier by
    Levin and Lubinsky \cite[p. 246]{LL2015a},
    \cite[Theorem 1]{LL2015b}.
    A special case of $m=1,\,V=[-1,1],\,
    q=\iy,\,p=1$, and
    $N=0$ was discussed by Taikov
    \cite[Eq. (3)]{T1965}, \cite[Lemma 1]{T1993}
    and Gorbachev \cite[Theorem 2]{G2005}.
    More references can be found in
    \cite{GT2017, G2018}.

    Condition C3.1.1) of approximation of entire
     functions of exponential type $V$
     by trigonometric polynomials was verified
     by using univariate and multivariate Levitan's polynomials
     (see \cite[Sect. 2]{GT2017} for $m=1$ and
      \cite[Sect. 3]{G2018} for $m>1$).

      In addition, we choose
      \ba
       B_n^\prime
      &=&\Bigg{\{}f_n\in\TT_{n,V}\setminus\{0\}:
      \Bigg{(}\|d_N(f_n)\|_{L_\iy(Q_{\pi n}^m)}
      =\vert d_N(f_n)(0)\vert=1\Bigg{)}\nonumber\\
      &\wedge&
      \left.\left(\sup_{T\in\TT_{n,V}
\setminus\{0\}}\frac{\|d_N(T)\|_{L_\iy(Q_{\pi n}^m)}}
{\|T\|_{L_p(Q_{\pi n}^m)}}
=\|f_n\|_{L_p(Q_{\pi n}^m)}^{-1}\right)\right\},\qquad n\in\N.
      \ea
      Since the set $\TT_{n,V}$ consists of $2\pi n$-periodic
       trigonometric polynomials, we have
\bna\label{E7.8d}
       \sup_{T\in\TT_{n,V}
\setminus\{0\}}\frac{\|d_N(T)\|_{L_\iy(Q_{\pi n}^m)}}
{\|T\|_{L_p(Q_{\pi n}^m)}}
=\sup_{T\in\TT_{n,V}
\setminus\{0\}}\frac{\vert d_N(T)(0)\vert}
{\|T\|_{L_p(Q_{\pi n}^m)}},\qquad n\in\N.
\ena
Therefore, by \eqref{E7.8d} and by the standard compactness argument
(see the proof of Theorem 1.3 in \cite{G2018} for details),
we conclude that
$B_n^\prime\ne \emptyset$ and in addition,
\ba
\sup_{T\in\TT_{n,V}
\setminus\{0\}}\frac{\|d_N(T)\|_{L_\iy(Q_{\pi n}^m)}}
{\|T\|_{L_p(Q_{\pi n}^m)}}
=\sup_{f_n\in B_n^\prime}\frac{\|d_N(f_n)\|_{L_\iy(Q_{\pi n}^m)}}
{\|f_n\|_{L_p(Q_{\pi n}^m)}},\qquad n\in\N.
\ea
Thus Condition C9) is satisfied for the chosen $B_n^\prime,\,n\in\N$.
Next, equality \eqref{E7.7} of C12.1)
     holds true since for any $f_n\in B_n^\prime$ we have
     ($n\in\N,\,s\in\N$)
     \ba
     \left\|D_{n}(f_{n})\right\|_{F^{(1)}\left(\Ome_{n}\right)}
     =\|d_N(f_n)\|_{L_\iy(Q_{\pi n}^m)}
     =\left\vert d_N(f_n)(0)\right\vert
     =\|d_N(f_n)\|_{L_\iy(Q_{\pi s}^m)}
=\left\|D_{n}(f_{n})\right\|_{F^{(1)}\left(\Ome_{s}\right)}.
 \ea
     Further, the compactness condition of $\{B_n^\prime\}_{n=1}^\iy$,
      stated in \eqref{E7.5}  and \eqref{E7.6} of C12.1)
     (with a sequence $\{n_k\}_{k=1}^\iy$ that is independent of $s$),
     is satisfied as well.
     Indeed, let $f_n\in B_n^\prime,\,n\in\N$. Then for $p\in(0,\iy)$
     \ba
     &&\sup_{n\in\N}\|f_n\|_{L_\iy(Q_{\pi n}^m)}
     \le  \sup_{n\in\N}
     \C\left(d_0,\TT_{n,V},L_\iy(Q_{\pi n}^m),L_p(Q_{\pi n}^m)\right)
     \|f_n\|_{L_p(Q_{\pi n}^m)}\\
     &&=  \sup_{n\in\N}
     \C\left(d_0,\TT_{n,V},L_\iy(Q_{\pi n}^m),L_p(Q_{\pi n}^m)\right)
     / \C\left(d_N,\TT_{n,V},L_\iy(Q_{\pi n}^m),L_p(Q_{\pi n}^m)\right)
     =C_8<\iy,
     \ea
     by the choice of $B_n^\prime$ and
     by crude estimates of the constants in multivariate
      Bernstein-Nikolskii type inequalities
      (see \cite[Eq. (4.2)]{G2018}).
      Therefore,
      \bna\label{E7.8e}
      \cup_{n=1}^\iy B_n^\prime
      &\subseteq&
      \{T\in \cup_{n=1}^\iy\TT_{n,V}:
      \|T\|_{L_\iy(Q_{\pi n}^m)}\le C_8\}\nonumber\\
      &\subseteq&
      \{f\in B_V:
      \|f\|_{L_\iy(\R^m)}\le C_8\},
      \ena
      since $\TT_{n,V}\subseteq B_V,\,n\in\N$,
      and the compactness property of
      the right-hand side set in \eqref{E7.8e} is proved
     in \cite[Lemma 2.3]{G2018}
     (a special case is discussed in
     \cite[Theorem 3.3.6]{N1969}).
     Note also that by \cite[Eq. (4.2)]{G2018},
     \ba
     \C\left(d_N,\TT_{n,V},
     L_\iy(Q_{\pi n}^m),L_p(Q_{\pi n}^m)\right)
     <\iy,\qquad n\in\N.
     \ea

     So \eqref{E7.10a} for $0<p\le q\le\iy$ and
     \eqref{E7.10b} for $0<p\le\iy$ and $q=\iy$ follow from
     statements (i) and (vi) of
     Corollary \ref{C7.1}, respectively.
    \end{example}

    \begin{example} \label{Ex7.3}
    (Multivariate
    inequalities of different metrics).
    Let
    $F^{(1)}(\R^m)=C_0(\R^m),\,
    \be_1=1,\,
    F^{(2)}(\R^m)=L_p(\R^m),\,
    \be_2=\min\{1,p\},\,
   B=B_V\cap L_p(\R^m),\,
   0<p\le\iy,\,
   \Ome_n=nV^*,\,
   \Ome_n^*=\tau nV^*,\,\tau\in(0,1),\,
    B_n=\PP_{n,m}\left\vert_{\Omega_n}\right.=\PP_{n,m},\,
    B_n^*=\PP_{n,m}\left\vert_{\Omega_n^*}\right.=\PP_{n,m}$.
    Next, let
    $d_N=\sum_{\vert\al\vert=N}b_\al D^\al$,
    where
    $b_\al\in\CC^1$ ($\al\in \Z^m_+,\,\vert\al\vert=N$)
     and $N\in\Z^1_+$, are independent of
    $n$, and let
    $D=d_N\left\vert_{\R^m}\right.=d_N,\,
    D_n=d_N\left\vert_{\Omega_n}\right.=d_N,\,
    D_n^*=d_N\left\vert_{\Omega_n^*}\right.=d_N,\,n\in\N$.

    Inequalities \eqref{E7.8b} were proved by the author
    \cite{G2019b} in the following form:
    \bna\label{E7.12a}
    &&\limsup_{n\to\iy}
\sup_{P\in\PP_{n,m}
\setminus\{0\}}\frac{\vert d_N(P)(0)\vert}
{\|P\|_{L_p(nV^*)}}
\le \sup_{f\in (B_{V}\cap L_p(\R^m))
    \setminus\{0\}}\frac{\vert d_N(f)(0)\vert}
{\|f\|_{L_p(\R^m)}}\nonumber\\
&&= \sup_{f\in (B_{V}\cap L_p(\R^m))
    \setminus\{0\}}\frac{\|d_N(f)\|_{L_\iy(\R^m)}}
{\|f\|_{L_p(\R^m)}}
    \le \liminf_{n\to\iy}
\sup_{P\in\PP_{n,m}
\setminus\{0\}}\frac{\vert d_N(P)(0)\vert}
{\|P\|_{L_p(\tau nV^*)}}
\ena
for $\tau\in(0,1)$.
By letting $\tau\to 1-$ in \eqref{E7.12a},
 we arrive at the following asymptotic result
 for $0<p\le\iy$
    (see \cite[Theorem 1.2]{G2019b}):
\beq\label{E7.13}
\sup_{f\in (B_{V}\cap L_p(\R^m))
    \setminus\{0\}}\frac{\|d_N(f)\|_{L_\iy(\R^m)}}
{\|f\|_{L_p(\R^m)}}
    = \lim_{n\to\iy} n^{-N-m/p}
\sup_{Q\in\PP_{n,m}
\setminus\{0\}}\frac{\vert d_N(Q)(0)\vert}
{\|Q\|_{L_p(V^*)}}.
\eeq
     For $m=1,\,V=[-1,1]$, $V^*=[-1,1]$,
     and $d_N(f)=d^Nf/dx^N$,
     equality \eqref{E7.13} was established in
    \cite[Theorem 1.1]{G2017}. In case when
    $V$ and $V^*$ are the unit ball
    $\BB_1^m$ and $d_{2N}=\Delta^N$ is the polyharmonic operator,
    \eqref{E7.13} was proved in  \cite[Corollary 3.9]{G2019a}
    by reducing the multivariate problem to a univariate one
    in weighted spaces (see Example \ref{Ex7.4}).

    Condition C3.1) of approximation of entire
     functions from $B_V\cap L_p(\R^m)$
     by algebraic polynomials was verified in
     \cite[Lemma 4.4]{G1982} and \cite[Lemma 2.4]{G2019b}
     by using approximation methods developed by Bernstein
     \cite{B1946} (see also \cite[Sect. 5.4.4]{T1963}
     and \cite[Appendix, Sect. 83]{A1965}) and by the author
     \cite{G1982, G2019b}.

     In addition, we choose
      \ba
       B_n^\prime
      =\left\{f_n\in\PP_{n,m}\setminus\{0\}:
      \left(\vert d_N(f_n)(0)\vert=1\right)\wedge
      \left(\sup_{P\in\PP_{n,m}
\setminus\{0\}}\frac{\vert d_N(P)(0)\vert}
{\|P\|_{L_p(nV^*)}}
=\|f_n\|_{L_p(nV^*)}^{-1}\right)\right\},\,\, n\in\N.
      \ea
      By the standard compactness argument
(see the proof of Theorem 1.2 in \cite{G2019b} for details),
we conclude that
$B_n^\prime\ne \emptyset$ and in addition,
\ba
\sup_{P\in\PP_{n,m}
\setminus\{0\}}\frac{\vert d_N(P)(0)\vert}
{\|P\|_{L_p(nV^*)}}
=\sup_{f_n\in B_n^\prime}\frac{\vert d_N(f_n)(0)\vert}
{\|f_n\|_{L_p(nV^*)}},\qquad n\in\N.
\ea
Thus Condition C9) is satisfied for the
chosen $B_n^\prime,\,n\in\N$.
 Next, it is easy to verify that equality \eqref{E7.7} of C12.1)
     holds true in this case as well
     since for any $f_n \in B_n^\prime$
     we have
     \ba
     \left\|D_{n}(f_{n})\right\|_{F^{(1)}\left(\Ome_{n}\right)}
     =\left\|d_N(f_n)\right\|_{C_0\left(nV^*\right)}
     =\left\|d_N(f_n)\right\|_{C_0\left(sV^*\right)}
=\left\|D_{n}(f_{n})\right\|_{F^{(1)}\left(\Ome_{s}\right)},
 \quad n\in\N,\quad s\in\N.
 \ea
     Further, the compactness condition of
     $\{B_n^\prime\}_{n=1}^\iy$, stated in
     \eqref{E7.5}  and \eqref{E7.6} of C12.1)
     (with a sequence $\{n_k\}_{k=1}^\iy$
     that is independent of $s$),
     is satisfied. The proof of this fact given in
     \cite{G2019b} is based on results
     of Bernstein \cite{B1946} (for $m=1$)
     and the author \cite[Lemma 2.2]{G2019b}
     (for $m>1$).
     Note also that by a crude estimate of
     the constant in a multivariate inequality
     of different metrics
     (see \cite[Lemma 2.7 (a)]{G2019b}),
     $\C\left(d_N,\PP_{n,m},
     C_0(nV^*),L_p(nV^*)\right)
     <\iy$ for $n\in\N$.

     So \eqref{E7.12a} for $p\in(0,\iy]$ follows from
     Corollary \ref{C7.1} (v).
    \end{example}

    \begin{example} \label{Ex7.4}
    (Univariate
    inequalities of different weighted metrics).
    Let $d\mu(x)=\vert x\vert^{2\nu+1}dx,\,\nu\ge -1/2$,
    be the univariate weighted measure on $\R^1$. Next, let
    $F^{(1)}(\R^1)=C_0(\R^1),\,
    \be_1=1,\,
    F^{(2)}(\R^1)
    =L_{p,\mu}(\R^1),\,
    \be_2=\min\{1,p\}$,
    \ba
    B=B_{1,e}\cap L_{p,\mu}(\R^1)
    :=\{f\in B_{[-1,1]}:f(-x)=f(x),\,x\in\R^1\}\cap L_{p,\mu}(\R^1),
    \ea
    $0<p\le\iy,\,
     \Ome_n=[-n,n],\,
    \Ome_n^*=[-\tau n,\tau n],\,\tau\in(0,1),\,n\in\N$.

    Next, let
    $Be(\vphi)(x):=\vphi^{\prime\prime}
    (x)+((2\nu+1)/x)\vphi^\prime(x)$ be the Bessel operator
    and let
    $D=(Be)^N\left\vert_{\R^m}\right.=(Be)^N,\,
    D_n=(Be)^N\left\vert_{\Omega_n}\right.=(Be)^N,\,
     D_n^*=(Be)^N\left\vert_{\Omega_n^*}\right.=(Be)^N,\,N\in\Z^1_+
    ,\,n\in\N$.
    Further, let
    \ba
    \PP_{n,1,e}:=
    \{P\in \PP_{n,1}:P(-x)=P(x),\,x\in\R^1\}
    \ea
    be a subset of $\PP_{n,1}$ of all even
     univariate polynomials with complex coefficients, and let
     $B_n=\PP_{n,1,e}\left\vert_{\Omega_n}\right.$ and
    $B_n^*=\PP_{n,1,e}\left\vert_{\Omega_n^*}\right.,\,n\in\N.$

    Inequalities \eqref{E7.8b} were proved by the author
    \cite{G2019a} in the following form:
    \bna\label{E7.13a}
    \lim_{n\to\iy}
\sup_{P\in\PP_{2\lfloor n/2 \rfloor ,1,e}
\setminus\{0\}}\frac{\vert (Be)^N(P)(0)\vert}
{\|P\|_{L_{p,\mu}([-n,n])}}
    &\le &\sup_{f\in (B_{1,e}\cap L_{p,\mu}(\R^1))
    \setminus\{0\}}\frac{\vert(Be)^N(f)(0)\vert}
{\|f\|_{L_{p,\mu}(\R^1)}}\nonumber\\
    &\le & \lim_{n\to\iy}
\sup_{P\in\PP_{2\lfloor n/2 \rfloor ,1,e}
\setminus\{0\}}\frac{\vert (Be)^N(P)(0)\vert}
{\|P\|_{L_{p,\mu}([-\tau n,\tau n])}}
\ena
for $\tau\in(0,1)$.
By letting $\tau\to 1-$ in \eqref{E7.13a},
 we arrive at the following asymptotic result
 for $0<p\le\iy$
    (see \cite[Theorem 4.1]{G2019a}):
    \ba
\sup_{f\in (B_{1,e}\cap L_{p,\mu}(\R^1))
    \setminus\{0\}}\frac{\vert(Be)^N(f)(0)\vert}
{\|f\|_{L_{p,\mu}(\R^1)}}
    = \lim_{n\to\iy} n^{-2N-(2\nu+2)/p}
\sup_{P\in\PP_{2\lfloor n/2 \rfloor ,1,e}
\setminus\{0\}}\frac{\vert (Be)^N(P)(0)\vert}
{\|P\|_{L_{p,\mu}([-1,1])}}.
\ea
Verifications of Conditions C3.1)
and C12.1) in weighted spaces are provided similarly
 to Example \ref{Ex7.3} (see \cite{G2019a} for details).

 So \eqref{E7.13a} for $p\in(0,\iy]$ follows from
     Corollary \ref{C7.1} (v).
    \end{example}

    \begin{remark}\label{R7.5}
    Special cases of relations \eqref{E7.8a}, \eqref{E7.8b},
    and \eqref{E7.10}
    that are discussed
    in Examples
    \ref{Ex7.2}, \ref{Ex7.3}, and \ref{Ex7.4}
    and published in
    \cite{GT2017, G2017, G2018, G2019a, G2019b}
    follow from statements (ii), (v), and (vi) of
    Corollary \ref{C7.1}.
    In other words, the relations presented in
    the aforementioned publications can be proved
     by verifying the validity of the corresponding
     conditions C3.1) (or C3.1.1)) and C12.1).
     Despite the fact that these conditions were
     not formulated explicitly in
     these papers,
     the proofs of
     special cases of \eqref{E7.8a}, \eqref{E7.8b},
      and \eqref{E7.10}
     were actually based on the verification of
     these conditions.
     This observation along with the general approach
     to sharp constants of approximation theory
     developed in \cite{G1992, G2000} were the major
     motivations for writing this paper.
    \end{remark}

    \begin{remark}\label{R7.6}
    The author \cite[Corollary 4.5]{G2019a} proved
    the following limit relation ($p\in[1,\iy)$):
    \beq\label{E7.13c}
    \sup_{f\in \left(B_{\BB^m_1}\cap L_p(\R^m)\right)
    \setminus\{0\}}\frac{\|\Delta^N(f)\|_{L_\iy(\R^m)}}
{\|f\|_{L_p(\R^m)}}
    = \lim_{n\to\iy} n^{-2N-(m-1)/p}
\sup_{P\in\PP_{n,m}
\setminus\{0\}}
\frac{\|\delta^N(P)\|_{L_\iy\left(S^{m-1}\right)}}
{\|P\|_{L_p\left(S^{m-1}\right)}},
    \eeq
    where $S^{m-1}$ is the $(m-1)$-dimensional
    unit sphere in $\R^m,\,\de$ is the Laplace-Beltrami
    operator on $S^{m-1},\,\Delta$ is the
    Laplacian on $\R^m$, and $B=B_{\BB^m_1}$ is
    the set of all entire functions of spherical
    type.
    However, unlike Example \ref{Ex7.3}, the proof of
    \eqref{E7.13c} was based rather on an invariance
    theorem \cite[Corollary 3.11]{G2019a}
     and Example \ref{Ex7.4} than
    on Corollary \ref{C7.1} because
    constructing the sequence $\{f_n\}_{n=1}^\iy$
    from Condition C3.1.1 on the sphere
    turned out to be a difficult task.
    In case of $N=0$, \eqref{E7.13c}  was established in
    \cite[Theorem 1.1 (i)]{DGT2019} by a different method.
    For more general manifolds, it is even not easy
    to choose the "right" set $B$.
    \end{remark}

    \begin{remark}\label{R7.7a}
    Limit equalities for approximating elements are
    proved only for $F^{(1)}(\R^m)=L_\iy(\R^m)$
    in Example \ref{Ex7.2} and for
    $F^{(1)}(\R^m)=C_0(\R^m)$
    in Examples \ref{Ex7.3} and \ref{Ex7.4}.
    In recent paper \cite{G2021c},
    the author proved limit equalities in
    Example \ref{Ex7.3} in cases of a cube and a ball
    with $C_0(\R^m)$ replaced by $L_\iy(\R^m)$.
    The problem of replacing $L_\iy(\R^m)$
    and/or $C_0(\R^m)$ by $L_q(\R^m),\,q\in(0,\iy)$,
    in Examples \ref{Ex7.2}--\ref{Ex7.4}
    is still open.
    Actually, inequalities of the form
    $\C\le\liminf_{n\to\iy} \C_n$ can be proved
    for $L_q(\R^m),\,q\in(0,\iy),$ as well
     (see, for example,
    \cite[Theorem 1.4 and Conjecture 1.10]{GT2017}
    and \cite[Theorem 1.3 and Remark 1.8]{G2017}).
    However, we cannot prove for examples
    \ref{Ex7.2}--\ref{Ex7.4}
    relations of the form
    $\C\ge\limsup_{n\to\iy} \C_n$
    and $\C=\lim_{n\to\iy} \C_n$.
    The problem in these examples with
    $L_\iy(\R^m)$
    and/or $C_0(\R^m)$ replaced by
    $L_q(\R^m),\,q\in(0,\iy)$,
    is in proving equality (7.9),
    which seems to be a difficult task.
    \end{remark}

    \subsection{Approximation of Individual Elements
  (Problem B)}\label{S7.2}
  Here, we discuss a special case of Corollary \ref{C2.14} and
  three examples.\vspace{.12in}\\
 \textbf{General result.}
 We assume that all assumptions of Section \ref{S5.2} are held
 in the following special cases.
 First,
 we assume that $F^{(2)}_*=F(\R^m)
 =F(\R^m,\mu)$ is
 a $\be_2$-normed functional space
 with the monotone and $\boldsymbol{\om}$-continuous $\be_2$-norm
 (see Definitions \ref{D6.1} and \ref{D6.2}), $\be_2\in(0,1]$, and
 $F^{(2)}_{*n}=F(\Ome_n)$
 or $F^{(2)}_{*n}=F(\Ome_n^*),
 \,n\in\N $.
 Next, we define the restriction operators $L_n^{(2)}
 =R(\R^m,\Ome_n^*),\,
 l_n^{(2)}=R(\R^m,\Ome_n),
 \,n\in\N$, and
 \ba
 l_{n,s}^{(2)}=\left\{\begin{array}{ll}
 R(\Ome_n,\Ome_s),&n\ge s,\\
 R(\Ome_s,\Ome_n),&n< s,\end{array}\right.
 \quad n\in\N,\quad s\in\N.
 \ea

We also assume that $B,\,B_n$, and $B_n^*$
are nontrivial subspaces of $F(\R^m),\
 F(\Ome_n)$, and $F(\Ome_n^*)$,
 respectively, $n\in\N$.
 Next, for $\vphi\in F(\R^m)$ we set
 $\vphi_n:=L_n^{(2)}(\vphi)
 =\vphi\left|_{\Omega_n^*}\right.\in F(\Omega_n^*)$
 or
 $\vphi_n:=l_n^{(2)}(\vphi)
 =\vphi\left|_{\Omega_n}\right.\in F(\Omega_n)$.
 In addition, we assume that
 $\vphi\in  F(\R^m)
 \setminus B$ and
 $\vphi\left|_{\Omega_n}\right.
 \in F(\Omega_n)
 \setminus B_n$ or
 $\vphi_n\left|_{\Omega_n^*}\right.\in F(\Omega_n^*)
 \setminus B_n^*$.
 Note that to the end of the example we use the
 notation $\vphi$ instead of
 $\vphi\left|_{\Omega_n}\right.$ or
 $\vphi\left|_{\Omega_n^*}\right.$.

 The "norms" are defined by \eqref{E2.52} as follows.
 \beq\label{E7.13b}
  \|h\|_{F^{(2)}}:=\|\vphi-h\|_{F(\R^m)},\quad h\in F(\R^m);
  \qquad
 \|h_n\|_{F^{(2)}_n}:=\|\vphi-h_n\|_{F(\Omega)},
 \quad h_n\in F(\Omega),
 \eeq
 where $\Omega=\Omega_n$ or $\Omega=\Omega_n^*$.

 Furthermore, for $\vphi\in F(\R^m)$  we consider the errors of
 best approximation
 $E(\vphi,B,F(\R^m))$, $E(\vphi,B_n,F(\Ome_n))$, and
 $E(\vphi,B_n^*,F(\Ome_n^*)),\,n\in\N$ (see \eqref{E2.8*}
 or \eqref{E6.2a1}
 for the definition);
 by the assumptions on $\vphi$ and its restriction to
 $\Omega_n$ or $\Omega_n^*,\,n\in\N$,
 all these approximation errors are positive.

 Finally, let
$B_n^\prime\subseteq B_n\setminus\{0\}$ be a set
from Condition C9) such that
the following relation holds:
\ba
E(\vphi,B_n,F(\Ome_n))
=E(\vphi,B_n^\prime,F(\Ome_n)).
\ea

 Some conditions on families
 $\left\{B,F(\R^m)\right\},
 \left\{B_n,F(\Ome_n)\right\}$,
 and $\left\{B_n^*,F(\Ome_n^*)\right\},
 n\in\N$,
 are introduced below; they
 are needed for the validity of relations \eqref{E2.54}
 and/or \eqref{E2.55} in function spaces.\vspace{.12in}\\
 \emph{Condition C3.2*).} For every $f\in B\setminus\{0\}$
 there exists a sequence $f_n^*\in B_n^*,\,n\in\N$,
such that
\beq \label{E7.14}
\lim_{n\to\iy}
\left\|f-f_n^*\right\|_{F(\Ome_n^*)}=0.
\eeq
\emph{Condition C10.2).} (Compactness condition)
For any sequence $\{f_n\}_{n=1}^\iy$ such that
$f_n\in B_n^\prime,\,n\in\N$,
there exist a sequence $\{n_k\}_{k=1}^\iy$ of natural
numbers and an element $f\in B$
such that
\beq
\lim_{k\to\iy}\left\|f_{n_k}-f\right\|_{F
\left(\Ome_{n_k}\right)}=0.
\label{E7.15}
\eeq
\emph{Condition C12.2).} (Weak compactness condition)
For any sequence $\{f_n\}_{n=1}^\iy$
such that
$f_n\in B_n^\prime,\,n\in\N$,
there exists $f\in B$ such that
for any fixed $s\in\N$,
there exists a sequence $\{n_k\}_{k=1}^\iy$ of natural numbers
 ($n_k=n_k(s),\,k\in\N$)
such that
\beq
\lim_{k\to\iy}\left\|f_{n_k}-f
\right\|_{F\left(\Ome_{s}\right)}=0,\label{E7.15a}
\eeq
 The following corollary holds true.
 \begin{corollary}\label{C7.6}
 (i) If Condition C3.2*)
 is satisfied, then
 \beq\label{E7.16}
 \limsup_{n\to\iy} E\left(\vphi,B_n^*,F(\Ome_n^*)\right)
 \le E\left(\vphi,B,F(\R^m)\right).
 \eeq
 (ii) If
Condition C10.2) is satisfied, then
 \beq\label{E7.17}
 \liminf_{n\to\iy} E\left(\vphi,B_n,F(\Ome_n)\right)
 \ge E\left(\vphi,B,F(\R^m)\right).
\eeq
(iii) If
Condition C12.2) is satisfied,
 then \eqref{E7.17} holds
true.\\
(iv)
Let
Condition C3.2*) and either C10.2) or C12.2)  be satisfied.
 Then
 \beq\label{E7.18}
 \limsup_{n\to\iy} E\left(\vphi,B_n^*,F(\Ome_n^*\right))
 \le E\left(\vphi,B,F(\R^m)\right)
 \le \liminf_{n\to\iy} E\left(\vphi,B_n,F(\Ome_n\right)).
\eeq
 \end{corollary}
 \proof
 Note first that in this proof all elements in discussed
 Conditions C1) through C13) are replaced without further notice
 by the corresponding elements that are used in Conditions
 C3.2*), C10.2), C12.2) and relations \eqref{E7.16},
 \eqref{E7.17}, \eqref{E7.18}. For example, in Condition
  C3*)  $B_n$ is replaced by $B_n^*$.

 By the assumptions about the functional spaces and operators
 introduced above, Conditions C4*), C6), C6b), C11), and C11a)
 of Section \ref{S3.2} are satisfied. Indeed,
 Conditions C6) and C11) are satisfied by the direct
 definitions of the operators whose existence is required
 in these conditions.
 In addition,
 Conditions C4*), C6b),  and C11a) are satisfied
  by Proposition \ref{P6.4}.
  For example, for every $h\in F^{(2)}$ we obtain by
  \eqref{E7.13b} and by Proposition \ref{P6.4} (ii)
\ba
\liminf_{n\to\iy}\left\|l_n^{(2)}(h)\right\|_{ F^{(2)}_n}
&=&\liminf_{n\to\iy}\left\|\vphi_n-l_n^{(2)}(h)\right\|_{ F(\Omega_n)}
=\liminf_{n\to\iy}\left\|l_n^{(2)}(\vphi-h)\right\|_{ F(\Omega_n)}\nonumber\\
&=&\|\vphi-h\|_{ F(\R^m)}
= \|h\|_{F^{(2)}}.
\ea
Therefore, Condition C6b) is satisfied. Similarly we can prove that
Conditions C4*) and C11a) are satisfied (see also Remark \ref{R2.17a}).

  Next, Conditions C3.2*), C10.2), and C12.2) of this section
  are equivalent to Conditions C3*),
  equation \eqref{E2.21} of C10), and
  equation \eqref{E2.30} of C12) from Section
  \ref{S3.2}, respectively.
  Therefore, statements (i), (ii), and (iii)
  of Corollary \ref{C7.6}
   follow immediately
  from Corollary \ref{C2.14}, while statement (iv) is
  a direct corollary of (i), (ii), and (iii).
  \hfill$\Box$
  \begin{remark}\label{R7.7}
  Statements (i) and (ii) of Corollary \ref{C7.6} with
  slightly different conditions were obtained in
  \cite[Corollary 3.1]{G1992} for normed spaces and
  in
  \cite[Sect. 1.4.3]{G2000}
  for $\be_2$-normed spaces.
  \end{remark}
  \noindent
\textbf{Examples.} There are numerous publications
  \cite{B1913, B1938, B1946(trig), B1946, B1947(lim), T1963,
  A1965, R1968,  G1982, G1991,
   G1992,
  G2000, G2002, G2008} about special cases of Corollary
    \ref{C7.6} and their modifications;
    here, we call them limit theorems of approximation theory
    (see also Section \ref{S2.3} B).
    Examples \ref{Ex7.8}, \ref{Ex7.9}, and \ref{Ex7.10} below
    discuss relations \eqref{E7.18}
    proved in these papers in various situations.

    \begin{example}\label{Ex7.8}
    (Limit theorems for multivariate periodic approximations).
    Let
    $F(\R^m)=L_p(\R^m),\linebreak
    \be_2=\min\{1,p\}$,
     and let $\vphi\in L_p(\R^m)\cap M_{C,N}$
     be a fixed function$,\,
     0<p\le\iy$.
     Next, let
     $ B=B_V\cap L_p(\R^m),\,
     0<p\le\iy,\,
    \Ome_n=\Ome_n^*=Q_{\g(n)}^m,\,
    B_n=B_n^*=\TT_{n,V}\cap M_{C_9C,N+2}\left\vert_{Q_{\g(n)}^m}\right.
    =\TT_{n,V}\cap M_{C_9C,N+2},\,
    n\in\N.$
    Here, $C>0$ and $N\ge 0$ are fixed constants;
    $C_9>0$ is independent
    of $n$ and $\vphi$. In addition,
    a number sequence $\{\g(n)\}_{n=1}^\iy$ satisfies
    the conditions
    \beq\label{E7.18a}
    0<\g(n)\le\pi n,\quad n\in\N,\qquad \lim_{n\to\iy}\g(n)=\iy,\qquad
    \lim_{n\to\iy}(\g(n))^{2+m/p}n^{-1/(m+2)}=0.
    \eeq

Relations \eqref{E7.18} were proved by the author \cite{G1992, G2000}
in the following form:
there exists a constant $C_{10}>0$, which is
independent of $n$ and $\vphi$, such that for all $C_9\ge C_{10}$,

    \beq\label{E7.19}
    \lim_{n\to\iy}E\left(\vphi,\TT_{n,V}\cap M_{C_9C,N+2},
    L_p\left(Q_{\g(n)}^m\right)\right)
    =E\left(\vphi,B_V\cap L_p(\R^m),L_p(\R^m)\right).
    \eeq
    For $p\in[1,\iy]$ equality \eqref{E7.19} was proved
    in \cite[Theorem 5.1]{G1992} in a more general situation;
    the case of $p\in(0,1)$ was discussed in
    \cite[Sect. 11.4.4.5]{G2000}. For $m=1,\,V=[-\sa,\sa]$, and
    $p=\iy$, \eqref{E7.19} was established by Bernstein
    \cite[Theorem 5]{B1946(trig)}.

    Condition C3.2*) of approximation of entire
     functions of exponential type $V$
     by multivariate trigonometric polynomials was verified in
      \cite[Lemma 5.1]{G1992} (see also \cite[Sect. 3]{G2018}).

      In addition, we choose
      $B_n^\prime=B_n=\TT_{n,V}\cap M_{C_9C,N+2}$.
Since $\TT_{n,V}\subseteq B_V$,
      the compactness property of $\cup_{n=1}^\iy B_n^\prime$, stated in C12.2),
     can be found in \cite[Lemma 2]{G1981} and
     \cite[Lemma 2.3]{G2018}
      (a special case is discussed in \cite[Theorem 3.3.6]{N1969}).

     So \eqref{E7.19} for $p\in(0,\iy]$ follows from Corollary
     \ref{C7.6} (iv).
    \end{example}

    \begin{example} \label{Ex7.9}
    (Limit theorems for multivariate polynomial approximations).
    Let
    $F(\R^m)=L_p(\R^m),\,
    \be_2=\min\{1,p\},\,
   B=B_V\cap L_p(\R^m),\,
   0<p\le\iy,\,
   \Ome_n=nV^*,\,
   \Ome_n^*=\tau nV^*,\,\tau\in(0,1),\,
    B_n=\PP_{n,m}\left\vert_{\Omega_n}\right.=\PP_{n,m},\,
    B_n^*=\PP_{n,m}\left\vert_{\Omega_n^*}\right.=\PP_{n,m},\,
    n\in\N.$

    Relations \eqref{E7.18}
    were proved in
    \cite{G1982}. The following result is a corollary of
    \eqref{E7.18}:
    if $\vphi\in L_\iy(\R^m)$, then
    \bna\label{E7.20}
    &&E\left(\vphi,B_{ V}\cap L_\iy(\R^m),L_\iy(\R^m)\right)
    \le \liminf_{n\to\iy}E\left(\vphi,\PP_{n,m},
    L_\iy\left(nV^*\right)\right)\nonumber\\
    &&\le \limsup_{n\to\iy}E\left(\vphi,\PP_{n,m},
    L_\iy\left(nV^*\right)\right)
   \le \lim_{\tau\to 1^-}
    E\left(\vphi,B_{\tau V}\cap L_\iy(\R^m),L_\iy(\R^m)\right),
    \ena
    and if $\vphi\in L_p(\R^m)\cap L_\iy(\R^m)
    ,\,p\in(0,\iy)$, then
     \beq\label{E7.21}
     \lim_{n\to\iy}E\left(\vphi,\PP_{n,m},
    L_p\left(nV^*\right)\right)
    =E\left(\vphi,B_{V}\cap L_p(\R^m),L_p(\R^m)\right).
    \eeq
    For $p\in[1,\iy]$ relations \eqref{E7.20} and
    \eqref{E7.21} were proved
    in \cite[Theorems 5.1 and 5.2]{G1982};
    the case of $p\in(0,1)$ was discussed
    in a more general and precise form in
    \cite[Theorem 11.3]{G2000}. For $m=1,\,V=[-\sa,\sa]$,
    and $V^*=[-1/\sa,1/\sa]$,
    inequalities \eqref{E7.20} were established by Bernstein
    \cite[Theorem 6]{B1946}. Equality \eqref{E7.21}
    for $m=1,\,V=[-\sa,\sa],\,V^*=[-1/\sa,1/\sa]$,
     and $p\in[1,\iy)$
    was proved by Raitsin \cite{R1968}.
     More precise relations than \eqref{E7.20} were proved in
    \cite{G1991}.

    Condition C3.2*) of approximation of entire
     functions of exponential type $V$
     by multivariate algebraic polynomials was verified in
      \cite[Lemma 4.4]{G1982}
      (see also \cite[Lemma 2.4]{G2019b}).

      In addition, we choose
      \ba
       B_n^\prime
      =\left\{f_n\in\PP_{n,m}:
\|f_n\|_{L_p(nV^*)}
\le \left(2^{\be_2}+1\right)^{1/\be_2}\|\vphi\|_{L_p(\R^m)}\right\},
\qquad n\in\N.
      \ea
      Then
      \beq\label{E7.21a}
      E(\vphi,\PP_{n,m},L_p(nV^*))
      = E(\vphi,B_n^\prime,L_p(nV^*)),\qquad n\in\N.
      \eeq

      Indeed, given $\vep\in\left(0,\|\vphi\|_{L_p(\R^m)}\right)$
      and $n\in\N$, let $f_{n,\vep}\in \PP_{n,m}$ satisfy the inequalities
      \ba
      \|\vphi\|_{L_p(\R^m)}\ge E\left(\vphi,\PP_{n,m},L_p
      \left(nV^*\right)\right)
      > \|\vphi-f_{n,\vep}\|_{L_p\left(nV^*\right)}-\vep.
      \ea
      Then using the $\be_2$-triangle inequality,
      we see that $f_{n,\vep}\in B_n^\prime$,
      and, in addition,
      \ba
      E\left(\vphi,\PP_{n,m},L_p\left(nV^*\right)\right)
      =\inf_{\vep\in\left(0,\|\vphi\|_{L_p(\R^m)}\right)}
      \|\vphi-f_{n,\vep}\|_{L_p\left(nV^*\right)}
      \ge E\left(\vphi,B_n^\prime,L_p\left(nV^*\right)\right).
      \ea
      Therefore, \eqref{E7.21a} holds true.

      The compactness property of $\cup_{n=1}^\iy B_n^\prime$, stated in C12.2),
     was verified by Bernstein \cite{B1946} for $m=1$
     and by the author \cite[Lemma 2.2]{G2019b} for $m>1$.

     So \eqref{E7.18} for $p\in(0,\iy]$ follows from Corollary
     \ref{C7.6} (iv).
    \end{example}

    \begin{example}\label{Ex7.10}
    (Limit theorems for univariate polynomial approximations
    with exponential
    \linebreak weight).
    Let  $Q:\R^1\to [0,\iy)$ be an even differentiable
    convex function and let $\be_2$-norms
    $\|\cdot\|_{L_{p}(W_p,\R^1)},\,p\in(0,\iy],\,\be_2=\min\{1,p\}$,
     be defined by
    \eqref{E6.2c} and \eqref{E6.2d}.

 In addition, let $a_n=a_n(Q)$ be the $n$th
 Mhaskar-Rakhmanov-Saff number defined as the
 positive root of the equation
 \ba
 n=\frac{2}{\pi}\int_0^1\frac{a_nx\,Q^\prime(a_nx)}
 {\sqrt{1-x^2}}dx
 \ea
 and let
 \ba
 b_n=b_n(Q):=\frac{2}{\pi}\int_0^1\frac{Q^\prime(a_nx)\sqrt{1-x^2}}
 {x}dx+n/a_n.
 \ea
 Next, let
    $F^{(2)}_*=L_p(\R^1),\,
    B=B_{[-\sa,\sa]}\cap L_p(\R^1),\,
    F^{(2)}_{*n}=L_p\left(W_p((\sa/b_n)\cdot),\R^1\right),\,
    \sa>0,\,
    \Omega_n=\Omega_n^*=\R^1,\,
    B_n=\PP_{n,1},\,
    B_n^*=\PP_{\lfloor n/\tau\rfloor,1},\,\tau\in(0,1),\,
    n\in\N,\,
    p\in(0,\iy]$.

Relations like \eqref{E7.18}
    were proved in
    \cite{G2008}. The following result is a corollary of
    these relations: let $W_p$ satisfy some additional conditions
    (see \cite[Definition 1.4.6]{G2008});
    if $\vphi\in L_\iy(\R^1)$, then
    \bna\label{E7.22}
    &&E\left(\vphi,B_{[-\sa,\sa]}\cap L_\iy(\R^1),L_\iy(\R^1)\right)\nonumber\\
    &&\le \liminf_{n\to\iy}E\left(\vphi((b_n/\sa)\cdot),
    \PP_{n,1},
    L_\iy\left(W_\iy,\R^1\right)\right)\nonumber\\
    &&\le \limsup_{n\to\iy}E\left(\vphi((b_n/\sa)\cdot),
    \PP_{n,1},
    L_\iy\left(W_\iy,\R^1\right)\right)\nonumber\\
   &&\le \lim_{\tau\to 1^-}
    E\left(\vphi,B_{[-\tau\sa,\tau\sa]}\cap L_\iy(\R^1),L_\iy(\R^1)\right),
    \ena
    and if $\vphi\in L_p(\R^1)
    ,\,p\in(0,\iy)$, then
     \beq\label{E7.23}
     \lim_{n\to\iy}(b_n/\sa)^{1/p}E\left(\vphi((b_n/\sa)\cdot),
    \PP_{n,1},
    L_p\left(W_p,\R^1\right)\right)
    =E\left(\vphi,B_{[-\sa,\sa] }\cap L_p(\R^1),L_p(\R^1)\right).
    \eeq

    For $p\in(0,\iy)$ equality \eqref{E7.23} was proved
    in \cite[Theorem 2.1.2]{G2008}, while a stronger version
    of relations \eqref{E7.22} for $p=\iy$
    was established in
    \cite[Theorem 2.1.1]{G2008}. For the Hermite weight
     similar results were obtained in \cite{G2002}.

    Condition C3.2*) of polynomial approximation
    with exponential weights of entire
     functions of exponential type follows from
      \cite[Theorem 2.2.1]{G2008}.
      The sets $B_n^\prime,\,n\in\N$, can be chosen
      like in Example \ref{Ex7.9}.
      The compactness property of
      $\cup_{n=1}^\iy B_n^\prime$, stated in C12.2),
      is based on \cite[Theorem 2.3.2]{G2008}
     and verified in \cite[Sect. 7.1]{G2008}.

     So \eqref{E7.18} for $p\in(0,\iy]$ follows from Corollary
     \ref{C7.6} (iv).
    \end{example}

    \begin{remark}\label{R7.10}
    Special cases of relations \eqref{E7.18}
    that are discussed
    in Examples
    \ref{Ex7.8}, \ref{Ex7.9}, and \ref{Ex7.10}
    and published in
    \cite{B1946, R1968, G1982, G2000, G2008}
    follow from statement (iv) of
    Corollary \ref{C7.6}.
    In other words, the relations presented in
    the aforementioned publications can be proved
     by verifying the validity of the corresponding
     conditions C3.2*)  or/and C12.2).
Despite the fact that these conditions were
     not formulated explicitly in
     these papers,
     the proofs of
     special cases of \eqref{E7.18}
     were actually based on the verification of
     these conditions.
     This observation was the major
     motivation for general approach developed in
     \cite {G1992, G2000}.
    \end{remark}

    \subsection{Approximation on Classes of Elements
  (Problem C)}\label{S7.3}
  Here, we discuss a special case of Corollary \ref{C2.19} and
  two examples.\vspace{.12in}\\
 \textbf{General result.}
 We assume that all assumptions of Section \ref{S5.3} are held
 in the following special cases.
 First,
 we assume that $F^{(1)}_*=F(\R^m)
 =F(\R^m,\mu)$ is
 a $\be_1$-normed functional space
 with the monotone and $\boldsymbol{\om}$-continuous
 $\be_1$-norm
 (see Definitions \ref{D6.1} and \ref{D6.2}),
 $\be_1\in(0,1]$, and
 $F^{(1)}_{*n}=F(\Ome_n)$
 or $F^{(1)}_{*n}=F(\Ome_n^*),
 \,n\in\N $.
 Next, we define the restriction operators $L_n^{(1)}
 =R(\R^m,\Ome_n^*),\,
 l_n^{(1)}=R(\R^m,\Ome_n),
 \,n\in\N$, and
 \ba
 l_{n,s}^{(1)}=\left\{\begin{array}{ll}
 R(\Ome_n,\Ome_s),&n\ge s,\\
 R(\Ome_s,\Ome_n),&n< s,\end{array}\right.
 \quad n\in\N,\quad s\in\N.
 \ea

We also assume that $B,\,B_n$, and $B_n^*$
are nontrivial subspaces of $F(\R^m),\,
 F(\Ome_n)$, and $F(\Ome_n^*)$,
 respectively, $n\in\N$.
 In addition, let  $G$ be a subspace of $F(\R^m),\,
 G_n$ be a subspace of $F(\Ome_n)$,
 and $G_n^*$ be a subspace of $F(\Ome_n^*),\,n\in\N$.
 Next, for $f\in B,\,f_n\in B_n$, and $f_n^*\in B_n^*$
  we consider the errors of
 best approximation
 $E(f,G,F(\R^m)),\,E(f_n,G_n,F(\Ome_n))$, and
 $E(f_n^*,G_n^*,F(\Ome_n^*)),\,n\in\N$ (see \eqref{E2.8*}
 or \eqref{E6.2a1}
 for the definition).
 Note that each of these
 errors of
 best approximation are $\be_1$-seminorms
 on the corresponding spaces.

 The "norms" are defined by \eqref{E2.60} and \eqref{E2.60a}
 as follows.
 \beq\label{E7.38aa}
 \|h\|_{F^{(1)}}=E(h,G,F(\R^m)),\quad h\in F(\R^m),\quad
 \|h_n\|_{F^{(1)}_n}=E(h_n,\tilde{G},F(\Omega)),
 \quad h_n\in F(\Omega),
 \eeq
 where $\Omega=\Omega_n$ or $\Omega=\Omega_n^*$
 and $\tilde{G}=G_n$ or $\tilde{G}=G_n^*
 ,\,n\in\N$.
 Note that $\tilde{G}$ in \eqref{E7.38aa} can depend on $h_n$.

 Finally, let
$B_n^\prime\subseteq B_n\setminus\{0\}$ be a set
from Condition C8) such that
the following relation holds:
\ba
\sup_{f_n\in B_n}E(f_n,G_n,F(\Ome_n))
=\sup_{f_n\in B_n^\prime}E(f_n,G_n,F(\Ome_n)).
\ea

 Below we introduce  some conditions on families
 $\left\{B,F(\R^m)\right\},
 \left\{B_n,F(\Ome_n)\right\}$,
 and $\left\{B_n^*,F(\Ome_n^*)\right\},\linebreak
 n\in\N$,
 that are needed for the validity of relations \eqref{E2.64}
 and/or \eqref{E2.65} in function spaces.\vspace{.12in}\\
 \emph{Condition C3.3).} For every $f\in B\setminus\{0\}$
 there exists a sequence $f_n^*\in B_n^*,\,n\in\N$,
such that
\beq \label{E7.24}
\lim_{n\to\iy}
\left\|f-f_n^*\right\|_{F(\Ome_n^*)}=0.
\eeq
\emph{Condition C4.3).} For every $h\in B$,
\beq \label{E7.25}
\liminf_{n\to\iy}
E\left(h\left\vert_{\Omega_n^*}\right.,G_n^*,F(\Ome^*_n)\right)
=\liminf_{n\to\iy}
E\left(h,G_n^*,F(\Ome^*_n))\ge E(h,G,F(\R^m)\right).
\eeq
\emph{Condition C6.3*a).} For every $h\in F(\R^m)$,
\beq \label{E7.25a}
\limsup_{n\to\iy}
E\left(h\left\vert_{\Omega_n}\right.,G_n,F(\Ome_n)\right)
=\limsup_{n\to\iy}
E\left(h,G_n,F(\Ome_n))\le E(h,G,F(\R^m)\right).
\eeq
\begin{remark}\label{R7.11}
Condition C4.3) can be verified by using statements
(ii) or (iii) of Corollary \ref{C7.6},
while Condition C6.3*a) can be verified by using statement
(i) of Corollary \ref{C7.6}.
\end{remark}
\noindent
\emph{Condition C10.3*a).} (Compactness condition)
For any sequence $\{f_n\}_{n=1}^\iy$ such that
$f_n\in B_n^\prime,\,n\in\N$,
there exist a sequence $\{n_k\}_{k=1}^\iy$ of natural
numbers and an element $f\in B$
such that
\beq \label{E7.26}
\lim_{k\to\iy}\left\|f_{n_k}-f\right\|_{F
\left(\Ome_{n_k}\right)}=0.
\eeq
 The following corollary holds true.
 \begin{corollary}\label{C7.12}
 (i) If Conditions C3.3) and C4.3)
 are satisfied, then
 \beq\label{E7.28}
 \liminf_{n\to\iy} \sup_{f_n\in B_n^*}
 E\left(f_n,G_n^*,F(\Ome_n^*)\right)
 \ge \sup_{f\in B}
 E\left(f,B,F(\R^m)\right).
 \eeq
 (ii) Let $\sup_{f_n\in B_n}
 E\left(f_n,G_n,F(\Ome_n)\right)$ be a
 finite number for every $n\in \N$.
 If Conditions C6.3*a) and C10.3*a) are
  satisfied, then
 \beq\label{E7.29}
 \limsup_{n\to\iy} \sup_{f_n\in B_n}
 E\left(f_n,G_n,F(\Ome_n)\right)
 \le \sup_{f\in B}
 E\left(f,B,F(\R^m)\right).
\eeq
(iii) Let $\sup_{f_n\in B_n}
 E\left(f_n,G_n,F(\Ome_n)\right)$ be a
 finite number for every $n\in \N$.
If
Conditions C3.3), C4.3), C6.3*a),
and C10.3*a) are satisfied, then
 \beq\label{E7.30}
 \liminf_{n\to\iy} \sup_{f_n\in B_n^*}
 E\left(f_n,G_n^*,F(\Ome_n^*)\right)
 \ge \sup_{f\in B}
 E\left(f,B,F(\R^m)\right)
 \ge \limsup_{n\to\iy} \sup_{f_n\in B_n}
 E\left(f_n,G_n,F(\Ome_n)\right).
\eeq
 \end{corollary}
 \proof
 Note first that in this proof all elements in discussed
 Conditions C1) through C13) and C10*a)
  are replaced without further notice
 by the corresponding elements that are used in Conditions
 C3.3), C4.3), C6.3a), C10.3*a) and relations \eqref{E7.28},
 \eqref{E7.29}, \eqref{E7.30}. For example, in Condition
  C3) $B_n$ is replaced by $B_n^*$.

 Next, Conditions C6) and C11) of Section \ref{S3.2}
  are satisfied by the direct
 definitions of the operators whose existence is required
 in these conditions.
  Further, Conditions C3.3), C4.3), C6.3*a),
  and C10.3*a) of this section
  are equivalent to Conditions C3), C4), C6*a)
   and C10*a) of Section
  \ref{S3.2}, respectively.
  Therefore, statements (i) and (ii)
  of Corollary \ref{C7.12}
   follow immediately
  from Corollary \ref{C2.19}, while statement (iii) is
  a direct corollary of (i) and (ii).
  \hfill$\Box$
  \begin{remark}
  Statements (i) and (ii) of Corollary \ref{C7.12} with
  slightly different conditions were obtained in
  \cite[Corollary 3.1]{G1992} for normed spaces and
  in
  \cite[Sect. 1.4.3]{G2000}
  for $\be_1$-normed spaces.
  \end{remark}
  \noindent
  \textbf{Examples.} There is a number of publications
  \cite{B1946(Nik), B1947(const), T1963,
  G1992, G2000, G2003} about special cases of Corollary
    \ref{C7.6} and their modifications;
    here, we call them
    relations between constants of approximation theory
    (see also Section \ref{S2.3} C).
    Examples \ref{Ex7.15} and \ref{Ex7.16} below
    discuss relations like \eqref{E7.28}
    proved in these papers in various situations.

    \begin{example}\label{Ex7.15}
    (Periodic approximation on $H_\om(\R^m)\cap \tilde{L}_{\iy,2\pi n,m}$).
    Let
    $F(\R^m)=L_\iy(\R^m),\,
    \be_1=1,\,
    B=H_{\om,0}(\R^m)\cap L_\iy(R^m),\,
    B_n^*=H_{\om,0}(\R^m)\cap \tilde{L}_{\iy,2\pi n,m},\,
    G=B_V\cap L_\iy(\R^m),\,
    G_n^*=\TT_{n,V}\cap M_{C_9\om(2),3},\,
    \Ome_n^*=Q_{\pi n/2}^m,\,
   n\in\N.$

   Let us choose $M=M(n)=\pi n/2,\,n\in\N$,
   and let $f\in H_{\om,0}(\R^m)\cap L_\iy(R^m)$.
Then by Lemma \ref{L6.3c}, there exists a function
$f_n^*=\tilde{f}:\R^m\to \R^m$
such that $f_n^*=f$ on $Q_{\pi n/2}^m$ and
$f_n^*\in H_{\om,0}(\R^m)\cap \tilde{L}_{\iy,2\pi n,m}$.
Since $f-f_n^*=0$ on $Q_{\pi n/2}^m$,
Condition C3.3) is satisfied.

  Next, Condition C4.3) is satisfied as well. Indeed,
  let a number sequence $\{\g(n)\}_{n=1}^\iy$
  satisfy \eqref{E7.18a}.
  Then
  for all large enough $n,\,\g(n)\le \pi n/2$,
  since by \eqref{E7.18a}, $\g(n)=o(n)$ as $n\to\iy$.
  In addition, for any
  $h\in H_{\om,0}(\R^m)\cap L_\iy(\R^m)$,
  we have $h\in M_{\om(2),1} \cap L_\iy(\R^m)$,
   by \eqref{E6.5aa}.
   Then using relation \eqref{E7.19}
   from Example
   \ref{Ex7.8} for $p=\iy,\,C=\om(2)$,
   and $N=1$,
    we obtain
   \ba
  && \liminf_{n\to\iy}E\left(h,\TT_{n,V}\cap M_{C_9\om(2),3},
    L_\iy\left(Q_{\pi n/2}^m\right)\right)\\
    &&\ge \lim_{n\to\iy}E\left(h,\TT_{n,V}\cap M_{C_9\om(2),3},
    L_\iy\left(Q_{\g(n)}^m\right)\right)\\
   && = E\left(h,B_V\cap L_\iy(\R^m),L_\iy(\R^m)\right),
    \qquad
    h\in H_{\om,0}(\R^m)\cap L_\iy(\R^m).
   \ea
   Thus Condition C4.3) is satisfied.
   We recall that $C_9$ is the
   same constant as in Example
   \ref{Ex7.8}, i.e.,
   there exists $C_{10}>0$ such that for any
   $C_9\ge C_{10}$, \eqref{E7.19}
   is valid for $\vphi=h$; in addition,
   $C_9$ and $C_{10}$ are independent of $n$ and $h$.

   Therefore using Corollary \ref{C7.12} (i), we obtain
   \bna\label{E7.31}
   &&\liminf_{n\to\iy}
   \sup_{f_n\in H_{\om,0}(\R^m)\cap \tilde{L}_{\iy,2\pi n,m}}
   E\left(f_n,\TT_{n,V}\cap M_{C_9\om(2),3},
   L_\iy\left(Q_{\pi n}^m\right)\right)\nonumber\\
   &&\ge\sup_{f_n\in H_{\om,0}(\R^m)\cap \tilde{L}_{\iy,2\pi n,m}}
   E\left(f_n,\TT_{n,V}\cap M_{C_9\om(2),3},
   L_\iy\left(Q_{\pi n/2}^m\right)\right)\nonumber\\
  && \ge  \sup_{f\in H_{\om,0}(\R^m)\cap L_\iy(\R^m)}
   E\left(f,B_V\cap L_\iy(\R^m) ,L_\iy\left(\R^m\right)\right).
   \ena
   Furthermore, we prove that for any
   $f_n\in H_{\om,0}(\R^m)\cap \tilde{L}_{\iy,2\pi n,m},\,
   n\in\N$, and for large enough $C_9\ge C_{10}$,
   the following equality holds true:
   \beq\label{E7.31a}
    E\left(f_n,\TT_{n,V}\cap M_{C_9\om(2),3},
   L_\iy\left(Q_{\pi n}^m\right)\right)
   = E\left(f_n,\TT_{n,V},
   L_\iy\left(Q_{\pi n}^m\right)\right).
   \eeq
   Indeed, let us set
   \ba
   \TT_{n,V,1}:=\left\{g_n\in \TT_{n,V}:
   E\left(f_n,\TT_{n,V},
   L_\iy\left(Q_{\pi n}^m\right)\right)
   \ge \|f_n-g_n\|_{L_\iy\left(Q_{\pi n}^m\right)}-1\right\}.
   \ea
   If $g_n\in  \TT_{n,V,1}$, then by \eqref{E6.5aa},
   \bna\label{E7.31b}
   \vert g_n(x)\vert
   &\le&  E\left(f_n,\TT_{n,V},
   L_\iy\left(Q_{\pi n}^m\right)\right)+1+\vert f_n(x)\vert
   \nonumber\\
   &\le& E\left(f_n,\TT_{n,V},
   L_\iy\left(Q_{\pi n}^m\right)\right)+1+
   \om(2)(1+\vert x\vert),\qquad x\in\R^m.
   \ena
   Since by Jackson's theorem
   (see \cite[Theorem 5.3.2]{N1969}),
   \bna\label{E7.31c}
  && \sup_{n\in\N}
   \sup_{f_n\in H_{\om,0}(\R^m)\cap \tilde{L}_{\iy,2\pi n,m}}
   E\left(f_n,\TT_{n,V},
   L_\iy\left(Q_{\pi n}^m\right)\right)\nonumber\\
   &&=\sup_{n\in\N}
   \sup_{f_n\in H_{\om,0}(\R^m)\cap \tilde{L}_{\iy,2\pi n,m}}
   E\left(f_n(n\cdot),\TT_{nV},
   L_\iy\left(Q_{\pi}^m\right)\right)
   <\iy,
   \ena
   we see from \eqref{E7.31b} and \eqref{E7.31c} that
   $\TT_{n,V,1}\subseteq M_{C_{11},1}$, where
   $C_{11}$ is independent of $n$.

   Then setting
   $C_9=\max\{C_{10},C_{11}/\om(2)\}$,
   we obtain
   \bna\label{E7.31d}
   && E\left(f_n,\TT_{n,V}\cap M_{C_9\om(2),3},
   L_\iy\left(Q_{\pi n}^m\right)\right)
   \le  E\left(f_n,\TT_{n,V}\cap M_{C_{11},1},
   L_\iy\left(Q_{\pi n}^m\right)\right)\nonumber\\
  && \le  E\left(f_n,\TT_{n,V,1},
   L_\iy\left(Q_{\pi n}^m\right)\right)
   =  E\left(f_n,\TT_{n,V},
   L_\iy\left(Q_{\pi n}^m\right)\right).
   \ena
   Thus \eqref{E7.31a} follows from \eqref{E7.31d}.

   Finally combining \eqref{E7.31} and \eqref{E7.31a},
   we obtain
   \bna\label{E7.31e}
   &&\liminf_{n\to\iy}
   \sup_{f_n\in H_{\om}(\R^m)\cap \tilde{L}_{\iy,2\pi n,m}}
   E\left(f_n,\TT_{n,V},
   L_\iy\left(Q_{\pi n}^m\right)\right)\nonumber\\
   &&=\liminf_{n\to\iy}
   \sup_{f_n\in H_{\om,0}(\R^m)\cap \tilde{L}_{\iy,2\pi n,m}}
   E\left(f_n,\TT_{n,V},
   L_\iy\left(Q_{\pi n}^m\right)\right)\nonumber\\
  && \ge  \sup_{f\in H_{\om,0}(\R^m)\cap L_\iy(\R^m)}
   E\left(f,B_V\cap L_\iy(\R^m) ,L_\iy\left(\R^m\right)\right)\nonumber\\
  && =\sup_{f\in H_{\om}(\R^m)\cap L_\iy(\R^m)}
   E\left(f,B_V\cap L_\iy(\R^m) ,L_\iy\left(\R^m\right)\right).
   \ena

Next, by Lemma \ref{L6.3b}, for any
   $f_n\in H_\om(\R^m)\cap \tilde{L}_{\iy,2\pi n,m}$,
   \beq\label{E7.32}
   E\left(f_n,\TT_{n,V},L_\iy\left(Q_{\pi n}^m\right)\right)
   =E\left(f_n,B_V\cap L_\iy(\R^m) ,L_\iy\left(\R^m\right)\right),
   \qquad n\in\N.
   \eeq
   Taking into account a trivial inclusion
   $H_\om(\R^m)\cap \tilde{L}_{\iy,2\pi n,m}
   \subseteq H_\om(\R^m)\cap L_\iy(\R^m)$,
   we obtain from \eqref{E7.32}
   \bna\label{E7.33}
   &&\limsup_{n\to\iy}
   \sup_{f_n\in H_\om(\R^m)\cap \tilde{L}_{\iy,2\pi n,m}}
   E\left(f_n,\TT_{n,V},L_\iy\left(Q_{\pi n}^m\right)
   \right)\nonumber\\
  && \le  \sup_{f\in H_\om(\R^m)\cap L_\iy(\R^m)}
   E\left(f,B_V\cap L_\iy(\R^m) ,L_\iy
   \left(\R^m\right)\right).
   \ena
   Combining \eqref{E7.31e} and \eqref{E7.33}, we arrive at
   \bna\label{E7.34}
   &&\lim_{n\to\iy}
   \sup_{f_n\in H_\om(\R^m)\cap \tilde{L}_{\iy,2\pi n,m}}
   E\left(f_n,\TT_{n,V},L_\iy
   \left(Q_{\pi n}^m\right)\right)\nonumber\\
  && =  \sup_{f\in H_\om(\R^m)\cap L_\iy(\R^m)}
   E\left(f,B_V\cap L_\iy(\R^m),
   L_\iy\left(\R^m\right)\right).
   \ena
   The corresponding values of the sharp constants
   in \eqref{E7.34} were found only for $m=1$
   and for a concave modulus of continuity $\om$
   (see \cite{K1991, D1975}).
   For the H\"{o}lder class $H^\la(\R^m)=H_\om(\R^m)$ with
   $\om(\tau)=\tau^\la,\,\tau\in(0,\iy),\,0<\la\le 1$,
   \eqref{E7.34} immediately implies the relation
   \ba
   &&\lim_{n\to\iy} n^\la
   \sup_{f\in H^\la(\R^m)\cap \tilde{L}_{\iy,2\pi,m}}
   E\left(f,\TT_{nV},L_\iy
   \left(Q_{\pi}^m\right)\right)\nonumber\\
  && =  \sup_{f\in H^\la(\R^m)\cap L_\iy(\R^m)}
   E\left(f,B_V\cap L_\iy(\R^m),
   L_\iy\left(\R^m\right)\right).
   \ea
   A more precise result with
   $H^\la(\R^m)\cap L_\iy(\R^m)$ replaced by $H^\la(\R^m)$
   was obtained in \cite[Theorem 5.2]{G1992}.
   The corresponding univariate result was proved by Bernstein
   \cite{B1946(Nik)}. Bernstein \cite{B1947(const)} attempted
   to extend this result to the class $W^sH^\la(\R^1)$
   but his proof was incomplete (see also Example \ref{Ex7.16}).
   \end{example}

   \begin{example}\label{Ex7.16}
    (Periodic approximation on $W^sH_\om(\R^1)$).
    \end{example}
    The author \cite[Theorem 1]{G2003} solved
    a relatively long-standing problem of finding
    the exact value of
    $\sup_{f\in W^sH_\om(\R^1)}E(f,B_{[-\sa,\sa]},L_\iy(\R^1)),
    \,s\in \Z_+^1$,
    where $\om$ is a concave modulus of continuity.
    The proof is based on a solution of the corresponding
     periodic problem (see Korneichuk \cite{K1971},
     \cite[Theorem 7.2.2]{K1991}) and on a limit relation
     between periodic and nonperiodic approximations
     like \eqref{E7.31} (see \cite[Theorem 2]{G2003}).
     Note that for $s=0$ the problem was solved by
     Dzjadik \cite{D1975}. More references can be found
      in \cite{G2003}.

    \begin{remark}  \label{R7.17}
    It appears that the most challenging and complicated
    part of condition checking is checking Condition C3)
    and its versions C3*), C3.1), C3.1.1), C3.2*), and
    C3.3).
    Despite a number of various examples presented
    in Section \ref{S7}, we could not find a pattern of
    constructing sequences $\left\{f_n\right\}_{n=1}^\iy$
    for an entire function $f$ of exponential type from
    $B_V\cap L_p(\R^m)$. We briefly discuss below
    approximation methods that are used for  checking
    some of those conditions.

    For example, Condition C3.1.1) for approximation $f$
     by multivariate trigonometric polynomials
     is checked by using the Levitan-type polynomials
     (see Example \ref{Ex7.2} and \cite{GT2017, G2018}),
     while a different
     approximation technique
     is used for checking Condition C3.2*)
      (see Example \ref{Ex7.8} and \cite{G1992}).
      Condition C3.1.1) for approximation $f$
     by algebraic polynomials
     is checked by using the Chebyshev expansion for
     the univariate case  and
     the Paley-Wiener type theorem for the multivariate one
     (see Example \ref{Ex7.3}
     and \cite{B1946, G1982}).
     Weighted polynomial approximation to $f$ with
     exponential weights uses Lagrange
     interpolation polynomials and properties of
     harmonic functions (see Example \ref{Ex7.10} and
     \cite{G2008, G2021}).
     \end{remark}
     \noindent
\textbf{Acknowledgements} We are grateful to all anonymous referees
 for valuable suggestions.

\end{document}